\documentclass[12pt,a4paper]{article}
\usepackage{xspace}
\usepackage{authblk}
\def\isfinal{1}
\def\ifnofinal{\ifnum\isfinal=0}

\usepackage{etoolbox}
\patchcmd{\thebibliography}{*}{}{}{}
\pretocmd\thebibliography{\csname c@secnumdepth\endcsname=-2 }{}{}

\usepackage[T1]{fontenc}
\usepackage[french,english]{babel}
\usepackage{bold-extra}

\usepackage[margin=6em]{geometry}

\usepackage{amssymb,amsthm,aliascnt,pdfsync,stmaryrd,mathtools,tabularx}
\usepackage{wasysym}

\usepackage{nicematrix}
\usepackage{chemarrow,tikz-cd,denpastyle}

\usepackage[colorlinks=true,citecolor=orange,linkcolor=blue,bookmarksnumbered=true]{hyperref}
\hypersetup{bookmarksdepth=2}
\usepackage[inline]{enumitem} 
\setlist[enumerate]{label=(\roman*)} 

\usepackage{tocloft}
\usepackage{splitidx}              
\makeindex

\makeatletter
\def\namedlabel#1#2{\begingroup
   \def\@currentlabel{#2}%
   \label{#1}\endgroup
}
\makeatother

\usepackage{url}
\makeatletter
\newcommand{\address}[1]{\gdef\@address{#1}}
\newcommand{\email}[1]{\gdef\@email{#1}}
\newcommand{\@endstuff}{\par\vspace{\baselineskip}\noindent\small
\begin{tabular}{@{}l}\scshape\@address\\\textit{E-mail address:} \@email\end{tabular}}
\AtEndDocument{\@endstuff}
\makeatother

\usepackage[switch]{lineno}

\usepackage{soul} 

\usepackage[backend=biber,doi=false,url=false,isbn=false,url=true]{biblatex}
\DeclareFieldFormat[article]{volume}{\textbf{#1}}
\DeclareFieldFormat[article]{number}{-#1}
\DeclareFieldFormat[book]{volume}{\textbf{#1}}
\renewbibmacro{in:}{%
  \ifentrytype{article}{}{\printtext{\bibstring{in}\intitlepunct}}}

\renewbibmacro*{volume+number+eid}{%
  \printfield{volume}
  \printfield[parens]{number}%
  \setunit{\addcomma\space}%
  \printfield{eid}}
\DeclareFieldFormat{postnote}{#1}
\DeclareFieldFormat{multipostnote}{#1}
\addto\extrasfrancais{

}
\addto\extrasenglish{

}

\newcommand{\ubar}[1]{\underline{#1}} 
 
\newcommand{\im}{\mathrm{im}} 
\newcommand{\can}{\mathrm{can}}
\newcommand{\triv}{\mathrm{triv}} 
\renewcommand{\lg}{\mathrm{lg}}

\DeclareMathOperator{\Tor}{Tor}
\DeclareMathOperator{\Ext}{Ext}
\DeclareMathOperator{\Rhom}{Rhom}
\DeclareMathOperator{\ho}{hom}

\DeclareMathOperator{\End}{End}
\newcommand{\Rlim}{\mathrm{R}\mathop{\varprojlim}}
\newcommand{\Llim}{\mathrm{L}\mathop{\varinjlim}}

\newcommand{\red}{\mathrm{red}}

\newcommand{\fl}{\mathrm{fl}}
\newcommand{\Cl}{\mathrm{Cl}}
\newcommand{\HC}{\mathrm{HC}}

\newcommand{\pt}{\mathrm{pt}}
\newcommand{\perf}{\mathrm{perf}}
\newcommand{\trans}[2]{T_{#1\leftarrow #2}}
\newcommand{\EQ}[2]{\mathrm{Eq}(#1, #2)}
\DeclareMathOperator{\proj}{-proj}
\DeclareMathOperator{\mof}{-mod}
\DeclareMathOperator{\bimod}{-mod-}

\newcommand{\tor}{\mathscr{T}}

\DeclareMathOperator{\Mod}{-Mod}

\newcommand{\Db}{\rmD^{\rmb}}

\DeclareMathOperator{\Spm}{Spm}

\renewcommand{\pre}[1]{\prescript{#1}{}}

\newcommand{\Fl}{\mathrm{Fl}}

\newcommand{\Rfin}{R}
\newcommand{\Raff}{\Phi}

\newcommand{\Scom}{\bfS^{\wedge}}
\newcommand{\Waff}{\ha W}
\newcommand{\Wext}{\til W}

\newcommand{\Deltafin}{\Delta}
\newcommand{\Deltaaff}{\ha\Delta}
\newcommand{\haffd}{\frakh^{1}}
\newcommand{\haff}{\til\frakh}
\newcommand{\hfin}{\frakh}
\newcommand{\shift}[1]{\mathsmaller{\langle#1\rangle}}
\newcommand{\BB}[2]{\ifthenelse{\equal{#1}{}}{\bfB^{\shift{#2}}}{\prescript{}{#1}\bfB_{#2}}}
\newcommand{\Bscr}[3]{\prescript{}{#1}\scrB^{#3}_{#2}}

\newcommand{\BM}[1]{\rmH^{#1}_{\bullet}}

\makeatletter
\def\subsubsection{\@startsection{subsubsection}{3}%
	\z@{.7\baselineskip}{-.5em}%
  {\normalfont\itshape}}
	
\makeatother

\title{Derived equivalences for trigonometric double affine Hecke algebras}
\date{\today}
\author{\sc Wille Liu}

\address{{\sc Institute of Mathematics, Academia Sinica}.  \\
	6F, Astronomy-Mathematics Building, 
	No. 1, Sec. 4, Roosevelt Road, 
	Taipei, Taiwan.
}
\email{\href{mailto:wliu@sinica.edu.tw}{wliu@sinica.edu.tw}}

\begin{document}

\maketitle
\abstract{The trigonometric double affine Hecke algebra $\bfH_c$ for an irreducible root system depends on a family of complex parameters $c$. Given two families of parameters $c$ and $c'$ which differ by integers, we construct the translation functor from $\bfH_{c}\Mod$ to $\bfH_{c'}\Mod$ and prove that it induces equivalence of derived categories. This is a trigonometric counterpart of a theorem of Losev on the derived equivalences for rational Cherednik algebras. }

\section*{Introduction}
\subsection*{Trigonometric double affine Hecke algebras}
Let $(\hfin_{\bbR}, \Rfin)$ be an irreducible root system, where $\hfin_{\bbR}$ is a euclidean vector space and $\Rfin\subset \hfin_{\bbR}^*$ is the set of roots, the latter of which is assumed to span $\hfin_{\bbR}^*$. Let $ W $ be the Weyl group and $T^{\vee}$ the algebraic torus whose character lattice is the coroot lattice of $(\hfin_{\bbR}, \Rfin)$. The \emph{trigonometric double affine Hecke algebra} $\bfH_c$, introduced by Cherednik in his works on Knizhnik--Zamolodchikov equations around 1995, is a flat deformation of the skew tensor product $\bbC  W \ltimes\calO(T^{\vee}\times\hfin)$, where $\hfin = \hfin_{\bbR}\otimes_{\bbR}\bbC$ and $\calO(T^{\vee}\times\hfin)$ is the ring of regular functions on the algebraic variety $T^{\vee}\times\hfin$. It depends on a family of parameters $c\in \frakP$, where $\frakP = \mathop{\mathrm{Map}}( \Rfin /  W , \bbC)$. There is a remarkable polynomial subalgebra $\calO(\hfin)\subset \bfH_c$, generated by \emph{Dunkl operators}.\par

The representation theory of $\bfH_c$ has been studied from different perspectives. Of great importance is the category of \emph{integrable $\bfH_c$-modules}, denoted by $\rmO(\bfH_c)$. Recall that an $\bfH_c$-module is called integrable if it is finitely generated and the Dunkl operators act locally finitely on it. The category $\rmO(\bfH_c)$ is both noetherian and artinian, and its blocks are indexed by the set of orbits in $\hfin$ under the reflection action of the affine Weyl group $\Waff$. For $[\lambda]\in \hfin/\Waff$, let $\rmO_{\lambda}(\bfH_c)\subseteq \rmO(\bfH_c)$ denote the corresponding block. \par

Let $\bfK_t$ be the extended affine Hecke algebra for the dual root system $(\hfin_{\bbR}^*, \Rfin^{\vee})$, with a family of parameters $t:\Rfin /  W \to \bbC^{\times}$, and let $\bfK_t\mof_{\mathrm{fd}}$ be the category of finite-dimensional $\bfK_t$-modules. The blocks of $\bfK_t\mof_{\mathrm{fd}}$ are indexed by the set $T / W$, where $T$ is the dual torus of $T^{\vee}$; for $[\ell]\in T/W$, let $\bfK_t\mof_{\ell}\subseteq \bfK\mof_{\mathrm{fd}}$ denote the corresponding block. A main feature of the category $\rmO(\bfH_c)$ is the trigonometric \emph{Knizhnik--Zamolodchikov functor} (KZ functor), constructed in~\cite{VV04} in terms of monodromy representation of trigonometric KZ equations. It is a quotient functor of abelian categories $\bbV_c:\rmO(\bfH_c)\to \bfK_t\mof_{\mathrm{fd}}$. Moreover, for each $[\lambda]\in \hfin/\Waff$, it restricts to a functor on the block $\bbV^{\lambda}_c: \rmO_{\lambda}(\bfH_c)\to \bfK_t\mof_{\ell}$, where $\ell = \exp(2\pi i \lambda)\in T$ and $t = \exp(2\pi ic)$. \par

The abelian category $\rmO_{\lambda}(\bfH_c)$ has finite global dimension whereas $\bfK_t\mof_{\ell}$ may have infinite global dimension, thus \emph{singular}, in general. By analogy with algebraic geometry, the functor $\bbV^{\lambda}_c$ can be viewed as a desingularisation of $\bfK_t\mof_{\ell}$. If we are given another family of parameters $c'\in \frakP$ such that $c' - c\in \frakP_{\bbZ} = \mathop{\mathrm{Map}}( \Rfin /  W , \bbZ)$, then $c'$ gives rise to the same parameters $t$ for $\bfK_t$ and yields another quotient functor $\bbV^{\lambda}_{c'}:\rmO_{\lambda}(\bfH_{c'})\to \bfK_t\mof_{\ell}$. \par
This leads to the question of the relation between these two ``desingularisations'' $\bbV^{\lambda}_c$ and $\bbV^{\lambda}_{c'}$. More precisely, one expects that there exists an equivalence on the bounded derived categories $\Db(\rmO_{\lambda}(\bfH_{c}))\cong \Db(\rmO_{\lambda}(\bfH_{c'}))$ which intertwines $\bbV^{\lambda}_c$ and $\bbV^{\lambda}_{c'}$. Even better, one can expect that it extends to an equivalence $\Db(\bfH_{c}\Mod)\cong \Db(\bfH_{c'}\Mod)$. The aim of this article is to prove a slight variant of this statement, where the KZ functor $\bbV^{\lambda}_c$ is replaced by an algebraic version $\bfV^{\lambda}_{c}$ introduced in~\cite{liu22} (it is conjectured that $\bbV^{\lambda}_c$ and $\bfV^{\lambda}_{c}$ are isomorphic), see also~\autoref{subsec:KZ}. 

\subsection*{Derived equivalences}
For $c, c'\in \frakP$ such that $c' - c\in \frakP_{\bbZ}$, we will construct an $(\bfH_{c'}, \bfH_c)$-bimodule $\BB{c'}{c}$. Let $\trans{c'}{c} = \BB{c'}{c}\otimes^{\rmL}_{\bfH_c}\relbar:\Db(\bfH_{c}\Mod)\to \Db(\bfH_{c'}\Mod)$ be the derived tensor product. We call $\trans{c'}{c}$ the \emph{translation functor}. The main theorem is the following:
\begin{theointro}[=\autoref{theo:main}+\autoref{prop:PBP}+\autoref{theo:equivO}]
	The translation functor $\trans{c'}{c}$ is an equivalence of categories. Moreover, for each $[\lambda]\in \hfin / \Waff$, it restricts to an equivalence on subcategories $\Db(\rmO_{\lambda}(\bfH_c))\xrightarrow{\sim}\Db(\rmO_{\lambda}(\bfH_{c'}))$ which \emph{conserves the support} on $T^{\vee} / W$ and intertwines the algebraic KZ functors $\bfV^{\lambda}_c$ and $\bfV^{\lambda}_{c'}$. 
\end{theointro}
For the precise meaning of conservation of support on $T^{\vee} / W$, we refer the reader to~\autoref{subsec:equivO}.
This theorem can be viewed as an analogue of a statement about rational Cherednik algebras conjectured by Rouquier~\cite{rouquier08} and proven by Losev~\cite{losev17,losev21}. However, both the construction and the proof of equivalence require new elements for trigonometric double affine Hecke algebras. \par
Moreover, the bimodule $\BB{c'}{c}$ can be explicitly expressed as subspace of a certain nil-Hecke algebra.
\begin{exam*}
	The trigonometric DAHA of type $A_1$ with parameter $c\in \bbC$ is defined by
	\[
		\bfH_c = \bbC\langle x, s_1, s_0\rangle / (s_1^2 - 1, s_0^2 - 1, s_1 x + xs_1 - c, s_0 x - (1 - x)s_0 + c).
	\]
	If we embed $\bfH_c$ into the nil-Hecke algebra
	\[
		\bfH^{\nil} = \bbC\langle x, \vartheta_1, \vartheta_0\rangle / (\vartheta_1^2, \vartheta_0^2, \vartheta_1 x + x\vartheta_1 - 1, \vartheta_0 x - (1 - x)\vartheta_0 + 1)
	\]
	via the injective ring map 
	\[
		\rho_c: \bfH_c\to \bfH^{\nil},\quad \rho_c(x) = x,\; \rho_c(1 - s_1) = (2x - c)\vartheta_1,\;  \rho_c(1 - s_0) = (1 - 2x - c)\vartheta_0,
	\]
	then we have $\BB{c-1}{c} = \left\{1, \vartheta_1, \vartheta_0  \right\}\cdot \rho_c(\bfH_c)\subseteq \bfH^{\nil}$. One verifies easily that $\BB{c-1}{c}$ is an $(\bfH_{c-1},\bfH_c)$-bimodule via $(\rho_{c-1}, \rho_c)$. 
\end{exam*}

\subsection*{Geometric heuristic}
Although our construction of the translation functor is purely algebraic, its idea arises from the realisation of the trigonometric DAHA as cohomological convolution algebra, as in the context of affine Springer theory~\cite{VV10} or Coulomb branches~\cite{BFN2}, and deserves a mention. \par

Suppose that the root system $(\hfin_{\bbR}, \Rfin)$ is simply laced so that $\frakP \cong \bbC$ and $\bfH_c$ has only one parameter $c\in \bbC$. Let $(\haff_{\bbR}, \til\Rfin)$ be the affinisation of $(\hfin_{\bbR}, R)$, so that we have $\haff_{\bbR}^* = \hfin_{\bbR}^* \oplus \bbR \delta$, where $\delta$ is the primitive imaginary root. Let $\calK = \bbC(\!(\varpi)\!)$, let $\calG = G(\calK)$ be the loop group and let $\calB\subset \calG$ be an Iwahori subgroup with pro-unipotent radical $\calU\subset \calB$. Set $\bfN = \frakg$, $\calN = \bfN\otimes \calK$ and $\calN^+ = \Lie \calU$. Following~\cite{vasserot05,BFN2}, one considers the following affine Springer resolution: 
\[
	\pi: \calT = \calG\times^\calB \calN^+\to \calN,\quad [g : x]\mapsto \Ad_g(x).
\]
Consider the homogeneous version of the trigonometric DAHA: $\bfH_{\delta}$, whose underlying vector space is $\calO(\frakP \times \haff)\otimes \bbC\Waff$, where $\delta$ is the ``quantisation parameter''. The algebra $\bfH_c$ can be obtained from $\bfH_{\delta}$ via specialisation of parameters to $\delta \mapsto 1$ and $c\in\frakP$. There a $(\calG\rtimes \bbC^{\times}\times \bbC^{\times})$-action on these spaces making $\pi$ equivariant and $\bfH_{\delta}$ is isomorphic to the convolution algebra $\BM{\calG\rtimes \bbC^{\times}\times \bbC^{\times}}(\calT\times_{\calN}\calT)$, defined in an \textit{ad hoc} way. The coordinate ring of parameters $\calO(\frakP)$ is sent to the ring of coefficients $\rmH^{\bullet}_{\bbC^{\times}}(\pt)$. \par
One may twist this construction: for $d\in \bbZ$, let $\varpi^{d}\calT = \calG\times^{\calB}\varpi^d\calN^+$; then, the isomorphism $\calT\xrightarrow{\sim} \varpi^{d}\calT$ given by multiplication by $\varpi^d$ twists the $(\calG\rtimes \bbC^{\times}\times \bbC^{\times})$-action, which yields an isomorphism $\bfH_{\delta}\xrightarrow{\sim} \BM{\calG\rtimes \bbC^{\times}\times \bbC^{\times}}(\varpi^d\calT\times_{\calN}\varpi^d\calT)$ such that the restriction of this map to $\calO(\frakP)\to \rmH^{\bullet}_{\bbC^{\times}}(\pt)$ is shifted by $d$. \par
One will then define $\BB{}{d} =\BM{\calG\rtimes \bbC^{\times}\times \bbC^{\times}}(\varpi^d \calT\times_{\calN}\calT)$. It can be equipped with the structure of an $\bfH_{\delta}$-bimodule via the convolution product. When the parameters are specialised to $c\in \frakP$ and $\delta$ is set to be $1$, the quotient of $\BB{}{d}$ becomes an $(\bfH_{c+d},\bfH_{c})$-bimodule, denoted by $\BB{c + d}{c}$. One will then show by geometric methods, notably via the equivariant localisation, that the derived tensor product $\BB{c + d}{c}\otimes^{\rmL}_{\bfH_c}\relbar$ yields a derived equivalence $\Db(\bfH_{c}\Mod) \cong \Db(\bfH_{c+d}\Mod)$. \par

When the root system $(\hfin_{\bbR}, \Rfin)$ is not simply laced, in order to allow full generality for the parameters $c\in \frakP$, one needs to replace the adjoint representation $\bfN$. For type $BC$, one can use Kato's exotic nilpotent cone~\cite{kato08}. For type $F_4$ and $G_2$, one can set $G$ to be the corresponding Chevalley group over a field of characteristic $2$ and $3$ respectively and set $\bfN$ to be the semi-simplification of the adjoint representation (as it is told to the author by Kato).\par

\subsection*{Algebraic construction}
Rather than working geometrically in the context of affine Springer theory, we have opted for a purely algebraic approach for three reasons: \begin{enumerate*}\item the algebraic approach is uniform for all types of root systems \item due to the infinite-dimensional nature of the geometry, the geometric approach would require an \textit{ad hoc} equivariant localisation theorem, which is a considerable technical complication \item the imaginary weights of the $\calG$-representation $\calN$ are problamatic; the algebraic construction allows us get rid of the imaginary part. \end{enumerate*}\par

The spherical polynomial representation of $\bfH_{\delta}$ on $\calO(\frakP \times \haff)$ is faithful and yields an embedding of $\bfH_{\delta}$ into the nil-Hecke algebra $\bfH^{\nil}$ for the affine Weyl group $\Waff$. The geometric counterpart of this embedding is the map $\bfz^*: \BM{\calG\rtimes\bbC^{\times}\times\bbC}(\calT\times_{\calN}\calT)\hookrightarrow \BM{\calG\rtimes\bbC^{\times}\times\bbC}(\Fl\times\Fl)$, induced by the zero section $\bfz:\Fl\hookrightarrow \calT$, where $\Fl = \calG / \calB$ is the affine flag manifold; such a map has appeared notably in~\cite[5(iv)]{BFN2}. \par

Based on this embedding, we will construct a category $\bfA$, whose objects are chambers of a certain hyperplane arrangement on the affine space $\haffd_{\bbR} = \{ z\in  \haff_{\bbR}\;;\; \delta(z) = 1 \}$ and whose hom-spaces are certain subspaces of $\bfH^{\nil}$. The composition of morphisms is given by the multiplication in $\bfH^{\nil}$. For each $d\in \frakP_{\bbZ}$, there is an object $\kappa_d\in \bfA$ and an isomorphism $\bfH_{\delta}\cong \End_{\bfA}(\kappa_d)$. The bimodule $\BB{}{d}$ will be defined as the hom-space $\Hom_{\bfA}(\kappa_0, \kappa_d)$, which is naturally an $\bfH_{\delta}$-bimodule via the composition of morphisms. The translation functor is then defined to be $\trans{c'}{c}=\BB{c'}{c}\otimes^{\rmL}_{\bfH_c}\relbar$, where $\BB{c'}{c}$ is the specialisation of $\BB{}{d}$ as explained above. The construction of $\bfA$ is based on an abstraction of the calculation of Euler classes on the Steinberg-type varieties $\varpi^{d'}\calT\times_{\calN}\varpi^{d}\calT$ for various shifts $d, d'\in \frakP_{\bbZ}$. \par

Let us remark that similar constructions of the category $\bfA$ and the bimodule $\BB{c'}{c}$ have already been carried out by Webster in~\cite{webster19koszul,webster19} in the context of quantised Coulomb branches. However, our strategy to prove the equivalence is quite different from his. 
\subsection*{Strategy of the proof}
The essence of our approach to establish the derived equivalences consists of the following steps: 
\begin{enumerate} 
	\item We describe the completion of the algebra $\bfH_c$ as well as the translation bimodule $\BB{c'}{c}$ at every maximal ideal $\lambda$ of the ring $\calO(\haffd)$ (called \emph{spectral completion}) in terms of a category $\scrA^{c, \lambda}$, imitating the equivariant localisation in the aforementioned cohomological approach; the category $\scrA^{c, \lambda}$ is Morita equivalent to a semi-perfect Frobenius algebra over a commutative noetherian complete local ring.
	\item 
		We define in \autoref{sec:equiv} a central hyperplane arrangement on $\frakP$ for each coset of $\frakP_{\bbQ} = \Map(R/W, \bbQ)$ in $\frakP$; we prove in \autoref{prop:intra} that for two families of parameters $c$ and $c'$ which differ by integers and lie in the same facet, the categories $\bfH_c\Mod$ and $\bfH_{c'}\Mod$ are canonically equivalent.
	\item
		We prove in \autoref{prop:antipode} that for two families of parameters $c$ and $c'$ which differ by integers and lie in antipodal facets, the translation bimodule $\BB{c'}{c}$ is a tilting module and induces a derived equivalence between $\Db(\bfH_c\Mod)$ and $\Db(\bfH_{c'}\Mod)$; this makes use of the spectral completion of $\BB{c'}{c}$ and a duality (\autoref{prop:duality}) induced from the Frobenius structure on $\scrA^{c, \lambda}$.
	\item 
		In \autoref{prop:simple}, we apply an analogue of the degeneration technique of Harish-Chandra bimodules due to Bezrukavnikov--Losev~\cite{BL21} to show that when $c$ and $c'$ differ by integers and are situated in adjacent open chambers, the translation bimodule yields a derived equivalence.
	\item 
		In \autoref{subsec:translation-semigroup}, we study the algebra $\bigoplus_{n\in \bbN} \BB{}{nd}$ of sum of translation bimodules along the semigroup generated by $d\in \frakP_{\bbZ}$; a refinement of the degeneration technique allows us to prove in \autoref{prop:general} the derived equivalence when both $c$ and $c'$ lie in open chambers and differ by integers.
	\item 
		We treat in \autoref{prop:except} the case where $c$ and $c'$ lie in non-open chambers, again by the degeneration technique; this step relies crucially on the use of the algebraic KZ functor.
\end{enumerate}  

\subsection*{Organisation}
The article is organised as follows: \par
In~\autoref{sec:DAHA}, we recollect some previously known results about trigonometric double affine Hecke algebras. Only the results of \autoref{subsec:derived} are original. \par
In~\autoref{sec:HC}, we develop the basic theory of Harish-Chandra bimodules for trigonometric DAHAs. We prove that they form a Serre subcategory closed under $\Tor$ and $\Ext$. Moreover, we prove that the $\Tor$ and $\Ext$ of Harish-Chandra bimodules commute with spectral completion. \par
In~\autoref{sec:A}, we define and study the category $\bfA$ and an embedding of the trigonometric DAHA in it.  \par
In~\autoref{sec:Aloc}, we define and study the category $\scrA^{c, \lambda}$, which is a spectral-completed version of $\bfA$. The main result of this section is the comparison of $\bfA$ and $\scrA^{c, \lambda}$. \par
In~\autoref{sec:transl}, we construct the translation bimodule $\BB{}{d}$ as certain hom-space in $\bfA$ and prove that they are Harish-Chandra bimodules. We describe the spectral completion of $\BB{}{d}$ in terms of $\scrA^{c, \lambda}$. We then prove the noetherianity of the algebra $\bigoplus_{n\ge 0}\BB{}{nd}$.  \par
In~\autoref{sec:equiv}, we define the translation functor $\trans{c'}{c}$ and prove that it is a derived equivalence. We also prove that it preserves integrable modules and their support.

\section{Reminder on the trigonometric DAHA}\label{sec:DAHA}
In this section, we recollect results about the representations of trigonometric DAHA which will be needed later. Except \autoref{subsec:derived} which is original, all other results presented here can be found in certain variant forms in~\cite{VV04, liu22} for example. 

\subsection{Affine root systems}\label{subsec:Raff}
We fix the notation for affine root systems, which will be constantly used throughout this article. \par
Let $(\hfin_{\bbR}, \Rfin)$ be an irreducible root system, either reduced or not, where $\hfin_{\bbR}$ is a euclidean vector space and $\Rfin\subset \hfin_{\bbR}^*$ is the set of roots which spans $\hfin_{\bbR}^*$. The affinisation $(\haff_{\bbR}, \til\Rfin)$ is an affine root system on the $\bbR$-vector space $\haff_{\bbR} = \hfin_{\bbR} \oplus \bbR\partial$, whose dual can be written as $\haff_{\bbR}^* = \hfin_{\bbR}^* \oplus \bbR \delta$ with the pairing
\[
	\langle \alpha + t\delta, \lambda + r\partial\rangle = \langle \alpha, \lambda \rangle + rt,\quad \text{for $\alpha\in \hfin_{\bbR}^*,\;\lambda\in \hfin_{\bbR},\quad $ and $r, t\in \bbR$}.
\]
The set of affine roots is defined to be
\[
	\til\Rfin = \left\{ \alpha + n\delta\in \haff_{\bbR}^*\;;\; \alpha\in \Rfin_{\red},\; n\in \bbZ \right\}\cup \left\{ \alpha + (2n + 1)\delta\in \haff_{\bbR}^*\;;\; \alpha\in \Rfin \setminus \Rfin_{\red},\; n\in \bbZ \right\},
\]
where $\Rfin_{\red}\subseteq \Rfin$ is the subset of indivisible roots. However, for technical reasons, the following modification of $\til\Rfin$ will be more convenient for us:
\[
	\Raff = \left\{ \alpha + n\delta\in \haff_{\bbR}^*\;;\; \alpha\in \Rfin_{\red},\; n\in \bbZ \right\}\cup \left\{ \alpha/2+(n + 1/2)\delta\in \haff_{\bbR}^*\;;\; \alpha\in \Rfin \setminus \Rfin_{\red},\; n\in \bbZ \right\}.
\]
We fix a base (= set of simple roots) $\Deltafin = \left\{ \alpha_1, \ldots, \alpha_n \right\}\subset \Rfin$, which extends to $\Deltaaff = \Deltafin \cup\left\{ \alpha_0 \right\}$, where $\alpha_0\in \Phi$ is the affine root defined as follows: let $\theta\in \Rfin$ be the largest root with respect to $\Deltafin$; then $\alpha_0 = \delta - \theta$ if $\theta\in \Rfin_{\red}$ and $\alpha_0 = (\delta - \theta)/2$ if $\theta\in \Rfin \setminus\Rfin_{\red}$. \par

The affine Weyl group of $(\haff_{\bbR}, \Raff)$ is the Coxeter group $\Waff$ generated by the simple reflections $\left\{ s_\alpha \right\}_{\alpha\in \Deltaaff}$. Let $(m_{\alpha,\beta})_{\alpha,\beta\in \Deltaaff}$ denote the Coxeter matrix, so that $m_{\alpha,\beta} = \ord(s_{\alpha}s_{\beta})$. \par

\begin{exam}[type $BC_n$] \label{exam:BC} Let $n\in \bbZ_{\ge 1}$.  Let $\hfin_{\bbR} = \bbR^n = \bigoplus_{i = 1}^{n}\bbR\epsilon_i$ be the euclidean space with inner product $\langle\epsilon_i, \epsilon_j\rangle = \delta_{i,j}$. The root system of type $BC_n$ is the datum $(\hfin_{\bbR}, \Rfin)$, where $\Rfin =  \Rfin_{\natural}\cup \Rfin_{\sharp}\cup \Rfin_{\flat}$ with
\[
	\Rfin_{\sharp} = \left\{ \pm\epsilon_i\;;\; 1\le i \le n \right\},\; \Rfin_{\natural} = \left\{ \pm\epsilon_i\pm \epsilon_j\;;\; 1\le i < j \le n \right\},\; \Rfin_{\flat} = \left\{ \pm 2\epsilon_i\;;\; 1\le i \le n \right\}.
\]
The indivisible roots are $\Rfin_{\red} = \Rfin_{\natural}\cup \Rfin_{\sharp}$. The three subsets $\Rfin_{\sharp}$, $\Rfin_{\natural}$ and $\Rfin_{\flat}$ (when $n=1$, the two sets $\Rfin_{\sharp}, \Rfin_{\flat}$) are the $ W $-orbits in $\Rfin$ and they are characterised by the length of roots. We shall write $\Rfin /  W  = \left\{ \natural, \sharp, \flat \right\}$ ($\Rfin /  W  = \left\{ \sharp,\flat \right\}$ when $n = 1$). 
The (standard) base given by $\Deltafin = \left\{ \alpha_i\right\}_{1\le i \le n}$ with $\alpha_i = \epsilon_i - \epsilon_{i+1}$ for $i = 1, \ldots, n-1$ and $\alpha_n = \epsilon_n$. The largest root is $\theta = 2\epsilon_1$. \par

The affine root system of type $\widetilde{BC_n}$ is the datum $(\haff_{\bbR}, \Raff)$, where $\Raff = \Raff_{\natural}\cup\Raff_{\sharp}\cup\Raff_{\flat}$ and $\Raff_{\natural} = \Rfin_{\natural} \times \bbZ \delta$, $\Raff_{\sharp} = \Rfin_{\sharp}\times \bbZ \delta$, $\Raff_{\flat} = (1/2)\Rfin_{\flat}\times (1/2 + \bbZ) \delta$.
Let $\alpha_0 = (\delta - \theta)/2$. The standard affine base is $\Deltaaff = \Deltafin \cup \left\{ \alpha_0 \right\}$.  
 Note that $\alpha_n\in \Raff_{\sharp}$, $\alpha_0\in \Raff_{\flat}$ and $\alpha_i \in \Raff_{\natural}$ for $i = 1, \ldots, n-1$.
\end{exam}
\begin{exam} [type $F_4$]\label{exam:F4}
	Let $\hfin_{\bbR} = \bigoplus_{i = 1}^4 \bbR\epsilon_i$ be the euclidean space. The root system of type $F_4$ is the datum $(\hfin_{\bbR}, \Rfin)$, where $\Rfin =  \Rfin_{\natural}\cup \Rfin_{\sharp}$ with
	\[
		\Rfin_{\natural} = \left\{ \pm\epsilon_i\pm \epsilon_j\;;\; 1\le i < j \le 4 \right\},\;\Rfin_{\sharp} = \left\{ \pm \epsilon_i\;;\; 1 \le i \le 4\right\}\cup \left\{(\pm \epsilon_1 \pm \epsilon_2\pm \epsilon_3\pm \epsilon_4)/2 \right\}. 
	\]
	The standard base is given by $\left\{ \alpha_i \right\}_{1\le i \le 4}$ with $\alpha_1 = \epsilon_2 - \epsilon_3$, $\alpha_2 = \epsilon_3 - \epsilon_4$, $\alpha_3 = \epsilon_4$ and $\alpha_4 = (\epsilon_1 - \epsilon_2 - \epsilon_3 - \epsilon_4)/2$. The largest root is $\theta = \epsilon_1 + \epsilon_2$. The affine root system of type $\widetilde{F_4}$ is the datum $(\haff_{\bbR}, \Raff)$, where $\Raff = \Raff_{\natural}\cup\Raff_{\sharp}$ and $\Raff_{\natural} = \Rfin_{\natural} \times \bbZ \delta$, $\Raff_{\sharp} = \Rfin_{\sharp}\times \bbZ \delta$.
\end{exam}

\begin{exam} [type $G_2$]\label{exam:G2}
	Let $\hfin_{\bbR} = \mathop{\mathrm{span}}_{\bbR}(\epsilon_1 - \epsilon_2, \epsilon_2 - \epsilon_3)\subset \bigoplus_{i = 1}^3 \bbR\epsilon_i$. The root system of type $G_2$ is the datum $(\hfin_{\bbR}, \Rfin)$, where $\Rfin =  \Rfin_{\natural}\cup \Rfin_{\sharp}$ with
	\[
		\Rfin_{\natural} = \left\{\pm(2\epsilon_i - \epsilon_j - \epsilon_k)\;;\; \left\{ i, j, k \right\} = \left\{ 1, 2, 3 \right\}\right\},\; \Rfin_{\sharp} = \left\{ \pm(\epsilon_i- \epsilon_j)\;;\; 1\le i < j \le 3 \right\}. 
	\]
	The standard base is given by $\left\{ \alpha_i \right\}_{1\le i \le 2}$ with $\alpha_1 = \epsilon_1 - \epsilon_2$ and $\alpha_2 = -2\epsilon_1+ \epsilon_2 + \epsilon_3$. The largest root is $\theta = -\epsilon_1 - \epsilon_2 + 2\epsilon_3$. The affine root system of type $\widetilde{G_2}$ is the datum $(\haff_{\bbR}, \Raff)$, where $\Raff = \Raff_{\natural}\cup\Raff_{\sharp}$ and $\Raff_{\natural} = \Rfin_{\natural} \times \bbZ \delta$, $\Raff_{\sharp} = \Rfin_{\sharp}\times \bbZ \delta$.

\end{exam}

\subsection{Trigonometric DAHA}\label{subsec:daha-def}
\subsubsection{}
Fix an irreducible root system $(\hfin_{\bbR}, \Rfin)$ as above. We put $\hfin = \hfin_{\bbR} \otimes_{\bbR} \bbC$ and $\haff = \haff_{\bbR} \otimes_{\bbR} \bbC$.
Let $\bbC[\bfc_{*}; *\in \Rfin /  W ]$ be the polynomial ring whose variables are indexed by the $W$-conjugacy classes of roots. We will write $\frakP = \Spm \bbC[\bfc_{*}; *\in \Rfin /  W ]$ for its maximal spectrum, so that $\bbC[\bfc_{*}; *\in \Rfin /  W ] = \calO(\frakP)$. Let $\Wext =  W \ltimes P^{\vee}_R$ be the extended affine Weyl group, where $P^{\vee}_R\subseteq \hfin_{\bbR}$ is the coweight lattice. There is a bijection $\Rfin /  W \xrightarrow{\sim} \Raff/\Wext$ which sends $R_*$ to $\Phi_*$ for $*\in \left\{ \natural, \sharp,\flat \right\}$ as in~\autoref{exam:BC}--\autoref{exam:G2}. For $*\in \Rfin /  W $ and $\alpha\in \Raff_{\ast}$, we will write $\bfc_{\alpha} = \bfc_{\ast}$.  \par

\subsubsection{}
The (homogeneous) trigonometric double affine Hecke algebra for $(\haff_{\bbR}, \Rfin, \Deltaaff)$, denoted by $\bfH_{\delta}$, is the unital associative algebra over $\calO(\frakP)$ generated by the sets $\left\{ s_{\alpha} \right\}_{\alpha\in \Deltaaff}$ and $\left\{ x^{\mu} \right\}_{\mu\in \haff^*}$ subject to the following relations for $\mu,\nu\in \haff^*$, $a, b\in \bbC$ and $\alpha,\beta\in \Deltaaff$:
\begin{align*}
	x^{a\mu+ b\nu} = ax^\mu + bx^\nu,\quad \left[ x^{\mu},x^{\nu} \right] = 0,\quad  (s_{\alpha}s_{\beta})^{m_{\alpha,\beta}} = 1,\quad s_{\alpha}x^{\mu} - x^{s_{\alpha}(\mu)}s_{\alpha} = \bfc_{\alpha}\langle \mu, \alpha^{\vee}\rangle ,
\end{align*}
where $\alpha^{\vee} = \alpha / \langle \alpha, \alpha\rangle$. 

The family $\left\{ x^{\mu} \right\}_{\mu\in \haff^*}$ generates a polynomial subalgebra $\Sym(\haff^*) \cong \calO(\haff)$ while $\left\{ s_{\alpha} \right\}_{\alpha\in \Deltaaff}$ generates a subalgebra isomorphic to the group ring $\bbC \Waff$ of the affine Weyl group. We will simply write $\mu = x^{\mu}$ for $\mu\in \haff^*$. Moreover, the centre of $\bfH_{\delta}$ is the polynomial ring $\calO(\frakP)[\delta]$. There is a decomposition of $\bfH_{\delta}$ as vector space:
\[
	\bfH_{\delta} = \calO(\frakP\times \haff)\otimes \bbC\Waff.
\]
We set $\bfS_{\delta} = \calO(\frakP\times \haff)\subset \bfH_{\delta}$; it is a polynomial subalgebra.

\subsubsection{}\label{subsubsec:filtcan}
We put an $\bbN$-grading on $\bfH_{\delta}$ by setting $\deg \frakP^* = \deg \haff^* = 1$ and $\deg \bbC \Waff = 0$. Let $\bfH = \bfH_{\delta} / (\delta - 1)$ (resp. $\bfS = \bfS_{\delta} / (\delta - 1)$) be the quotient ring of $\bfH_{\delta}$ (resp. $\bfS_{\delta}$) by the two-sided ideal generated by $\delta - 1$. Let $\haffd$ be the vanishing locus of $\delta - 1$ in $\haff$. We may identify $\bfS$ with the coordinate ring $\calO(\frakP\times \haffd)$. The ring $\bfH$ is equipped with a filtration $F^{\can}_{\bullet}\bfH$ induced from the grading on $\bfH_{\delta}$, called the \emph{canonical filtration} of $\bfH$. Since $\bfH_{\delta}$ is graded-free over $\bbC[\delta]$, the homogeneous trigonometric DAHA $\bfH_{\delta}$ is isomorphic to the Rees algebra of the canonical filtration $F^{\can}_{\bullet}\bfH$ in a natural way. Similarly, the graded subring $\bfS_{\delta}\subseteq \bfH_{\delta}$ can be identified with the Rees algebra of the canonical filtration $F^{\can}_{\bullet}\bfS = \bfS\cap F^{\can}_{\bullet}\bfH$. \par

Given a family of parameters $c\in \frakP$, we will write $c_\ast = \bfc_\ast(c)$ for $\ast\in R/W$ and define the specialisation of parameters $\bfH_c = \bfH / (\bfc - c)$ and $\bfS_c = \bfS / (\bfc - c)$.
\subsubsection{}
The following is well known and can be established by passing to the associated graded ring $\gr^{\can}\bfH$, see~\cite[\S 2.2]{liu22} for example:
\begin{prop}\label{prop:dimglob}
	The rings $\bfH$ and $\bfH_c$ for $c\in \frakP$ have finite left and right global dimension. \hfill\qedsymbol
\end{prop}

\subsection{Integrable modules}\label{subsec:OH}
\begin{defi}
	A finitely generated $\bfH$-module $M$ is called \emph{integrable} if $M$ is, by restriction, a locally finite $\bfS$-module. 
\end{defi}
We let $\rmO(\bfH)\subset \bfH\mof$ denote the full subcategory of integrable $\bfH$-modules. 
\subsubsection{}
For each point $(c,\lambda)\in \frakP\times \haffd$, we let $\frakm_{c,\lambda}\subset \bfS$ denote the maximal ideal attached to $(c, \lambda)$.
If $M$ is an integrable $\bfH$-module, then there is a decomposition
\[
	M = \bigoplus_{(c,\lambda)\in \frakP\times \haffd} M_{c,\lambda},\quad M_{c,\lambda} = \bigcup_{k\ge 0} \left\{ m\in M\;;\; \frakm_{c,\lambda}^k m = 0 \right\}.
\]
The category $\rmO(\bfH)$ can be decomposed into blocks
\[
	\rmO(\bfH) = \bigoplus_{(c,[\lambda])\in \frakP\times\haffd / \Waff} \rmO_{c,\lambda}(\bfH),\quad \mathrm{Obj}(\rmO_{c,\lambda}(\bfH)) = \left\{ M\in \rmO(\bfH)\;;\; M = \bigoplus_{\lambda'\in [\lambda]} M_{c,\lambda'} \right\}.
\]
\subsubsection{}
Given $c\in \frakP$, we define $\rmO(\bfH_c)$ to be the full subcategory of $\bfH_c\mof$ consisting of modules $M$ which, regarded as $\bfH$-module, is in $\rmO(\bfH)$. The block decomposition
\[
	\rmO(\bfH_c) = \bigoplus_{[\lambda]\in \haffd / \Waff} \rmO_{\lambda}(\bfH_c)
\]
works in a similar way.

\subsection{Spectral completion of \texorpdfstring{$\bfH$}{H}}\label{subsec:completion}
\subsubsection{}
Define the following (non-unital) ring
\[
	\scrS = \bigoplus_{(c, \lambda)\in \frakP\times \haffd} \Scom_{c, \lambda},
\]
where $\Scom_{c, \lambda}$ is the completion of $\bfS$ at the defining ideal $\frakm_{c, \lambda}\subset \bfS$ of $(c, \lambda)\in \frakP \times \haffd$. For each $(c,\lambda)\in \frakP\times \haffd$, the unit element of $\Scom_{c, \lambda}$ is an idempotent in $\scrS$, denoted by $1_{c,\lambda}$. Then, $\left\{ 1_{c, \lambda} \right\}_{(c,\lambda)\in  \frakP\times\haffd}$ is an orthogonal family of central idempotents which is complete in the sense that $\scrS= \bigoplus_{(c, \lambda)\in \frakP\times\haffd}1_{c, \lambda}\scrS$.
\begin{lemm}\label{lemm:platitude}
	The ring $\scrS$ is faithfully flat over $\bfS$.
\end{lemm}
\begin{proof}
	The flatness is due to the noetherianity of $\bfS$, see~\cite[\S 3.4, Th\'eor\`eme 3.(iii)]{ACIII}. For the faithfulness, let $M$ be any $\bfS$-module and let $a\in M\setminus\left\{ 0 \right\}$. Choose a maximal ideal $\frakm\subset \bfS$ containing the annihilator $\ann_{\bfS}(a)$, so that the image of $a$ in the localisation $\bfS_{\frakm}\otimes_{\bfS}M$ is non-zero. Since the completion $\Scom_{\frakm}$ is faithfully flat over $\bfS_{\frakm}$ (the latter being noetherian and local), the image of $a$ in $\Scom_{\frakm}\otimes_{\bfS}M$ is non-zero. It follows that $\Scom_{\frakm}\otimes_{\bfS}M\neq 0$ and hence $\scrS\otimes_{\bfS}M\neq 0$ holds, so $\scrS$ is faithfully flat over $\bfS$.
\end{proof}


\subsubsection{}\label{subsubsec:scrH}
For each element $w\in \Waff$, the action $w: \haffd\xrightarrow{\sim} \haffd$ induces an isomorphism of complete local rings for each pair $(c, \lambda)\in \frakP\times \haffd$:
\[
	\pre{w}(\cdot):\Scom_{ c, \lambda}\to \Scom_{c, w\lambda},\quad f_{\lambda}\mapsto \pre{w}(f_{\lambda})\quad \text{for } f_{\lambda}\in \Scom_{c, \lambda}.
\]
For each affine root $\alpha\in \Raff$, define the Demazure operator $\vartheta_{\alpha}:\scrS\to \scrS$ to be the linear map satisfying
\[
	\vartheta_{\alpha}(f_{\lambda}) = \alpha^{-1}(f_{\lambda} - \pre{s_{\alpha}}(f_{\lambda})) \quad  \text{for $(c, \lambda)\in \frakP\times \haffd$ and $f_{\lambda}\in \Scom_{c, \lambda}$}.
\]
Note that the above expression makes sense: when $\alpha(\lambda) = 0$, we have $\alpha^{-1}(f_{\lambda} - \pre{s_{\alpha}}(f_{\lambda}))\in \Scom_{c, \lambda}$; when $\alpha(\lambda)\neq 0$, the element $\alpha$ is invertible in both $\Scom_{c, \lambda}$ and $\Scom_{c, s_\alpha\lambda}$, so we have $\alpha^{-1}f_{\lambda} - \alpha^{-1}\pre{s_{\alpha}}(f_{\lambda})\in \Scom_{c, \lambda}\oplus \Scom_{c, s_{\alpha}\lambda}$. \par

Let $\scrH = \scrS\otimes_{\bfS}\bfH$. We extend the ring structures on $\scrS$ and $\bfH$ to $\scrH$ by the following rules:
\begin{align*}
	(f\otimes 1)(g\otimes a) &= fg\otimes a \\
	(f\otimes s_\alpha)(g\otimes a) &= f\cdot \pre{s_\alpha}g \otimes s_\alpha a + \bfc_{\alpha} \vartheta_{\alpha}(g)\otimes a
\end{align*}
for $f,g\in \scrS$, $a\in \bfH$ and $\alpha\in \Deltaaff$. This ring structure induces an isomorphism:
\begin{equation}\label{equa:HSSH}
	\bfH\otimes_{\bfS}\scrS \cong \scrS\otimes_{\bfS}\bfH =\scrH,\quad a\otimes f\mapsto \sum_{(c, \lambda)\in \frakP\times \haffd}(1_{c,\lambda}\otimes a)(f\otimes 1).
\end{equation}
Note that the above summation must be finite.  There is a natural $\bfH$-bimodule structure on $\scrH$. \par

Each element of $\scrH$ can be expressed as finite sum of the following form in a unique way:
\[
	\sum_{(c,\lambda)\in \frakP\times\haffd}\sum_{w\in \Waff} 1_{c,w\lambda} w f_{c,\lambda,w}1_{c, \lambda}\quad\text{for } f_{c,\lambda,w}\in \Scom_{c, \lambda}.
\]

\subsubsection{}
The ring $\scrH$ admits a decomposition into direct sum (coproduct in the category of non-unital rings):
\begin{equation}\label{equa:decompH}
	\scrH = \bigoplus_{(c, [\lambda])\in \frakP\times \haffd/\Waff} \scrH^{c, \lambda},\quad \scrH^{c, \lambda} = \bigoplus_{\lambda',\lambda''\in [\lambda]}1_{c,\lambda'}\scrH1_{c,\lambda''}.
\end{equation}
Note that this decomposition is in accordance with that of the integrable modules given in \autoref{subsec:OH}.

\subsubsection{}\label{subsubsec:Hcl}
Given $c\in \frakP$, we define $\scrH_c = \scrH / (\bfc - c)$ to be the specialisation of parameters at $c$. Similarly to \eqref{equa:decompH}, there is a decomposition $\scrH_c = \bigoplus_{[\lambda]\in \haffd/\Waff} \scrH^{\lambda}_c$. 

\subsection{Properties of \texorpdfstring{$\scrH$}{H}-modules}
\subsubsection{Module categories.}
Let $\scrH\Mod$ denote the category of non-degenerate $\scrH$-modules (i.e. $\scrH$-modules $\scrM$ satisfying $\scrM = \bigoplus_{(c,\lambda)\in \frakP\times \haffd} 1_{c, \lambda} \scrM$). It admits a set of compact projective generators $\{\scrH 1_{c,\lambda}\}_{(c,\lambda)\in \frakP\times \haffd}$. Let $\scrH\mof\subset \scrH\Mod$ denote the full subcategory of compact objects and let $\scrH\mof_{\fl}\subseteq \scrH\mof$ be the full subcategory of objects of finite length.
The categories $\scrH_c\Mod$, $\scrH_c\mof$ and $\scrH_c\mof_{\fl}$ are defined similarly for each $c\in \frakP$. \par
The decomposition of ring \eqref{equa:decompH}, induces a decomposition of the module category:
\begin{equation}\label{equa:decompMod}
	\scrH\Mod\cong \bigoplus_{(c, [\lambda])\in \frakP \times\haffd / \Waff}\scrH^{c, \lambda}\Mod.
\end{equation}
A similar decomposition holds for $\scrH_c\Mod$.
\subsubsection{Centre.}\label{subsubsec:Z} Let $(c, [\lambda])\in \frakP\times\haffd$. We put $\scrZ^{c,\lambda} = (\Scom_{c, \lambda})^{\Waff_{\lambda}}$. It is a complete regular local ring since the stabiliser $\Waff_{\lambda} = \Stab_{\Waff}(\lambda)$ is a finite reflection groups (see~\cite[V, \S 5.3, Th 3]{bourbakiLie456}). For each $\lambda'\in [\lambda]$, choose $w\in \Waff$ such that $\lambda' = w\lambda$ and consider the embedding
\[
	\scrZ^{c,\lambda}\hookrightarrow \Scom_{c, \lambda}\xrightarrow{w}\Scom_{c, \lambda'},\quad z\mapsto z_{c,\lambda'}:= \pre{w}z,
\]
which is independent of the choice of $w$ and makes $\Scom_{c, \lambda'}$ a free $\scrZ^{c,\lambda}$-module of rank $\# \Waff_{\lambda}$. 
It defines a $\scrZ^{c,\lambda}$-action on the identity functor of $\scrH^{c,\lambda}\Mod$ by natural transformations:
\[
	\chi:\scrZ^{c,\lambda}\to \End(\id_{\scrH^{c,\lambda}\Mod}),\quad \chi(z)(m) = \sum_{\lambda'\in [\lambda]} z_{c,\lambda'} m.
\]
Consequently, the category $\scrH^{c,\lambda}\Mod$ comes with a natural $\scrZ^{c,\lambda}$-linear structure (i.e. it is enriched over $\scrZ^{c,\lambda}\Mod$); moreover, the Hom-spaces between objects from $\scrH^{c,\lambda}\mof$ (resp. $\scrH^{c,\lambda}\proj$) are finitely generated $\scrZ^{c,\lambda}$-modules (resp. free $\scrZ^{c,\lambda}$-modules of finite rank). \par
The same construction can be applied to $\scrH_c$, with $\scrZ^{c, \lambda}$ replaced by $\scrZ^{\lambda}_c = \scrZ^{c, \lambda} / (\bfc - c)$. 
\begin{lemm}\label{lemm:FracZ}
	Let $\tor\subset \scrH^{c,\lambda}\Mod$ be the full subcategory consisting of objects which are locally finite-dimensional as $\scrZ^{c,\lambda}$-modules. Then, the quotient category $\scrH^{c,\lambda}\Mod/\tor$ is equivalent to $(\Frac\scrZ^{c,\lambda})\Mod$. A similar statement holds for $\scrH_c^{\lambda}$ for $c\in \frakP$. 
\end{lemm}
\begin{proof}
	Use the compact projective generators $\left\{ \scrH^{c,\lambda} 1_{c, \lambda'}\right\}_{\lambda'\in [\lambda]}$ and observe that for each $\lambda'\in [\lambda]$, the morphism 
		\[
			\relbar\cdot s_{\alpha}1_{c, s_{\alpha}\lambda'}:  \scrH^{c,\lambda} 1_{c, \lambda'}\to \scrH^{c,\lambda} 1_{c, s_{\alpha}\lambda'}
		\]
		is injective with cokernel lying in $\tor$; it follows that this map becomes an isomorphism after passing to the quotient category $\scrH^{c,\lambda}\Mod/\tor$; therefore, the image of $\scrH^{c,\lambda}1_{c,\lambda}$ in $\scrH^{c,\lambda}\Mod/\tor$ is a compact projective generator. Let $A = \End_{\scrH^{c,\lambda}}(\scrH^{c,\lambda} 1_{c, \lambda})\otimes_{\scrZ^{c,\lambda}}\Frac\scrZ^{c,\lambda}$, so that $\scrH^{\lambda}_c\Mod/\tor$ is equivalent to $A\Mod$ and
		\[
			A \cong 1_{c, \lambda}\scrH^{c,\lambda} 1_{c, \lambda}\otimes_{\scrZ^{c,\lambda}}\Frac\scrZ^{c,\lambda}
		\]
		is isomorphic to the subalgebra of $\End_{\bbC}(\Frac\Scom_{c, \lambda})$ generated by $\Frac\Scom_{c, \lambda}$ and the action of $\Waff_{\lambda}$. Since $\Waff_{\lambda}$ is a reflection group on $\haffd_{\bbR}$, the Galois theory implies that this subalgebra is equal to $\End_{\Frac\scrZ^{c,\lambda}}(\Frac \Scom_{c, \lambda})$, which is a matrix algebra of rank $\# \Waff_{\lambda}$ over the field $\Frac\scrZ^{c,\lambda}$. The statement follows.
\end{proof}

\subsubsection{Coxeter complex.}\label{sssec:coxeter}
For a non-constant affine function $\alpha$ on $\haffd_{\bbR}$, let $H_{\alpha} = \left\{ h\in \haffd_{\bbR}\;;\; \alpha(h) = 0\right\}$ be its zero locus in $\haffd_{\bbR}$. The affine hyperplane arrangement $\{H_{\alpha}\}_{\alpha\in \Raff}$ yields a simplicial decomposition of $\haffd_{\bbR}$, which can be identified with the Coxeter complex of $\Waff$. The relative interior of simplices is called \emph{facets}. The facets of maximal dimension are called \emph{alcoves}. The fundamental alcove associated with the base $\Deltaaff$ is defined by $\nu_{0} = \bigcap_{\alpha\in \Deltaaff} \alpha^{-1}(\bbR_{>0})$, where the affine root $\alpha$ is viewed as affine function on $\haffd_{\bbR}$.
\subsubsection{Clan decomposition.}
Fix $(c, \lambda)\in \frakP\times \haffd$. Consider the subfamily of affine roots $\Raff_{c, \lambda} = \left\{ \alpha\in \Raff\;;\;\alpha(\lambda) = c_{\alpha}\right\}$ and the hyperplane arrangement $\left\{ H_{\alpha} \right\}_{\alpha\in \Raff_{c, \lambda}}$ on $\haffd_{\bbR}$. 
\begin{defi}
	The connected components of the complements $\haffd_{\bbR} \setminus\bigcup_{\alpha\in \Raff_{c, \lambda}} H_{\alpha}$ are called $(c, \lambda)$-clans. 
\end{defi}
Let $\Cl^{c,\lambda}(\haffd_{\bbR})$ denote the set of $(c, \lambda)$-clans in $\haffd_{\bbR}$. Note that $\Phi_{c, \lambda}$ and $\Cl^{c, \lambda}(\haffd_{\bbR})$ are finite sets.
Given a $(c, \lambda)$-clan $\frakC\subset \haffd_{\bbR}$, we define its \emph{salient cone} to be $\kappa = \left\{ h\in \hfin_{\bbR}\;;\; \frakC + h\subseteq \frakC \right\}$. Its dual cone can be described as follows:
\[
	\kappa^{\vee} = \sum_{\substack{\alpha\in \Raff_{c, \lambda} \\ \alpha(\frakC)\subseteq \bbR_{> 0}}}\bbR_{\ge 0} \ba\alpha + \sum_{\substack{\alpha\in \Raff_{c, \lambda} \\ \alpha(\frakC)\subseteq \bbR_{< 0}}}\bbR_{\le 0} \ba\alpha,\quad \text{where }\ba\alpha = \alpha\mid_{\frakh_{\bbR}}\in \hfin^*_{\bbR}.
\]
In other words, $\kappa$ tells in which directions the $(c,\lambda)$-clan $\frakC$ is unbounded. Each alcove of the Coxeter complex for $\Raff$ is contained in a unique $(c,\lambda)$-clan. \par
\begin{defi}\label{defi:generic}
	A $(c, \lambda)$-clans $\frakC\in \Cl^{c,\lambda}(\haffd_{\bbR})$ is said to be \emph{generic} if its salient cone $\kappa$ satisfies $\dim \kappa = \dim \hfin_{\bbR}$.
\end{defi}
In other words, a $(c, \lambda)$-clan $\frakC\in \Cl^{c,\lambda}(\haffd_{\bbR})$ is generic if $\ba\alpha$ does not vanish on $\kappa$ for any $\alpha\in \Raff_{c, \lambda}$. 

\subsubsection{Morita equivalence with complete noetherian algebra.}\label{subsubsec:morita}
Fix $(c, \lambda)\in \frakP\times \haff$ as in the previous paragraph. The action of $\Waff_{\lambda} = \Stab_{\Waff}(\lambda)$ on $\haffd_{\bbR}$ preserves $\Raff_{c, \lambda}$ and thus induces an action on $\Cl^{c,\lambda}(\haffd_{\bbR})$. Let $\Sigma\subset \Waff$ be a finite subset such that each $\Waff_{\lambda}$-orbit in $\Cl^{c,\lambda}(\haffd_{\bbR})$ contains at least one alcove from the family $\left\{ w^{-1} \nu_0 \right\}_{w\in \Sigma}$.
It is shown in~\cite{liu22} that the idempotent $1_\Sigma := \sum_{w\in \Sigma} 1_{c,w \lambda}\in \scrH^{c,\lambda}$ yields a Morita equivalence:
\begin{equation}\label{equa:morita}
	\scrH^{c,\lambda}\Mod \xrightarrow{\sim} 1_\Sigma\scrH^{c,\lambda}1_\Sigma\Mod,\quad M\mapsto 1_{\Sigma} M.
\end{equation}
Moreover, the natural map 
\[
	\scrZ^{c,\lambda}\to 1_\Sigma\scrH^{c,\lambda}1_\Sigma,\quad z\mapsto \sum_{w\in \Sigma}z_{c, w\lambda}
\]
is an isomorphism onto the centre of $1_\Sigma\scrH^{c,\lambda}1_\Sigma$. 
\begin{prop}\label{prop:Hmod}
	The following statements hold:
	\begin{enumerate}
		\item
			Let $M\in \scrH^{c,\lambda}\Mod$. Then, $M\in \scrH^{c,\lambda}\mof$ holds if and only if $1_{c, \lambda'}M$ is a finitely generated $\scrZ^{c,\lambda}$-module for each $\lambda'\in [\lambda]$. Moreover, $M\in \scrH^{c,\lambda}\mof_{\fl}$ holds if and only if $1_{c, \lambda'}M$ is finite-dimensional for each $\lambda'\in [\lambda]$.
		\item
			For $M, N\in \scrH^{c,\lambda}\mof$, the $\scrZ^{c, \lambda}$-module $\Ext^n_{\scrH^{c,\lambda}}(M, N)$ is finitely generated for every $n\in \bbZ$; moreover, it is finite-dimensional whenever $M\in\scrH^{c,\lambda}\mof_{\fl}$ or $N\in \scrH^{c,\lambda}\mof_{\fl}$ holds.
	\end{enumerate}
\end{prop}
\begin{proof}
	Making use of the Morita equivalence \eqref{equa:morita}, it suffices to show these properties for $1_{\Sigma}\scrH^{c,\lambda}1_{\Sigma}$. They hold since $1_\Sigma\scrH^{c,\lambda}1_\Sigma$ is of finite rank over the complete noetherian subalgebra $\scrZ^{c,\lambda}$.
\end{proof}
Similar statements hold for $\scrH^{\lambda}_c$ for every $(c, [\lambda])\in \frakP\times \haffd / \Waff$. 

\subsection{Spectral completion of modules}
\subsubsection{}
For $M\in \bfH\Mod$, define $\scrC(M) = \scrS\otimes_{\bfS}M$. We equip $\scrC(M)$ with an $\scrH$-module structure by the following rules:
\begin{align*}
	(f\otimes 1)(g\otimes m) &= fg\otimes m \\
	(f\otimes s_\alpha)(g\otimes m) &= f\cdot \pre{s_\alpha}g \otimes s_\alpha m + \bfc_{\alpha} \vartheta_{\alpha}(g)\otimes m
\end{align*}
for $f,g\in \scrS$, $m\in M$ and $\alpha\in \Deltaaff$.

\begin{defi}\label{defi:compl}
	The \emph{spectral completion} is the following functor:
	\[
		\scrC: \bfH\Mod\to \scrH\Mod,\quad M\mapsto \scrS\otimes_{\bfS} M.
	\]
\end{defi}
Alternatively, via the $\bfH$-bimodule structure on $\scrH$, we may write $\scrC(M) = \scrH\otimes_{\bfH} M$. It is an exact conservative functor by~\autoref{lemm:platitude}.  It induces $\scrC_c:\bfH_c\to \scrH_c\Mod$ for each $c\in \frakP$, which is also exact and conservative.
\subsubsection{} The decomposition \eqref{equa:decompMod} induces
\[
	\scrC = \bigoplus_{(c, [\lambda])\in \frakP\times \haffd/\Waff}\scrC^{c, \lambda},\quad \scrC^{c, \lambda}: \bfH\Mod\to \scrH^{c,\lambda}\Mod.
\]
Similarly, we have $\scrC_c = \bigoplus_{[\lambda]\in \haffd / \Waff} \scrC^{\lambda}_c$ with $\scrC_c^{\lambda}: \bfH_c\Mod\to \scrH_c^{\lambda}\Mod$.

\subsection{Comparison of derived categories}\label{subsec:derived}
\subsubsection{}
We denote by $\Db(\bfH) = \Db(\bfH\Mod)$ and $\Db(\scrH) = \Db(\scrH\Mod)$ the corresponding bounded derived category of modules. Let $\Db_{\rmO}(\bfH)$ denote the full subcategory of $\Db(\bfH)$ formed by the complexes $K$ such that $\rmH^i(K)\in \rmO(\bfH)$ for each $i\in \bbZ$; $\Db_{\fl}(\scrH)$, $\Db_{\rmO}(\bfH_c)$ and $\Db_{\fl}(\scrH_c)$ are defined similarly.
\begin{theo}\label{theo:DequivO}
	The spectral completion yields t-exact equivalences of categories:
	\[
		\scrC: \Db_{\rmO}(\bfH)\xrightarrow{\sim} \Db_{\fl}(\scrH),\quad 
		\scrC_c: \Db_{\rmO}(\bfH_c)\xrightarrow{\sim} \Db_{\fl}(\scrH_c).
	\]
	and for each $(c, [\lambda])\in \frakP\times \haffd / \Waff$:
	\[
		\scrC^{c,\lambda}: \Db_{\rmO_{c,\lambda}}(\bfH)\xrightarrow{\sim} \Db_{\fl}(\scrH^{c,\lambda}),\quad \scrC^{\lambda}_c: \Db_{\rmO_{\lambda}}(\bfH_c)\xrightarrow{\sim} \Db_{\fl}(\scrH^{\lambda}_c).
	\]
\end{theo}
\begin{proof}
	Let $c\in \frakP$. We prove only the statement for $\scrC_c$; the proof for $\scrC$, $\scrC^{c,\lambda}$ and $\scrC^{\lambda}_c$ is similar. \par
	We prove first that $\scrC_c(M)\in \scrH_c\mof_{\fl}$ for $M\in \rmO(\bfH_c)$. Given such $M$, we have $M = \bigoplus_{\lambda\in \haffd}M_{\lambda}$, where $M_{\lambda}$ is the generalised $\lambda$-eigenspace for the $\calO(\haffd)$-action. The condition $M\in \rmO(\bfH)$ implies that there is an isomorphism $\scrM := \scrC_c M \cong M$ of vector spaces and the action maps $\bfH_c\to \End_{\bbC}(M)$ and $\scrH_c\to \End_{\bbC}(\scrM)\cong\End_{\bbC}(M)$ have the same image; therefore, $\scrM$ is a compact $\scrH_c$-module. \par
	By \eqref{equa:decompMod}, we may decompose $\scrM$:
	\[
		\scrM = \bigoplus_{[\lambda]\in \haffd / \Waff} \scrM_{[\lambda]},\quad  \scrM_{[\lambda]} = \bigoplus_{\lambda'\in [\lambda]}M_{\lambda'}\in \scrH^{\lambda}_c\mof;
	\]
	the summation is finite due to the compactness of $\scrM$. The centre $\scrZ^{\lambda} = \calO(\haffd)^{\Waff_{\lambda}}$ acts locally finitely on $\scrM_{[\lambda]}$; therefore, $\scrM_{[\lambda]}$ lies in $\scrH_c^{\lambda}\mof_{\fl}$ by \autoref{prop:Hmod}. This proves that the spectral completion induces $\scrC_c: \Db_{\rmO}(\bfH_c)\to \Db_{\fl}(\scrH_c)$. \par

	We prove the full faithfulness. Given $M, N\in \rmO(\bfH_c)$, choose a free resolution $(P_{\bullet}, \delta_{\bullet})$ for $M$, where $P_{i} = \bfH_c^{\oplus n_i}$ and $\delta_i\in \Hom_{\bfH_c}(\bfH_c^{\oplus n_{i-1}}, \bfH_c^{\oplus n_{i}})$ for $i\ge 0$. We may write $\delta_i$ as right multiplication by a matrix $m_i\in \mathrm{M}_{n_i\times n_{i-1}}(\bfH_c)$.  Then, we have
	\[
		\Hom_{\bfH_c}(P_{\bullet}, N) = [ 0\to N^{\oplus n_0}\xrightarrow{m_0} N^{\oplus n_1}\xrightarrow{m_1} \cdots ]
	\]
	Via the generalised eigenspace decomposition $N = \bigoplus_{(c,\lambda)\in \frakP\times \haffd} N_{\lambda}$, we may write 
	\[
		m_i = \sum_{\lambda', \lambda''\in \haffd}m_{i,\lambda', \lambda''},\quad m_{i, \lambda', \lambda''}: N^{\oplus n_i}_{\lambda'}\to N^{\oplus n_{i+1}}_{\lambda''}\quad \text{for $\lambda',\lambda''\in \haffd$}.
	\]
	Then, it follows that for $i\in \bbN$, we have
	\begin{equation}\label{equa:ext1}
		\Ext^i_{\bfH_c}(M, N) = \bigoplus_{\lambda'\in \haffd}\left(\bigcap_{\lambda''\in \haffd} \ker(m_{i,\lambda',\lambda''}) / \sum_{\lambda''\in \haffd} \im(m_{i-1, \lambda'', \lambda'})\right).
	\end{equation}
	On the other hand, $(\scrC_c P_{\bullet}, \scrC_c\delta_{\bullet})$ is a projective resolution of $\scrC_c M$. We may write
	\[
		\Hom_{\scrH_c}(\scrC_c P_{\bullet}, \scrC_c N) = \left[0 \to \prod_{\lambda'\in \haffd}N_{\lambda'}^{\oplus n_0}\xrightarrow{(m_{0, \lambda', \lambda''})_{\lambda', \lambda''}}  \prod_{\lambda''\in \haffd}N_{\lambda''}^{\oplus n_1}\to\cdots\right].
	\]
	Therefore, 
	\begin{equation}\label{equa:ext2}
		\Ext^i_{\scrH_c}(\scrC_c M, \scrC_c N) = \prod_{\lambda'\in \haffd}\left(\bigcap_{\lambda''\in \haffd} \ker(m_{i,\lambda',\lambda''}) / \sum_{\lambda''\in \haffd} \im(m_{i-1, \lambda'', \lambda'})\right).
	\end{equation}
	The chain map $\Hom_{\bfH_c}(P_{\bullet}, N) \to \Hom_{\scrH_c}(\scrC_c P_{\bullet}, \scrC_c N)$ is given by the natural inclusion $\bigoplus_{\lambda'\in \haffd}\hookrightarrow\prod_{\lambda'\in \haffd}$. This inclusion induces an isomorphism $\Ext^i_{\bfH_c}(M, N)\cong \Ext^i_{\scrH_c}(\scrC_c M, \scrC_c N)$ for $i\in \bbZ$ since all but finitely many terms in the sum \eqref{equa:ext1} and product \eqref{equa:ext2} are zero by the finite dimensionality of $\Ext^i_{\scrH_c}(\scrC_c M, \scrC_c N)$ (see \autorefitem{prop:Hmod}{ii}). This proves the full faithfulness. \par
	It remains to show the essential surjectivity. We define a functor $F:\scrH_c\Mod\to \bfH_c\Mod$ as follows: given $M\in \scrH_c\Mod$, we equip $M$ with an $\bfH_c$-action by the formula: $h\cdot m = \sum_{\lambda\in \haffd} h 1_{c,\lambda}m$ for $h\in \bfH_c$ and $m\in M$; it is well-defined because $1_{c, \lambda}m$ is zero for all but finitely many $\lambda\in \haffd$. For $M\in \scrH_c\mof_{\fl}$, we have $F(M)\in \rmO(\bfH)$ by the non-degeneracy of $M$ as $\scrH_c$-module and by the finite dimensionality of $1_{\lambda}M = M_{\lambda}$ for each $\lambda\in \haffd$. It is easy to see that $\scrC_c F(M) \cong M$ for $M\in \rmO(\bfH_c)$; therefore, $\rmO(\bfH_c)$ lies in the essential image of $\scrC_c$. Now, given $\scrK\in \Db_{\fl}(\scrH_c)$, we may suppose up to cohomological shift that $\rmH^{n}(\scrK) = 0$ for $n < 0$; we may then prove by induction on $n\in \bbN$ that the truncation $\tau_{\le n} \scrK$ lies in $\Db_{\fl}(\scrH_c)$ via the exactness and full faithfulness of $\scrC_c$, whence the essential surjectivity.
\end{proof}
\begin{coro}\label{coro:fd}
	The objects in $\rmO(\bfH)$ are of finite length. Moreover, given any $M\in \rmO(\bfH)$, the generalised $(c,\lambda)$-weight space $M_{c,\lambda}\subseteq M$ is finite-dimensional for every $(c,\lambda)\in \frakP\times\haffd$. 
\end{coro}
\begin{proof}
	This is immediate from \autorefitem{prop:Hmod}{i} and \autoref{theo:DequivO}.
\end{proof}
\begin{coro}\label{coro:globdimscrH}
	The categories $\scrH\Mod$ and $\scrH_c\Mod$ for $c\in \frakP$ have finite global dimension. 
\end{coro}
\begin{proof}
	We prove only the statement for $\scrH\Mod$. Since $\scrH\Mod$ is compactly generated, it suffices to show that there exists $d\in \bbN$ such that $\Ext^d_{\scrH}(M, N) = 0$ holds for $n > d$ and for $M, N\in \scrH\mof$. For $M\in \scrH\mof$, we have 
	\[
		M \xrightarrow{\sim} \varprojlim_{M'\subseteq M} M / M',
	\]
	where $M'$ runs over submodules of $M$ such that $M/M'\in \scrH\mof_{\fl}$. It follows that
	\[
		\Ext^n_{\scrH}(M,N) \cong \varprojlim_{N'\subseteq N}\varinjlim_{M'\subseteq M}\Ext^n_{\scrH}(M/M',N/N')
	\]
	It results from \autoref{prop:dimglob} and \autoref{theo:DequivO} that there exists $d \in \bbN$ such that $\Ext^n_{\scrH}(M/M',N/N') = 0$ whenever $n > d$, whence $\Ext^n_{\scrH}(M,N) = 0$.
\end{proof}

\subsubsection{}
For $A \in \left\{ \bfH, \bfH_c, \scrH, \scrH_c\;;\; c\in \frakP \right\}$, let $\Db_{\perf}(A)$ be the full subcategory of $\Db(A) := \Db(A\Mod)$ formed by perfect complexes.
\begin{prop}
	For $A \in \left\{ \bfH, \bfH_c, \scrH, \scrH_c\;;\; c\in \frakP \right\}$, the natural functor $\Db(A\mof)\to \Db_{\perf}(A)$ is an equivalence of categories.
\end{prop}
\begin{proof}
	That the embedding $\Db(A\mof)\hookrightarrow \Db_{A\mof}(A\Mod)$ is an equivalence is a standard result in homological algebra. On the other hand, we have $\Db_{A\mof}(A) = \Db_{\perf}(A)$ by the finiteness of global dimension (\autoref{prop:dimglob} and \autoref{coro:globdimscrH}). 
\end{proof}
\begin{prop}\label{prop:Dfl}
	The natural functors $\Db(\scrH\mof_{\fl})\to \Db_{\fl}(\scrH)$ and $\Db(\scrH_c\mof_{\fl})\to \Db_{\fl}(\scrH_c)$ for $c\in \frakP$ are equivalences of categories.
\end{prop}
\begin{proof}
	We prove only the statement for $\scrH$. By a dual statement of \cite[13.17.4]{stacks-project}, this is a consequence of the following property: given $M\in \scrH\mof$ and $N\subseteq M$ with $N \in \scrH\mof_{\fl}$, there exists $M'\subseteq M$ satisfying $M'\cap N = 0$ and $M/M'\in \scrH\mof_{\fl}$ --- this follows from the Morita equivalence \eqref{equa:morita} and the Artin--Rees lemma applied to the centre $\scrZ^{c,\lambda}$ for every $(c,[\lambda])\in \frakP\times \haffd/\Waff$. 
\end{proof}
\begin{coro}\label{coro:DO}
	The natural functors $\Db(\rmO(\bfH_c))\to \Db_{\rmO}(\bfH_c)$ and $\Db(\rmO(\bfH))\to \Db_{\rmO}(\bfH)$ are equivalences of categories.
\end{coro}
\begin{proof}
	Consider the following commutative diagram of functors:
	\[
		\begin{tikzcd}[row sep = 15pt]
			\Db(\rmO(\bfH)) \arrow{r}\arrow{d}{\scrC} & \Db_{\rmO}(\bfH)\arrow{d}{\scrC} \\
		\Db(\scrH\mof_{\fl}) \arrow{r} & \Db_{\fl}(\scrH)
		\end{tikzcd}
	\]
	By \autoref{theo:DequivO} and \autoref{prop:Dfl}, the vertical arrows as well as the lower horizonal arrow are equivalences of categories; therefore, the upper horizontal arrow is also an equivalence. The statement for $\scrC_c$ is similar. 
\end{proof}
\begin{rema}
	Given an abelian category $\calA$ and a Serre subcategory $\calB\subseteq \calA$, the natural functor $\Db(\calB)\to \Db_{\calB}(\calA)$ may not be an equivalence in general.
\end{rema}

\subsection{Length filtration}\label{subsec:lg}
We define another filtration $\{ F^{\lg}_k\bfH\}_{k\in \bbN}$ on $\bfH$, called \emph{length filtration}, by setting
\[
	F^{\lg}_k \bfH = \bigoplus_{\substack{w\in \Waff \\ \ell(w) \le k}} w\bfS.
\]
For each $w\in \Waff$, let $\ba w$ denote its image in $\gr^{\lg}_{\ell(w)}\bfH$. Then we have 
\[
	\gr^{\lg}_k \bfH = \bigoplus_{\substack{w\in \Waff \\ \ell(w)= k}} \ba w\bfS,
\]
which is an $\bfS$-bimodule, finitely generated both as left and right $\bfS$-module.
\begin{prop}\label{prop:lgO}
	Let $M\in \bfH\mof$. Then $M$ lies in $\rmO(\bfH)$ if and only if for some (equiv. for every) good filtration $\left\{ F_{k}M \right\}_{k\in \bbZ}$ with respect to the length filtration of $\bfH$, we have $\dim F_{k}M < \infty$ for every $k\in \bbZ$. 
\end{prop}
\begin{proof}
	Suppose the good filtration $F_\bullet M$ satisfies the condition of finiteness. Then, as the action of $\bfS$ on $M$ preserves each degree $F_{k} M$, the action of $\bfS$ on $M$ is locally finite; hence $M\in \rmO(\bfH)$. Conversely, suppose that $M\in \rmO(\bfH)$ and $\left\{F_k M  \right\}_{k\in \bbZ}$ is a good filtration with respect to the length filtration $F^{\lg}_\bullet\bfH$. Since $\gr^F M$ is finitely generated over $\gr^{\lg}\bfH$, we can find a finite-dimensional generating graded subspace $V\subset \gr^F M$. Hence we have
	\[
		\gr^F_k M = \bigoplus_{j\in \bbZ} \gr^{\lg}_j \bfH(V\cap \gr^F_{k-j} M) = \bigoplus_{j\in \bbN}\bigoplus_{\ell(w) = j}[w]\bfS (V\cap \gr^F_{k-j} M).
	\]
	By the hypothesis that $M\in \rmO(\bfH)$, \autoref{coro:fd} implies that the space $\bfS (V\cap \gr^F_{k-j} M)$ is finite-dimensional. Since $\gr^F_{k-j} M$ vanishes for $j \gg 0$ by the finite generation, it follows that $\dim \gr^F_k M < \infty$ and, consequently, $\dim F_k M < \infty$. 
\end{proof}

\subsection{Algebraic Knizhnik--Zamolodchikov functor}\label{subsec:KZ}
We review the ``algebraic Knizhnik--Zamolodchikov (KZ) functor'' introduced in~\cite{liu22}. 
\subsubsection{}
Let $(c, \lambda)\in \frakP\times\haffd$. Consider the algebra $\scrH^{\lambda}_c$ introduced in \autoref{subsubsec:Hcl}. 
It has been shown in {\it loc. cit.} that there exists an idempotent $1_{\bfV}\in \scrH^{\lambda}_c$ such that the idempotent subalgebra $1_{\bfV}\scrH^{\lambda}_c1_{\bfV}$ is isomorphic to a block algebra of the affine Hecke algebra. The algebraic KZ functor is defined to be the idempotent truncation $M\mapsto 1_{\bfV}M$ for $M\in \scrH^{\lambda}_c\mof$. \par
Let us recall the affine Hecke algebra and its completion. Let $P = P_{\Rfin} \subseteq \hfin^*_{\bbR}$ be the weight lattice of the root system $(\hfin_{\bbR}, \Rfin)$. The \emph{extended affine Hecke algebra} (for the dual root system $(\hfin_{\bbR}^*, R^{\vee})$), denoted $\bfK$, is the associative algebra over the Laurent polynomial ring $\calO_{\bft} = \bbC[\bft^{\pm 1}_{\ast}, \ast\in R/W]$ generated by the two sets $\left\{ T_{\alpha} \right\}_{\alpha\in \Deltafin}$ and $\left\{ X^{\mu} \right\}_{\mu\in P}$ modulo the following relations for $\mu,\nu\in P$ and $\alpha,\beta\in \Delta$:
\[
	X^{0} = 1,\quad X^{\mu+ \nu} = X^\mu  X^\nu,\quad \underbrace{T_{\alpha}T_{\beta}T_{\alpha}\cdots}_{m_{\alpha,\beta}} =  \underbrace{T_{\beta}T_{\alpha}T_{\beta}\cdots}_{m_{\alpha,\beta}}, 
\]
and if $(\frakh_{\bbR}, R)$ is of type $BC_n$ with $n\ge 1$, with the notation of~\autoref{exam:BC}:
\begin{align*}
	&\begin{cases}(T_{\alpha} - \bft_{\natural})(T_{\alpha} + 1) = 0 &\text{for $\alpha\in \left\{ \alpha_1, \ldots, \alpha_{n-1} \right\}$}, \\ (T_{\alpha} - \bft_{\sharp}\bft_{\flat})(T_{\alpha} + 1) = 0 & \text{for $\alpha = \alpha_n$}, \end{cases}\\
	&T_{\alpha}X^{\mu} - X^{s_{\alpha}(\mu)}T_{\alpha} = \begin{cases} (\bft_{\natural} - 1)\frac{X^{\mu} - X^{s_{\alpha}(\mu)}}{1 - X^{-\alpha}} & \text{for $\alpha\in \left\{ \alpha_1, \ldots, \alpha_{n-1} \right\}$}, \\ ((\bft_{\sharp}\bft_{\flat} - 1) + (\bft_{\sharp} - \bft_{\flat})X^{\alpha})\frac{X^{\mu} - X^{s_{\alpha}(\mu)}}{1 - X^{-2\alpha}} & \text{for $\alpha = \alpha_n$};\end{cases}
\end{align*}
otherwise:
\[
	(T_{\alpha} - \bft_{\alpha})(T_{\alpha} + 1) = 0,\quad T_{\alpha}X^{\mu} - X^{s_{\alpha}(\mu)}T_{\alpha} = (\bft_{\alpha} - 1)\frac{X^{\mu} - X^{s_{\alpha}(\mu)}}{1 - X^{-\alpha}},\quad \text{for $\alpha\in \Delta$}.
\]

Given a family of parameters $t=(t_{*})_{*\in R/W}$ in $\bbC^{\times}$, we define $\bfK_t = \bfK \otimes_{\calO_{\bft}} \calO_{\bft} / (\bft = t)$ to be the specialisation of parameters. \par

\subsubsection{}
	The elements $\left\{ X^{\mu} \right\}_{\mu\in P}$ generates inside $\bfK$ a Laurent polynomial subalgebra isomorphic to the group ring $\bbC P$. Let $T = \Spm \bbC P$ be its maximal spectrum. For each $ W $-orbit $[\ell] \in T /  W $, let $\bfK_t\mof_{\ell}\subset \bfK_t\mof$ be the full subcategory consisting of finitely generated $\bfK_t$-modules $M$ which admit the following decomposition:
	\[
		M = \bigoplus_{\ell'\in [\ell]}M_{\ell'},\quad M_{\ell'} = \bigcup_{n \ge 0} \left\{ x\in M\;;\;  \frakm_{\ell'}^n x = 0 \right\}.
	\]
	According to a well-known result of Bernstein, the ring of invariant functions $\calO(T)^W$ coincides with the centre of $\bfK_t$. Let $\scrK_t^{\ell} = (\calO(T)^ W )^{\wedge}_{[\ell]}\otimes_{\calO(T)^ W } \bfK_t$. The exact functor
	\[
		(\calO(T)^ W )^{\wedge}_{[\ell]}\otimes_{\calO(T)^ W }\relbar : \bfK_t\mof\to \scrK_t^{\ell}\mof
	\]
	induces an equivalence on subcategories $\bfK_t\mof_{\ell}\xrightarrow{\sim} \scrK_t^{\ell}\mof_{\mathrm{fd}}$. \par
\subsubsection{}
Given $c\in \frakP$ and $\lambda = (\lambda^0, 1)\in \haffd$, consider the spectral completion $\scrH^{\lambda}_c$. Choose an element $\gamma$ in the coroot lattice $Q^{\vee} = \bbZ R^{\vee}$ such that $\langle \alpha, \gamma\rangle \ll 0$ for every $\alpha\in R^+$. In particular, $w^{-1}\nu_0 - w^{-1}\gamma$ lies in a generic $(c,\lambda)$-clan (\autoref{defi:generic}) for each $w\in W$. Set $1_{\bfV} = \sum_{[w]\in  W  /  W _{\ell}} 1_{c, w\lambda + \gamma}$, where $\ell = \exp(2\pi i\lambda^0)\in T$ and $ W_{\ell} = \Stab_{ W }(\ell)$. Moreover, we set $t_* =  \exp(2\pi i c_{*})$ for $*\in \Rfin/ W $. 
\begin{theo}[\cite{liu22}]\label{theo:KZ}
		Under the above assumptions, there exists an isomorphism $1_{\bfV}\scrH^{\lambda}_c1_{\bfV} \cong \scrK_t^{\ell}$. Moreover, the following functor 
		\[
			\bfV: \scrH^{\lambda}_c\Mod \to \scrK_t^{\ell}\Mod,\quad M\mapsto 1_{\bfV} M
		\]
		is a quotient functor and satisfies the double centraliser property: the restriction of $\bfV$ to the subcategory of compact projective objects $\scrH^{\lambda}_c\proj\subset \scrH^{\lambda}_c\Mod$ is fully faithful. \hfill\qedsymbol
	\end{theo}
	The functor $\bfV$ above is called the algebraic KZ functor. 

	\subsubsection{}
Let $\frakP_{\bbZ} = \mathop{\mathrm{Map}}(R/W, \bbZ)\subseteq \frakP$ be the set of maps from $R/W$ to $\bbZ$. Suppose we are given $c, c'\in \frakP$ such that their difference $d := c' - c$ lies in $\frakP_{\bbZ}$; set $t = (t_*)_{*\in R/W}$ with $t_\ast =  e^{2\pi i c_{\ast}} = e^{2\pi i c'_{\ast}}$. By~\autoref{theo:KZ}, we obtain the respective algebraic KZ functors:
	\[
		\scrH^{\lambda}_c\Mod \xrightarrow{\bfV_c} \scrK_t^{\ell}\Mod\xleftarrow{\bfV_{c'}}\scrH^{\lambda}_{c'}\Mod.
	\]
	In~\autoref{sec:transl}, we will construct the translation functor $\trans{c'}{c}$ which intertwines the derived functors of $\bfV_c$ and $\bfV_{c'}$.

\section{Harish-Chandra bimodules}\label{sec:HC}
In this section, we introduce the notion of Harish-Chandra (HC) bimodules for the trigonometric DAHA $\bfH$ and establish their basic properties.
\subsection{Category of HC bimodules}
Let $(A, \left\{ F_k A \right\}_{k\in \bbN})$ and $(A', \left\{ F_k A' \right\}_{k\in \bbN})$ be filtered rings such that $\gr^F A$ and $\gr^F A'$ are left and right noetherian. Let $M$ be an $(A, A')$-bimodule. We say that a filtration $\left\{ F_k M \right\}_{k\in \bfZ}$ on $M$ is \emph{good} if it is exhaustive and separable, $(F_j A)(F_k M)(F_l A') \subseteq F_{j+k+l} M$ holds for each $k,j,l\in \bbZ$ and $\gr^F M$ is finitely generated as $(\gr^F A, \gr^F A')$-bimodule; we say that $\left\{ F_k M \right\}_{k\in \bfZ}$ is \emph{excellent} if it is good and $\gr^F M$ is finitely generated both as left $(\gr^F A)$-module and as right $(\gr^F A')$-module.
\begin{defi}\label{defi:HC}
	Let $M$ be a $(\bfH, \bfH)$-bimodule. A filtration $\{F_{k} M\}_{k\in \bbZ}$ on $M$ is called a \emph{Harish-Chandra filtration} (HC filtration) if it is a good filtration with respect to the canonical filtration $F^{\can}_k\bfH$ and the action of $\ad z$ on $\gr^{F}M$ is nilpotent for every $z\in \rmZ(\gr^F \bfH)$. We say $M$ is a \emph{Harish-Chandra bimodule} (HC bimodule) if it admits a HC filtration.
\end{defi}

\begin{prop}\label{prop:HCexcellent}
	Every HC filtration is excellent. Consequently, if $M$ is a HC $\bfH$-bimodule. Then $M$ is finitely generated both as left and right $\bfH$-module.
\end{prop}
\begin{proof}
See~\cite[3.4.3]{losev12}. 
\end{proof}
Let $\HC(\bfH)$ denote the category of HC $\bfH$-bimodules. Given $M\in \HC(\bfH)$, if $\left\{ F_k M \right\}_{k\in \bbZ}$ and $\left\{ F'_k M \right\}_{k\in \bbZ}$ are good filtrations on $M$, then there exists $a\gg 0$ such that $F_{k-a} M \subseteq F'_{k}M \subseteq F_{k+a}M$ for every $k\in \bbZ$, see~\cite[D.1.3]{hotta2007d}. Consequently, every good filtration on $M$ is HC and therefore excellent by~\autoref{prop:HCexcellent}.
\begin{prop}
	$\HC(\bfH)$ is a Serre subcategory of the category of $\bfH$-bimodules.
\end{prop}
\begin{proof}
	It is obvious that the induced filtration of a HC filtration on a sub-bimodule or a quotient bimodule is again a HC filtration; thus $\HC(\bfH)$ is closed under formation of sub-objects and quotient objects. Suppose that $M$ is an $\bfH$-bimodule, $M'\subset M$ is a sub-$\bfH$-bimodule and $M'' = M / M'$ such that $M'$ and $M''$ are in $\HC(\bfH)$. Pick a good filtration $\left\{ F_k M \right\}_{k\in \bbZ}$ for $M$ and let $\left\{ F_k M' \right\}_{k\in \bbZ}$ and $\left\{ F_k M'' \right\}_{k\in \bbZ}$ be the induced filtrations. We have a short exact sequence
	\begin{equation}\label{equa:grB}
		0\to \gr^F M'\to \gr^F M\to \gr^F M''\to 0.
	\end{equation}
	Since $M'$ and $M''$ are HC and since $\left\{ F_k M' \right\}_{k\in \bbZ}$ and $\left\{ F_k M'' \right\}_{k\in \bbZ}$ are good (using~\eqref{equa:grB}), it follows that they are also HC filtrations. We see that the filtration $\{F_k M\}$ is also HC by using~\eqref{equa:grB} again.
\end{proof}

\begin{prop}\label{prop:extHC}
	Let $M,N\in \HC(\bfH)$. Then, $\Tor^{\bfH}_i(M, N), \Ext^i_{\bfH}(M, N)\in \HC(\bfH)$ hold for each $i\in \bbZ$.
\end{prop}
\begin{proof}
	Let $\{ F_k M \}_{k\in \bbZ}$ and $\{ F_k N \}_{k\in \bbZ}$ be HC filtrations on $M$ and $N$ respectively. Consider the Rees bimodules 
	\[
		M_{\delta} = \bigoplus_{k\in \bbZ} \delta^k F_k M,\quad N_{\delta} = \bigoplus_{k\in \bbZ} \delta^k F_k N.
	\]
	They are finitely generated graded $\bfH_{\delta}$-bimodules flat over $\bbC[\delta]$. The extension group $\Ext^k_{\bfH_{\delta}}(M_{\delta}, N_{\delta})$ is a finitely generated graded $\bfH_{\delta}$-bimodule. We have an exact sequence
	\[
		\Ext^k_{\bfH_{\delta}}(M_{\delta}, N_{\delta})\xrightarrow{\delta\cdot}\Ext^k_{\bfH_{\delta}}(M_{\delta}, N_{\delta})\to \Ext^k_{\bfH_{\delta}}(M_\delta, N_{\delta}/(\delta)) .
	\]
	By the flatness of $M_{\delta}$ over $\bbC[\delta]$, there is an isomorphism:
	\[
		\Ext^k_{\bfH_{\delta}}(M_\delta, N_{\delta}/(\delta)) \cong \Ext^k_{\bfH_{\delta}/(\delta)}\left( M_{\delta}/(\delta), N_{\delta}/(\delta)\right)= \Ext^k_{\gr^{\can}\bfH}(\gr^F M, \gr^F N).
	\]
	We obtain an embedding
	\[
		\Ext^k_{\bfH_{\delta}}(M_{\delta}, N_{\delta}) / (\delta) \hookrightarrow \Ext^k_{\gr^{\can}\bfH}(\gr^F M, \gr^F N).
	\]
	 The $\bbZ$-grading on $\Ext^k_{\bfH_{\delta}}(M_{\delta}, N_{\delta})$ induces a filtration on the quotient
	\[
		\Ext^k_{\bfH_{\delta}}(M_{\delta}, N_{\delta}) / (\delta - 1) \cong \Ext^k_{\bfH}(M, N),
	\]
	denoted by $F$. It is easy to show that the associated graded $\gr^F\Ext^k_{\bfH}(M, N)$ is a quotient of $\Ext^k_{\bfH_{\delta}}(M_{\delta}, N_{\delta})/(\delta)$. By the assumption, for every $z\in Z(\gr^{\can}\bfH)$, the operator $(\ad z)^d$ vanishes on $\gr^F M$ and $\gr^F N$ when $d \gg 0$; hence, it is also the case on $\Ext^k_{\gr^{\can}\bfH}(\gr^F M, \gr^F N)$ and on the subquotient bimodule $\gr^F\Ext^k_{\bfH}(M, N)$. It follows that $F$ is a HC filtration on $\Ext^k_{\bfH}(M, N)$ and the latter is a Harish-Chandra $\bfH$-bimodule. The proof for $\Tor^{\bfH}_k(M, N)$ is similar. 
\end{proof}

Let $\Db_{\HC}(\bfH)$ be the triangulated subcategory of $\Db(\bfH\otimes \bfH^{\op}\Mod)$ consisting of complexes $K$ such that $\rmH^k(K)\in \HC(\bfH)$ for every $k\in \bbZ$.
\begin{coro}\label{coro:extHC}
	Let $K,L\in \Db_{\HC}(\bfH)$. Then $K\otimes^{\rmL}_{\bfH}L$ and $\RHom_{\bfH}(K, L)$ are also in $\Db_{\HC}(\bfH)$. 
\end{coro}
\begin{proof}
	It results immediately from~\autoref{prop:dimglob} and~\autoref{prop:extHC}. 
\end{proof}

\subsection{Functors from HC bimodules}
For $B\in (\bfH\otimes \bfH^{\op})\Mod$, we have the derived tensor product functor and the derived Hom functor: 
\[
	B\otimes^{\rmL}_{\bfH}\relbar:\Db(\bfH\Mod)\to \Db(\bfH\Mod),\quad \RHom_{\bfH}(B, \relbar):\Db(\bfH\Mod)\to \Db(\bfH\Mod). 
\]

\begin{prop}\label{prop:ExtTorO}
	For $B\in \HC(\bfH)$, the derived functors $B\otimes^{\rmL}_{\bfH}\relbar$ and $\RHom_{\bfH}(B, \relbar)$ preserve the subcategories $\Db_{\perf}(\bfH)$ and $\Db_{\rmO}(\bfH)$. 
\end{prop}
\begin{proof}
	We prove the statements for $\RHom_{\bfH}(B, \relbar)$ and leave those for $B\otimes^{\rmL}_{\bfH}\relbar$ to the reader. By~\autoref{coro:extHC}, we have $\RHom_{\bfH}(B, \bfH)\in \Db_{\HC}(\bfH)$. Since HC bimodules are finitely generated left $\bfH$-modules, we have $\RHom_{\bfH}(B, \bfH)\in \Db_{\perf}(\bfH)$ when the right $\bfH$-module structure is forgotten. Since $\bfH$ generates $\Db_{\perf}(\bfH)$ as thick subcategory, we see that $\RHom_{\bfH}(B, M)\in \Db_{\perf}(\bfH)$ for each $M\in \Db_{\perf}(\bfH)$. \par
	Let $M\in \rmO(\bfH)$. We pick a good filtration $\left\{ F_k M \right\}_{k\in \bbZ}$ for $M$ and, by~\autoref{lemm:filtB} below, an excellent filtration $\left\{ F_k B \right\}_{k\in \bbZ}$ for $B$ with respect to the length filtration $F^{\lg}_\bullet \bfH$ introduced in~\autoref{subsec:lg}. Set 
	\[
		\bfH_{\zeta} = \bigoplus_{k\in \bbN}\zeta^k F^{\lg}_k \bfH,\quad B_{\zeta} = \bigoplus_{k\in \bfZ}\zeta^k F_k B,\quad M_{\zeta} = \bigoplus_{k\in \bfZ}\zeta^k F_k M.
	\]

By~\autoref{prop:lgO}, $F_k M$ is a finite-dimensional $\bfS$-module for each $k\in \bbZ$. Moreover, $B_{\zeta}$ is a $(\bfH_{\zeta}, \bfH_{\zeta})$-bimodule and finite as left $\bfH_{\zeta}$-module. By taking a resolution of $B_{\zeta}$ by finite graded-free left $\bfH_{\zeta}$-modules, we deduce easily that $\Ext^k_{\bfH_{\zeta}}(B_{\zeta}, M_{\zeta})$ is a graded vector space and finite-dimensional in each degree and for each $k\in \bbZ$; moreover, it has the structure of graded left $\bfH_{\zeta}$-module coming from the right $\bfH_{\zeta}$-module structure on $B_{\zeta}$. By the finite dimensionality, $\bfS$ acts locally finitely on $\Ext^k_{\bfH_{\zeta}}(B_{\zeta}, M_{\zeta})$; this implies that the $\bfS$-action on the quotient $\Ext^k_{\bfH_{\zeta}}(B_{\zeta}, M_{\zeta}) / (\zeta - 1) \cong \Ext^k_{\bfH}(B, M)$ is also locally finite. It follows that $\Ext^k_{\bfH}(B, M)\in \rmO(\bfH)$ for each $k\in \bbZ$.
\end{proof}

\begin{lemm}\label{lemm:filtB}
	Let $B\in \HC(\bfH)$. Then, there exists an excellent filtration $\left\{ F_k B \right\}_{k\in \bbZ}$ for $B$ with respect to the length filtration $F^{\lg}_{\bullet}\bfH$.
\end{lemm}
\begin{proof}
The proof is modeled on that of~\cite[5.4.3]{losev12}. We first show that for each finite left $\bfS$-submodule $M\subseteq B$, the product $N = M\bfS$ remains finite as left $\bfS$-module. Let $\left\{ F'_k B \right\}_{k\in \bfZ}$ be a HC filtration (\autoref{defi:HC}) on $B$ and set $F'_k M = M \cap F'_k B$, $F'_k N = N \cap F'_k B$. Then $\gr^{F'} N = (\gr^{F'} M)(\gr^{\can} \bfS)$ holds, where $\gr^{\can}\bfS$ is the associated graded of restriction of the canonical filtration (\autoref{subsubsec:filtcan}) to $\bfS$. Since the adjoint action of $(\gr^{\can} \bfS)^ W $ on $\gr^{F'} B$ is nilpotent and since $\gr^{\can} \bfS$ is finite over $(\gr^{\can} \bfS)^ W $, it follows that $\gr^{F'} N$ is also finite over $\gr^{\can} \bfS$. Hence, $N$ is finite over $\bfS$. \par
Recall the length filtration $F^{\lg}_{\bullet}\bfH$ from~\autoref{subsec:lg}. Now, let $V\subseteq B$ be a finite-dimensional generating subspace as left $\bfH$-module. Let $\left\{ \ba b_j \right\}_{j = 1}^r$ be a spanning set for $\gr^{\lg}\bbC \Waff$ as left $\gr^{\lg}(\bbC Q^{\vee})^{ W }$-module. We assume $\ba b_j$ is homogeneous of degree $d'_j > 0$ except that $\ba b_1 = 1$ and we choose a lifting $b_j$ of $\ba b_j$ in $(\bbC Q^{\vee})^ W $. Let $\left\{ \ba z_k \right\}_{k=1}^s$ be a set of homogeneous $\bbC$-algebra generators for $\gr^{\lg} (\bbC Q^{\vee})^ W $ of degree $\deg\ba z_k = d_k$ and we choose a lifting $z_k$ of $\ba z_k$ in $(\bbC Q^{\vee})^ W $. Define
	\[
	F_n B = \sum_{j=1}^r\sum_{\substack{n_0, n_1, \ldots, n_s\in \bbN \\ n_0 + d_1n_1 + \cdots + d_sn_s + d'_j \le n}} (F^{\lg}_{n_0}\bfH)\left((\ad z_s)^{n_s}\cdots (\ad z_1)^{n_1}V\right)b_j\bfS
	\]
	Let's show that $\gr^F B$ is finitely generated as left $\gr^{\lg} \bfH$-module. Since the sum
	\[
		V' = \sum_{j=1}^r\sum_{n_1, \ldots, n_s\in \bbN} \bfS\left((\ad z_s)^{n_s}\cdots (\ad z_1)^{n_1}V\right)b_j\bfS
	\]
	is a finite left $\bfS$-module by the first paragraph, we can find a finite subset of homogeneous polynomials $I\subset \bfS$ such that the finite-dimensional subspace
	\[
		V'' = \sum_{f\in I}\sum_{j=1}^r\sum_{n_1, \ldots, n_s\in \bbN} \bbC \left((\ad z_s)^{n_s}\cdots (\ad z_1)^{n_1}V\right)b_jf\subseteq V'
	\]
	spans $V'$ as left $\bfS$-module. Let $F_nV'' = V''\cap F_nB$; then $\gr^F V''$ clearly generates $\gr^FB$ as left $\gr^{\lg}\bfH$-module. Thus $\gr^F B$ is a finite left $\gr^{\lg}\bfH$-module. \par
	Let's show that $(F_k B)(F^{\lg}_l \bfH)\subseteq (F_{k+l} B)$. Clearly, it is enough to show that for $h\in F_k\bfH$, $j \in\{ 1, \ldots, r\}$ and $v_0\in V$, we have
	\[
		((\ad z_s)^{n_s} \cdots (\ad z_1)^{n_1}v_0)b_j h \in F_{k+n} B,\quad \text{where $n = n_1d_1 + \cdots + n_s d_s + d'_j$}.
	\]
	We prove this statement by induction on $n$. When $n = 0$, we must have $d'_j = 0$, so $b_j = b_1 = 1$. In other words, we need to prove that $v_0 h \in F_k B$. We may assume that $h$ is of the form $z b_{j'}$, where $z = z_1^{m_1}\cdots z_s^{m_s}$ is a monomial and $j'\in \left\{ 1, \ldots,r \right\}$ such that $m_1d_1 + \cdots + m_sd_s + d'_j \le k$. One can easily prove the following formula by induction:
	\[
		v_0h = \sum_{\substack{0\le m'_1\le m_1 \\ \cdots \\ 0\le m'_s\le m_s}}(-1)^{m'_1 + \cdots + m'_s}\binom{m_1}{m'_1}\cdots\binom{m_s}{m'_s} z_1^{m_1 - m'_1} \cdots z_s^{m_s - m'_s}\left((\ad z_s)^{m'_s}\cdots(\ad z_1)^{m'_1}v_0\right)b_{j'}. 
	\]
	Since $z_1^{m_1 - m'_1} \cdots z_s^{m_s - m'_s} \in F_{d_1(m_1 - m'_1) + \cdots + d_s(m_s - m'_s)} \bfH$, we deduce $v_0h\in F_k B$. This is the initial step of the induction. Suppose $n> 0$ and the statement has been proven for smaller $n$. If $j\neq 1$, we have $d'_j > 0$; the induction hypothesis applied to $b_jh$ in place of $h$ yields $((\ad z_s)^{n_s} \cdots (\ad z_1)^{n_1}v_0)(b_jh)\in F_{k+n} B$. Suppose now that $j = 1$. Let $i_0 = \max\left\{ i\;;\; n_i > 0 \right\}$. Then
	\begin{align*}
		((\ad z_{i_0})^{n_{i_0}} \cdots (\ad z_1)^{n_1}v_0)h = & z_{i_0}((\ad z_{i_0})^{n_{i_0}-1} \cdots (\ad z_1)^{n_1}v_0) h \\
		&- ((\ad z_{i_0})^{n_{i_0}-1} \cdots (\ad z_1)^{n_1}v_0)z_{i_0} h.
	\end{align*}
	Then, the induction hypothesis yields the result. It follows that $F_\bullet B$ is an $\bfH$-bimodule filtration with respect to $F^{\lg}_\bullet \bfH$ and is good with respect to the left $\bfH$-action. Analogously, we can construct an $\bfH$-bimodule filtration $F'_\bullet B$ with respect to $F^{\lg}_\bullet \bfH$ and is good with respect to the right $\bfH$-action. Since both filtrations are good over $\bfH \otimes \bfH^{\op}$ with respect to the filtration $F_{\bullet}^{\zeta}\otimes F_{\bullet}^{\zeta}$, there exists $a \gg 0$ such that $F'_{n-a} B \subseteq F_n B \subseteq F'_{n+a}B$ for every $n\in \bbZ$, see~\cite[D.1.3]{hotta2007d}. It follows that $F_\bullet B$ is also good for the right $\bfH$-action, thus excellent.
\end{proof}

\subsection{Spectral completion of HC bimodules}
\begin{lemm}\label{lemm:BSSB}
	For each $\bfS$-bimodule $B$ which is finite both as left and as right $\bfS$-module, there are natural isomorphisms of $\bfS$-bimodules
	\[
		B\otimes_{\bfS}\scrS\cong \bigoplus_{\substack{(c,\lambda),(c',\lambda')\in \frakP\times\haffd}}\varprojlim_{k,l} B / \left(B\frakm^k_{c,\lambda} + \frakm_{c',\lambda'}^{l} B\right)\cong \scrS\otimes_{\bfS} B.
	\]
\end{lemm}
\begin{proof}
	The proof, which we leave to the reader, is an easy exercise of commutative algebra involving the Chinese remainder theorem.
\end{proof}
Let $\eta_{B}: B\otimes_{\bfS}\scrS\xrightarrow{\sim} \scrS\otimes_{\bfS}B$ denote the isomorphism given by~\autoref{lemm:BSSB}. The naturality of $\eta_{B}$ implies that $\eta_B$ can be extended to those $B$ which can be written as colimit of modules satisfying the condition of~\autoref{lemm:BSSB}. Moreover, this isomorphism satisfies the property $\eta_{B\otimes_{\bfS}B'} = (\eta_{B}\otimes \id_{B'}) \circ (\id_B\otimes \eta_{B'})$. \par
By~\autoref{lemm:filtB}, every HC bimodule $B\in \HC(\bfH)$ admits a filtration $\left\{ F_k B \right\}$ such that $F_k B$ is an $\bfS$-bimodule satisfying the condition of~\autoref{lemm:BSSB} for each $k\in \bbZ$. 

The multiplication of $\scrH$ can be alternatively described as follows:
\begin{align*}
	&\scrH\otimes_{\bfS}\scrH = \scrS\otimes_{\bfS}\bfH\otimes_{\bfS}\scrS\otimes_{\bfS}\bfH\\
	&\xrightarrow{\id_{\scrS}\otimes \eta_{\bfH}\otimes \id_{\bfH}}\scrS\otimes_{\bfS}\scrS\otimes_{\bfS}\bfH\otimes_{\bfS}\bfH \xrightarrow{\mu_{\scrS}\otimes \mu_{\bfH}}\scrS\otimes_{\bfS}\bfH = \scrH,
\end{align*}
where $\mu_{\scrS}$ and $\mu_{\bfH}$ are the multiplication map of $\scrS$ and $\bfH$, respectively. The $\scrH$-module structure on $\scrC M$ for $M\in \bfH\Mod$ can be described in a similar way.

\begin{prop}\label{prop:Hbimod}
	Let $B\in \HC(\bfH)$. Then, the $\bfH$-bimodule structure on $B$ induces an $\scrH$-bimodule structure on the spectral completion $\scrC B$.
\end{prop}
\begin{proof}
	The left $\scrH$-module structure on $\scrC B$ is described in the previous paragraph. The right $\scrH$-module structure is defined as follows:
	\begin{align*}
		&\scrC B\otimes_{\bfS}\scrH = \scrS\otimes_{\bfS}B\otimes_{\bfS}\scrS\otimes_{\bfS}\bfH \\
		\xrightarrow{\id_{\scrS}\otimes \eta_B\otimes \id_{\bfH}}&\scrS\otimes_{\bfS}\scrS\otimes_{\bfS}B\otimes_{\bfS}\bfH \xrightarrow{\mu_{\scrS}\otimes \rho_B} \scrS\otimes_{\bfS}B = \scrC B,
	\end{align*}
	where $\rho_B$ is the right $\bfH$-action on $B$. It is easy to verify that the left and right actions on $B$ commute.
\end{proof}
For every $P\in \left(\scrH \otimes \scrH^{\op}\right)\Mod$, let $\Rhom_{\scrH}(P, \relbar)$ be the derived functor of the following functor:
\[
	\ho_{\scrH}(P, \relbar) : \scrH\Mod \to \scrH\Mod,\quad  N\mapsto \bigoplus_{(c, \lambda)\in \frakP\times \haffd} \Hom_{\scrH}(P 1_{c, \lambda}, N).
\]
Its left adjoint functor is given by $P\otimes^{\rmL}_{\scrH} \relbar: \Db(\scrH) \to \Db(\scrH)$.

\begin{prop}\label{prop:HCTor}
	Let $B\in\HC(\bfH)$. Then, the following statements hold:
	\begin{enumerate}
		\item\label{prop:HCTor-i}
			The functor $\scrC (B)\otimes^{\rmL}_{\scrH}\relbar$ preserves the subcategories $\rmD^b_{\perf}(\scrH)$ and $\rmD^b_{\fl}(\scrH)$.
		\item\label{prop:HCTor-ii}
			There is a natural isomorphism
			\[
				\scrC(B\otimes^{\rmL}_{\bfH} M)\cong \scrC(B)\otimes^{\rmL}_{\scrH}\scrC(M) 
			\]
			for $M\in \bfH\Mod$.
	\end{enumerate}
\end{prop}
\begin{proof}
	For $M\in \bfH\Mod$, we have isomorphisms
	\begin{align*}
		&\scrC(B)\otimes_{\scrH}\scrC(M) \cong B\otimes_{\bfH}\scrH\otimes_{\scrH}\scrH\otimes_{\bfH}M \\
		&\cong B\otimes_{\bfH}\scrH\otimes_{\bfH}M\cong \scrH\otimes_{\bfH}B\otimes_{\bfH}M = \scrC(B\otimes_{\bfH}M).
	\end{align*}
	They are bi-functorial in $B$ and $M$. Taking the derived functors, we obtain the quasi-isomorphism:
	\begin{equation}\label{equa:tensor}
		\scrC(B\otimes^{\rmL}_{\bfH} M)\xrightarrow{\sim} \scrC(B)\otimes^{\rmL}_{\scrH}\scrC(M).
	\end{equation}
	This improves \ref{prop:HCTor-ii}
	Given any $\scrM\in \scrH\mof_{\fl}$, we have $\scrM\cong \scrC(M)$ for some $M\in \rmO(\bfH)$ by~\autoref{theo:DequivO}. Hence,~\eqref{equa:tensor} yields $\scrC(B)\otimes^{\rmL}_{\scrH}\scrM \cong \scrC(B\otimes^{\rmL}_{\bfH} M)$. By~\autoref{prop:ExtTorO}, $B\otimes^{\rmL}_{\bfH} M$ lies in $\Db_{\rmO}(\bfH)$, so $\scrC(B\otimes^{\rmL}_{\bfH} M)$ is in $\Db_{\fl}(\scrH)$ by~\autoref{theo:DequivO} again. Therefore, the functor $\scrC(B)\otimes^{\rmL}_{\scrH}\relbar$ preserves $\Db_{\fl}(\scrH)$. \par
	Given any $(c, \lambda)\in \frakP\times \haffd$, the projectivity of $\scrH 1_{c, \lambda}$ yields $\scrC(B)\otimes^{\rmL}_{\scrH}\scrH 1_{c, \lambda}  \cong \scrC(B)1_{c, \lambda}$. Let $\left\{ b_i \right\}_{i\in I}$ be a finite generating set of $B$ as left $\bfH$-module. Then, there exists a finite subset $\Sigma\subset\frakP\times\haffd$ such that $b_i 1_{c, \lambda}\in \scrC(B)1_{c, \lambda}$ is decomposed as $b_i 1_{c, \lambda} = \sum_{(c', \lambda')\in \Sigma} 1_{c',\lambda'}b_i 1_{c,\lambda}$ for each $i\in I$. Then, the family $\left\{ 1_{c',\lambda'}b_i 1_{c,\lambda} \right\}_{i\in I,\; (c',\lambda')\in \Sigma}$ yields the following surjective morphism of $\scrH$-modules:
	\[
		\bigoplus_{i\in I}\bigoplus_{(c',\lambda')\in \Sigma} \scrH 1_{c',\lambda'} \to \scrC(B) 1_{c, \lambda},\quad 
	\]
	Since the left-hand side is compact, so is right-hand side. As the family $\left\{ \scrH 1_{c, \lambda} \right\}_{(c,\lambda)\in \frakP\times\haffd}$ generates $\Db_{\perf}(\scrH)$ as thick subcategory, we see that the functor $\scrC(B)\otimes^{\rmL}_{\scrH}\relbar$ preserves $\Db_{\perf}(\scrH)$. This proves \ref{prop:HCTor-i}.
\end{proof}
	\begin{prop}\label{prop:HCExt}
Let $B\in\HC(\bfH)$. Then, the following statements hold:
	\begin{enumerate}
		\item\label{prop:HCExt-i}
			The functor $\Rhom_{\scrH}(\scrC B, \relbar)$  preserves the subcategories $\rmD^b_{\perf}(\scrH)$ and $\rmD^b_{\fl}(\scrH)$.
		\item\label{prop:HCExt-ii}
			There is a natural isomorphism
			\[
				\scrC(\RHom_{\bfH}(B, N))\cong \Rhom_{\scrH}(\scrC(B), \scrC(N))
			\]
			for $N\in \Db_{\perf}(\bfH)$. 
	\end{enumerate}
\end{prop}
\begin{proof}
	We prove \ref{prop:HCExt-ii}. The morphism is constructed as follows: for $N\in \bfH\Mod$, we have
	\[
		\Hom_{\bfH}(B, N)\to \Hom_{\scrH}(\scrC (B), \scrC(N))
	\]
	by the functoriality of $\scrC$; it induces
	\[
		\scrC\Hom_{\bfH}(B, N) \to \bigoplus_{c, \lambda}1_{c, \lambda}\Hom_{\scrH}(\scrC (B), \scrC(N)) = \hom_{\scrH}(\scrC (B), \scrC(N));
	\]
	the natural morphism
	\begin{equation}\label{equa:CRHom}
		\scrC\RHom_{\bfH}(B, N)\to \Rhom_{\scrH}(\scrC (B), \scrC(N))
	\end{equation}
	is obtained by passing to the derived functors.

	Since $\RHom_{\bfH}(B, \bfH)\in \Db_{\HC}(\bfH)$ by~\autoref{prop:extHC}, the spectral completion $\scrC(\RHom_{\bfH}(B, \bfH))$ is a complex of $\scrH$-bimodules by~\autoref{prop:Hbimod}. Consider first the case $N\in \rmO(\bfH)$. We have for $(c, \lambda)\in \frakP\times\haffd$ 
	\begin{align*}
		&1_{c,\lambda}\scrC(\RHom_{\bfH}(B, N))  \cong \Llim_{l}\RHom_{\bfS}\left(\bfS / \frakm_{c,\lambda}^l, \RHom_{\bfH}(B, N)\right) &\text{(\autoref{lemm:limS})}\\
		& \cong \Llim_{l}\RHom_{\bfH}\left(B\otimes^{\rmL}_{\bfS}(\bfS / \frakm_{c,\lambda}^l), N\right)& \text{(adjunction)}\\
		& \cong \Llim_{l}\RHom_{\scrH}\left(\scrC(B\otimes^{\rmL}_{\bfS}(\bfS / \frakm_{c,\lambda}^l)), \scrC(N)\right) & \text{(\autoref{theo:DequivO})} \\
		& \cong \Llim_{l}\RHom_{\bfS^{\wedge}_{c, \lambda}}\left(\bfS^{\wedge}_{c, \lambda} / \frakm_{c,\lambda}^l, \RHom_{\scrH}(\scrC(B)1_{c, \lambda},\scrC(N))\right) & \text{(adjunction)} \\
		& \cong \RHom_{\scrH}(\scrC(B)1_{c, \lambda}, \scrC(N)) & \text{(\autoref{lemm:limS})}.
	\end{align*}
	Thus, \eqref{equa:CRHom} is a quasi-isomorphism for $N\in \rmO(\bfH)$. Consider now the case $N = \bfH$. For $(c, \lambda),(c',\lambda')\in \frakP\times \haffd$, we have
	\begin{align*}
		&1_{c,\lambda}\scrC(\RHom_{\bfH}(B, \bfH))1_{c',\lambda'}\\
		&\cong \Rlim_{l}1_{c,\lambda}\scrC(\RHom_{\bfH}(B, \bfH))1_{c',\lambda'}\otimes^{\rmL}_{\bfS^{\wedge}_{c',\lambda'}} (\bfS/\frakm^l_{c',\lambda'})  & \text{(\autoref{lemm:limS})}\\
		&\cong \Rlim_{l}1_{c,\lambda}\scrC(\RHom_{\bfH}(B, \bfH\otimes^{\rmL}_{\bfS} (\bfS/\frakm^l_{c',\lambda'})))\\
	&\cong \Rlim_{l}\RHom_{\scrH}(\scrC(B)1_{c, \lambda}, \scrH 1_{c',\lambda'}\otimes_{\bfS^{\wedge}_{c',\lambda'}} (\bfS/\frakm^l_{c',\lambda'}))& \text{(case $N\in \rmO(\bfH)$)}\\
&\cong \RHom_{\scrH}(\scrC(B)1_{c, \lambda}, \Rlim_{l}\scrH 1_{c',\lambda'}\otimes_{\bfS^{\wedge}_{c',\lambda'}} (\bfS/\frakm^l_{c',\lambda'})))\\
&\cong \RHom_{\scrH}(\scrC(B)1_{c, \lambda}, \scrH 1_{c',\lambda'} )) = \RHom_{\scrH}(\scrC(B)1_{c, \lambda}, \scrH))1_{c',\lambda'}& \text{(Mittag-Leffler)}
	\end{align*}
	Thus, \eqref{equa:CRHom} is a quasi-isomorphism for $N = \bfH$. Since $\bfH$ generates $\Db_{\perf}(\bfH)$ as thick subcategory, this proves~\ref{prop:HCExt-ii}. \par
	We prove \ref{prop:HCExt-i}. Given any $\scrN\in \scrH\mof_{\fl}$, we have $\scrN\cong \scrC(N)$ for some $N\in \rmO(\bfH)$ by~\autoref{theo:DequivO}. Hence,~\ref{prop:HCExt-ii} yields $\Rhom_{\scrH}(\scrC(B),\scrN) \cong \scrC(\RHom_{\bfH}(B,N))$. By~\autoref{prop:ExtTorO}, $\RHom_{\bfH}(B,N)$ lies in $\Db_{\rmO}(\bfH)$, so $\scrC(\RHom_{\bfH}(B,N))$ lies in $\Db_{\fl}(\scrH)$ by~\autoref{theo:DequivO} again. Therefore, the functor $\Rhom_{\scrH}(\scrC(B), \relbar)$ preserves $\Db_{\fl}(\scrH)$. \par
Let $(c, \lambda)\in \frakP\times \haffd$. By \ref{prop:HCExt-ii}, we have
\[
	\Rhom_{\scrH}(\scrC(B), \scrH 1_{c, \lambda}) \cong \scrC(\RHom_{\bfH}(B, \bfH)) 1_{c, \lambda}\cong \scrC(\RHom_{\bfH}(B, \bfH))\otimes^{\rmL}_{\scrH}\scrH 1_{c, \lambda}.
\]
Since $\RHom_{\bfH}(B, \bfH)\in \Db_{\HC}(\bfH)$ by~\autoref{prop:extHC}, we have $\scrC(\RHom_{\bfH}(B, \bfH))\otimes^{\rmL}_{\scrH}\scrH 1_{c, \lambda} \in \Db_{\perf}(\scrH)$ by \autoref{prop:HCTor}. Since the family $\left\{ \scrH 1_{c, \lambda} \right\}_{(c,\lambda)\in \frakP\times\haffd}$ generates $\Db_{\perf}(\scrH)$ as thick subcategory, it follows that functor $\Rhom_{\scrH}(\scrC(B), \relbar)$ preserves $\Db_{\perf}(\scrH)$. This proves \ref{prop:HCExt-i}.
\end{proof}
\begin{lemm}\label{lemm:limS}
	Let $A$ be a complete noetherian local ring with maximal ideal $\frakm$. Then the following natural morphisms for $K\in \rmD_{\perf}(A)$ and $L\in \Db_{\fl}(A)$ are quasi-isomorphisms:
	\[
		K\to \Rlim_k (A/\frakm^k\otimes^{\rmL}_A K),\quad \Llim_k \RHom(A / \frakm^k, L)\to L.
	\]
\end{lemm}
\begin{proof}
	The first morphism is a quasi-isomorphism: as it is so for the regular $A$-module, also is it for every perfect complex $K\in \rmD_{\perf}(A)$ by d\'evissage. As for the second morphism, Nakayama's lemma implies that $\varinjlim_k \Hom_A(A / \frakm^k, L)\to L$ is an isomorphism for $L\in A\mof_{\fl}$. Given $L\in A\mof_{\fl}$, let $L\hookrightarrow I$ be an injection with $I$ injective. We may find such $I$ whose finitely generated submodules are of finite length. Then, by the compactness of $A / \frakm^k$, an easy argument of d\'evissage shows that $\varinjlim_k \Ext^n_A(A / \frakm^k, L) = 0$ for $n \ge 1$.
\end{proof}

\section{Category of divided-difference calculus}\label{sec:A}
In this section, we introduce a small $\bbC$-linear category $\bfA$, which realises the trigonometric DAHA $\bfH$ as the endomorphism ring of an object.
\subsection{Chambers and galleries}\label{subsec:Ch}
Recall $\Phi$ from~\autoref{subsec:Raff}. Consider the affine space $\fraka:=\frakP_{\bbR} \times \haffd_{\bbR}$ and the family of functions:
\[
	\Psi = \left\{ \alpha - \bfc_{\alpha}\in \frakP^*_{\bbQ} \times \haff^*_{\bbQ}\;;\; \alpha\in \Raff\right\},
\]
where $\haff^*_{\bbQ}$ is the $\bbQ$-linear span of $\Phi$ and $\frakP_{\bbQ}^*$ is the $\bbQ$-liear span of $\left\{ \bfc_* \right\}_{*\in R/W}$. 
For $\mu\in \Psi\cup \Raff$, let $H_{\mu}\subset \fraka$ be its zero locus. Then, $\left\{ H_{\mu} \;;\;\mu\in\Psi\cup \Raff\right\}$ defines a hyperplane arrangement on $\fraka$. The hyperplanes $H_{\mu}$ are called \emph{walls}. Put
\[
	\Ch(\fraka) = \pi_0\left(\fraka \setminus \bigcup_{\mu\in \Psi\cup \Raff}H_{\mu}\right),
\]
the elements of which are called \emph{chambers}. For each chamber $C\in \Ch(\fraka)$, let 
\[
	\Raff_C^{\pm} = \left\{ \alpha\in \Raff\;;\; \pm \alpha(C)\subseteq \bbR_{>0} \right\},\quad \Psi_{C}^{\pm} = \left\{ \mu\in \Psi\;;\; \pm \mu(C)\subseteq \bbR_{>0} \right\}.
\] \par

For $C_1, C_2\in \Ch(\fraka)$, we write $C_1\lceil_{\mu} C_2$ for $\mu\in \Psi$ if $C_1\neq C_2$ and they share a face with support $H_{\mu}$ and $\mu(C_1) < 0$, $\mu(C_2) > 0$; we write $C_1\|_{\alpha} C_2$ for $\alpha\in \Raff$ if $C_1\neq C_2$ and they share a face with support $H_{\alpha}$ and $\alpha(C_1) < 0$, $\alpha(C_2) > 0$. We write $C_1 \sim C_2$ if $C_1$ and $C_2$ lie in the same connected component of $\fraka \setminus \bigcup_{\alpha\in \Raff} H_{\alpha}$. \par

A sequence of chambers $G = (C_0, \ldots, C_n)$ is called a \emph{gallery} (from $C_0$ to $C_n$) if, for $i = 0, \ldots, n-1$, the chambers $C_i$ and $C_{i+1}$ share a face and $C_i \neq C_{i+1}$; we denote $\ell(G) = n$ and call it the \emph{length} of $G$. The \emph{distance} $d(C, C')$ is defined to be the minimum of the length of galleries from $C$ to $C'$. A gallery $G = (C_0, \ldots, C_n)$ is called \emph{minimal} if $\ell(G) = d(C_0, C_n)$. Given chambers $C, C'\in \Ch(\fraka)$, we define the \emph{interval} between $C$ and $C'$ to be 
\[
	[C, C'] = \left\{ C''\in \Ch(\fraka)\;;\; d(C, C') = d(C, C'') + d(C'', C') \right\}. 
\]
\begin{lemm}\label{lemm:minimal}
	The following statements hold:
	\begin{enumerate}
		\item
			Let $C, C', C''\in \Ch(\fraka)$ be chambers. Then, $C'\in [C, C'']$ holds if and only if the following conditions hold:
			\[
				\Phi^+_{C}\cap \Phi^+_{C''}\subseteq \Phi^+_{C'}\subseteq \Phi^+_{C}\cup \Phi^+_{C''},\quad \Psi^+_{C}\cap \Psi^+_{C''}\subseteq \Psi^+_{C'}\subseteq \Psi^+_{C}\cup \Psi^+_{C''}.
			\]
		\item
			A gallery $G = (C_0, \cdots, C_n)$ is minimal if and only if for each triplet $0\le i_1 < i_2 < i_3\le n$, we have $C_{i_2}\in [C_{i_1}, C_{i_3}]$.
	\end{enumerate}
\end{lemm}
\begin{proof}
	The first assertion results immediately from the fact that for $C, C'\in \Ch(\fraka)$, the distance $d(C, C')$ is the number of walls separating $C$ and $C'$. The second can be easily derived from the first one. 
\end{proof}

 If $G = (C_0, \ldots, C_n)$ and $G' = (C_n, \ldots, C_{n+m})$ are galleries, their \emph{composite} is defined to be $GG'= (C_0, \ldots, C_{n+m})$. The gallery \emph{opposite to} a gallery $G = (C_0, \ldots, C_n)$ is defined to be $G^* = (C_n, \ldots, C_0)$; for $w\in \Waff$, the transport of $G$ by $w$ is defined to be $wG = (wC_0, \ldots, wC_n)$. \par

 \par

 \subsection{The category \texorpdfstring{$\bfA^o$}{Ao}}\label{subsec:Ao}

\subsubsection{}\label{subsubsec:Ao}
Recall $\bfS = \calO(\frakP \times \haffd)$. Let $\bfA^o$ be the $\bbC$-linear category whose objects are $\Ch(\fraka)$ and whose hom-spaces $\Hom_{\bfA^o}(C, C')$ are subspaces of $\End_\bbC(\bfS)$ for $C, C'\in \Ch(\fraka)$ to be defined below. \par
For each $C, C'\in \Ch(\fraka)$ such that $C'\|_{\alpha} C$ for $\alpha\in \Raff$, let $\tau^o_{C, C'} = \vartheta_{\alpha}\in \End_{\bbC}(\bfS)$, where $\vartheta_{\alpha}:f\mapsto \alpha^{-1}(f- \pre{s_{\alpha}}f)$ is the Demazure operator. 

For chambers $C, C'\in \Ch(\fraka)$ such that $C\sim C'$, put $\tau^o_{C, C'} = 1$. For each gallery $G = (C_0, \ldots, C_n)$, put $\tau^o_G = \tau^o_{C_{n-1}, C_n}\cdots\tau^o_{C_{0}, C_1}\in \End_{\bbC}(\bfS)$. 
For $C, C'\in \Ch(\fraka)$ the morphisms are given by 
	\[
		\Hom_{\bfA^o}(C, C') = \sum_{w\in \Waff}\sum_{G = (C, C_1, \ldots, C_{n-1}, w^{-1}C')} w\,\tau^o_{G}\;\bfS\subseteq \End_{\bbC}(\bfS).
	\]
	The composition of morphisms is the usual composition of linear maps. Let $1_{C, C'} = 1_{\bfS}\in \Hom_{\bfA^o}(C, C')$ denote the identity map.
\subsubsection{}\label{subsubsec:grad}
The \emph{canonical filtration} on the hom-spaces of $\bfA^o$ is defined to be
\[
	F^{\can}_n \Hom_{\bfA^o}(C, C') = \left\{ f\in \Hom_{\bfA^o}(C, C')\;;\; \forall m\in \bbZ,\; f(\bfS_{\le m})\subseteq \bfS_{\le m + n}\right\}
\]
for $C,C'\in \Ch(\fraka)$ and $n\in \bbZ$, where $\bfS_{\le m}$ stands for polynomial functions on $\frakP\times \haffd$ of order $\le m$.
The \emph{length filtration} is defined to be
\[
	F^{\lg}_n \Hom_{\bfA^o}(C, C') = \sum_{w\in \Waff}\sum_{\substack{G = (C, C_1, \ldots, C_{n-1}, w^{-1}C') \\ \ell(G) \le n}} w\,\tau^o_{G}\;\bfS.
\]
Both filtrations are compatible with composition of morphisms. Note that each filtered piece $F^{\lg}_n\Hom_{\bfA^o}(C, C')$ is an $\bfS$-sub-bimodule of $\Hom_{\bfA^o}(C, C')$.

\subsubsection{}
	Let $C, C'\in \Ch(\fraka)$. For $w\in \Waff$, choose a minimal gallery $G_w$ from $C$ to $w^{-1}C'$ and let 
	\[
		\tau^o_{C, C', w} = w \,\tau^o_{G_w}\in F_{ d(C, w^{-1}C')}\Hom_{\bfA^o}(C, C').
	\]
	Let $\ba\tau^o_{C, C', w}$ denote the image of $\tau^o_{C, C', w}$ in the associated graded quotient 
	\[
		\gr^{\lg}_{ d(C, w^{-1}C')}\Hom_{\bfA^o}(C, C') = F^{\lg}_{ d(C, w^{-1}C')}\Hom_{\bfA^o}(C, C') / F^{\lg}_{ d(C, w^{-1}C') - 1}\Hom_{\bfA^o}(C, C').
	\]
	\begin{lemm}\label{lemm:basisAo}
		Given $C, C'\in \Ch(\fraka)$, the following statements hold:
		\begin{enumerate}
			\item
				if $G = (C, \ldots, w^{-1}C')$ is a non-minimal gallery, then $w\tau_G^o\in F_{\ell(G) - 1}\Hom_{\bfA^o}(C, C')$ holds.
			\item
				$\tau^o_{C, C', w}$ is independent of the choice of the minimal gallery $G_w$ for $w\in \Waff$;
			\item
				$\{ \ba\tau^o_{C, C', w}\}_{w\in \Waff}$ forms a basis for $\gr^{\lg}\Hom_{\bfA^o}(C, C')$ as free left and right $\bfS$-module; 
		\end{enumerate}
	\end{lemm}
	\begin{proof}
		According to~\cite[11.1.2]{kumar02}, the Demazure operators $\vartheta_{\alpha}\in \End(\bfS)$ are square-zero $\vartheta_{\alpha}^2 = 0$ and satisfy the braid relations for $\alpha,\beta\in \Deltaaff$, $\alpha\neq \beta$:
		\[
			\underbrace{\vartheta_{\alpha}\vartheta_{\beta}\vartheta_{\alpha}\cdots}_{m_{\alpha, \beta}} = \underbrace{\vartheta_{\beta}\vartheta_{\alpha}\vartheta_{\beta}\cdots}_{m_{\alpha, \beta}},\quad m_{\alpha, \beta} = \ord(s_{\alpha}s_{\beta}).
		\]
		For $w\in \Waff$, choose a reduced decomposition $w = s_{\beta_1} \cdots s_{\beta_k}$ for $w$ and put $\vartheta_{w} = \vartheta_{\beta_1}\cdots \vartheta_{\beta_k}$, which is, in view of the braid relations, independent of the choice of the reduced decomposition. Moreover, the family $\left\{ \vartheta_w \right\}_{w\in \Waff}$ is free in $\End_{\bbC}(\bfS)$ both as left and as right $\bfS$-modules, see~\cite[11.1.3]{kumar02}. \par
		Let $\kappa_0\in \Ch(\fraka)$ be the chamber which contains $\nu_0\times \left\{ 0 \right\}$, where $\nu_0\subseteq \haffd_{\bbR}$ is the fundamental alcove associated with the base $\Deltaaff$ (see \autoref{sssec:coxeter}). We deduce easily that $\tau^o_{C, C', w} = y'\vartheta_{y^{-1}w^{-1}y'}y^{-1}$, where $y, y'\in \Waff$ are such that $y\kappa_0\sim C$ and $y'\kappa_0 \sim C'$ hold. Now, given $w\in \Waff$ and a galery $G = (C, \ldots, w^{-1}C')$, we have either $w\tau^o_G = \tau^o_{C, C', w}$ or $\tau^o_G = 0$. The assertions follow easily from this.
	\end{proof}

	\subsection{The category \texorpdfstring{$\bfA$}{A}}\label{subsec:A}
\subsubsection{}
For $C, C'\in \Ch(\fraka)$, we define two invariants in $\bfS$:
\[
	\frakd(C, C') = \prod_{\mu\in \Psi_{ C'}^+\cap \Psi_{ C}^-}\mu,\quad \frake(C, C') = \prod_{\alpha\in \Raff_{C}^+\cap \Raff_{C'}^-}\alpha.
\]
For $C, C'\in \Ch(\fraka)$ such that $C\sim C'$, set 
\[
	\tau_{C, C'} = \tau^o_{C, C'}\;\frakd(C, C')\in \Hom_{\bfA^o}(C, C'),\quad f\mapsto f\;\frakd(C, C').
\]
For each $C, C'\in \Ch(\fraka)$ such that $C'\|_{\alpha} C$ for $\alpha\in \Raff$, set 
\[
	\tau_{C, C'} = \tau^o_{C,C'}\in \Hom_{\bfA^o}(C, C') ,\quad f\mapsto \vartheta_{\alpha}(f) = \alpha^{-1}(f- \pre{s_{\alpha}}f).
\]
Moreover, for $C\in \Ch(\fraka)$ and $w\in \Waff$, let $w_C\in \Hom_{\bfA^o}(C, wC)$ be given by $w_C: f\mapsto w(f)$. \par
For every gallery $G = (C_0, \ldots, C_n)$, we put $\tau_{G} = \tau_{C_{n-1}, C_n}\cdots \tau_{C_0, C_1}\in \Hom_{\bfA^o}(C_0, C_n)$. Similarly, we put $\frakd(G) = \prod_{0\le i\le n-1}\frakd(C_i, C_{i+1})$ and $\frake(G) = \prod_{0\le i\le n-1}\frake(C_i, C_{i+1})$. 

\subsubsection{}\label{subsubsec:A}
Let $\bfA\subset \bfA^o$ be the subcategory with the same objects whose morphisms are given by
\begin{equation}\label{equa:HomA}
	\Hom_{\bfA}(C, C') = \sum_{w\in \Waff}\sum_{G = (C, C_1, \ldots, C_{n-1}, w^{-1}C')} w\; \tau_{G}\;\bfS.
\end{equation}
For $n\in \bbZ$, let $F^{\lg}_{n}\Hom_{\bfA}(C, C')\subset \Hom_{\bfA}(C, C')$ be the sum in~\eqref{equa:HomA} with the second summation taken over galleries $G$ with $\ell(G) \le n$. \par
The canonical filtration $F^{\can}_\bullet\Hom_{\bfA}(C, C')$ for $C, C'\in \Ch(\fraka)$ is defined to be the induced filtration from $F^{\can}_{\bullet}\Hom_{\bfA^o}(C, C')$.

\subsection{Basis theorem}
Let $C, C'\in \Ch(\fraka)$. For $w\in \Waff$, choose a minimal gallery $G_w$ from $C$ to $w^{-1}C'$ and let 
\[
	\tau_{C, C', w} = w\,\tau_{G_w}\in F^{\lg}_{ d(C, w^{-1}C')}\Hom_{\bfA}(C, C').
\]
Let $\ba\tau_{C, C', w}$ denote its image in the quotient:
\[
	\gr^{\lg}_{d(C, w^{-1}C')}\Hom_{\bfA}(C, C') = F^{\lg}_{ d(C, w^{-1}C')}\Hom_{\bfA}(C, C') / F^{\lg}_{ d(C, w^{-1}C')-1}\Hom_{\bfA}(C, C'). 
\]

\begin{theo}\label{theo:basis}
	For $C, C'\in \Ch(\fraka)$ and $w\in \Waff$, the element $\ba\tau_{C, C', w}$ is independent of the choice of the minimal gallery $G_w$; moreover, $\{ \ba \tau_{C, C', w}\}_{w\in \Waff}$ forms a basis for $\gr^{\lg}\Hom_{\bfA}(C, C')$ as left and right free $\bfS$-module. 
\end{theo}
\begin{proof}
	We will only prove the right freeness, the proof for left freeness being similar. We prove by induction on $n\in \bbN$ the following statements for all chambers $C, C'\in \Ch(\fraka)$:
	\begin{enumerate}
		\item
			if $G = (C, C_1, \ldots, C_{n-1}, C')$ is a gallery, then 
			\[
				\tau_G - \tau^o_G\frakd(G) \in F^{\lg}_{n - 1}\Hom_{\bfA^o}(C, C'),\quad \tau_G\frake(G) - 1_{C, C'}\frakd(G) \in F^{\lg}_{n - 1}\Hom_{\bfA}(C, C');
			\]
		\item
			if $G = (C, C_1, \ldots, C_{n-1}, C')$ is a non-minimal gallery, then $\tau_G\in F^{\lg}_{ n-1} \Hom_{\bfA}(C, C')$;
		\item
			if $G = (C, C_1, \ldots, C_{n-1}, C')$ and $G' = (C, C'_1, \ldots, C'_{n-1}, C')$ are minimal galleries, then $\tau_G - \tau_{G'}\in F^{\lg}_{ n-1} \Hom_{\bfA}(C, C')$;
		\item
			the set $\{ \ba\tau_{C, C', w} \}_{\substack{w\in \Waff \\ d(C, w^{-1}C') = n}}$ forms a basis for 
			\[
				\gr^{\lg}_{n}\Hom_{\bfA}(C, C') = F^{\lg}_{ n} \Hom_\bfA(C, C') / F^{\lg}_{ n-1} \Hom_\bfA(C, C')
			\]
			as free right $\bfS$-module;
		\item
			the quotient
			\[
				F^{\lg}_{ n} \Hom_{\bfA^o}(C, C') / F^{\lg}_{ n} \Hom_\bfA(C, C')
			\]
			is right $\alpha$-torsion-free for every $\alpha\in \Raff$.
	\end{enumerate}
	For $n = 0$, all these five statements are trivial. Let $n \ge 1$ and suppose that the statements are proven for $0, \ldots, n-1$. We prove them for $n$:
	\begin{enumerate}
		\item
			Let $G^- = (C_1, \ldots, C_{n-1}, C')$. We have
			\[
				\tau_G - \tau^o_G\frakd(G) = (\tau_{G^-} - \tau^o_{G^-}\frakd(G^-))\tau_{C, C_1} + \tau^o_{G^-}(\frakd(G^-)\tau_{C, C_1} - \tau^o_{C, C_1}\frakd(G)).
			\]
			The term $\tau_{G^-} - \tau^o_{G^-}\frakd(G^-)$ lies in $F^{\lg}_{ n-2}\Hom_{\bfA^o}(C_1, C')$ by induction hypothesis. If $C\|_{\alpha} C_1$ for $\alpha\in \Raff$, then $\tau_{C, C_1} = \tau^o_{C, C_1} = \vartheta_{\alpha}$ is the Demazure operator, and thus
			\[
				\frakd(G^-)\tau_{C, C_1} = \tau^o_{C, C_1}\frakd(G^-) + \vartheta_{\alpha}\left(\frakd(G^-)  \right)s_{\alpha}.
			\]
			Since $\frakd(G) = \frakd(G^-)$ and $\vartheta_{\alpha}\left(\frakd(G^-)  \right)s_{\alpha}\in F^{\lg}_{ 0}\Hom_{\bfA^o}(C, C')$, the second term $\tau^o_{G^-}(\frakd(G^-)\tau_{C, C_1} - \tau^o_{C, C_1}\frakd(G))$ lies in $F^{\lg}_{ n-1}\Hom_{\bfA^o}(C, C')$. If $C\lceil C_1$ or $C_1\lceil C$, then $\frakd(G^-)\tau_{C, C_1} = \frakd(G^-)\frakd(C, C_1) = \frakd(G) = \tau^o_{C, C_1}\frakd(G)$, so $\tau^o_{G^-}(\frakd(G^-)\tau_{C, C_1} - \tau^o_{C, C_1}\frakd(G)) = 0$. It follows that $\tau_G - \tau^o_G\frakd(G)\in F^{\lg}_{n-1}\Hom_{\bfA^o}(C, C')$ in all cases. The other statement can be proven similarly. 
		\item
			By (i), we have
			\[
				\tau_G - \tau^o_G\frakd(G) \in F^{\lg}_{ n - 1}\Hom_{\bfA^o}(C, C').
			\]
			When $G$ is non-minimal, $\tau^o_G\in F^{\lg}_{ n - 1}\Hom_{\bfA^o}(C, C')$ holds by \autorefitem{lemm:basisAo}{i}; hence $\tau_G\in F^{\lg}_{ n - 1}\Hom_{\bfA^o}(C, C')$. 
			On the other hand, 
			\[
				\tau_G\frake(G)\in F^{\lg}_{ n - 1}\Hom_{\bfA}(C, C').
			\]
			Let $\ubar \tau_G$ be the image of $\tau_G$ in $F^{\lg}_{ n-1}\Hom_{\bfA^o}(C, C') / F^{\lg}_{ n-1}\Hom_{\bfA}(C, C')$. It follows that $\ubar\tau_G \frake(G) = 0$. By the induction hypothesis, the space 
			\[
				F^{\lg}_{ n-1}\Hom_{\bfA^o}(C, C') / F^{\lg}_{ n-1}\Hom_{\bfA}(C, C')
			\]
			is $\alpha$-torsion-free for each $\alpha\in \Raff$ and therefore $\frake(G)$-torsion-free; hence $\ubar \tau_G = 0$ and subsequently $\tau_G\in F^{\lg}_{ n-1}\Hom_{\bfA}(C, C')$ holds.
		\item
			By~(i), we have $\tau_G - \tau^o_G\frakd(G)\in F^{\lg}_{ n-1}\Hom_{\bfA^o}(C, C')$ and $\tau_{G'} - \tau^o_{G'}\frakd(G')\in F^{\lg}_{ n-1}\Hom_{\bfA^o}(C, C')$. Since $G$ and $G'$ are minimal, we have $\frakd(G) = \frakd(G') = \frakd(C, C')$. On the other hand, the minimality also implies $\tau^o_G = \tau^o_{G'}$ by \autorefitem{lemm:basisAo}{ii}. Therefore
			\[
				\tau_G - \tau_{G'} = (\tau_G - \frakd(G)\tau^o_G) - (\tau_{G'} - \frakd(G')\tau^o_{G'})\in F^{\lg}_{ n-1}\Hom_{\bfA^o}(C, C')
			\]
			holds. On the other hand, we have $\frake(G) = \frake(G') = \frake(C, C')$; therefore, (i) yields:
			\[
				(\tau_G - \tau_{G'})\frake(C, C') = \tau_G\frake(G) - \tau_{G'}\frake(G')\in  F^{\lg}_{n-1}\Hom_{\bfA}(C, C');
			\]
			hence the image of $\tau_G - \tau_{G'}$ in $F^{\lg}_{n-1}\Hom_{\bfA^o}(C, C') / F^{\lg}_{n-1}\Hom_{\bfA}(C, C')$ is annihilated by right multiplication by $\frake(C, C')$. The induction hypothesis implies that $F^{\lg}_{n-1}\Hom_{\bfA^o}(C, C') / F^{\lg}_{n-1}\Hom_{\bfA}(C, C')$ is right $\frake(C, C')$-torsion-free, so $\tau_G - \tau_{G'}\in F^{\lg}_{n-1}\Hom_{\bfA}(C, C')$ holds.
		\item
			By (ii) and (iii), the space $\gr^{\lg}_{n}\Hom_{\bfA}(C, C')$ is spanned by $\{ \ba\tau_{C, C', w} \}_{\substack{w\in \Waff \\ d(C, w^{-1}C') = n}}$ as right $\bfS$-module. This family is free --- indeed, suppose there is an $\bfS$-linear relation
			\[
				\sum_{w}\ba\tau_{C, C', w} f_w = 0,\quad f_w\in \bfS; 
			\]
			then, we have $\sum_{w}\tau_{C, C', w}f_w\in F^{\lg}_{ n-1}\Hom_{\bfA}(C, C')$; by (i), we have
			\[
				\sum_{w}\tau^{o}_{C, C', w}\frakd(C, w^{-1}C')f_w\in F^{\lg}_{ n-1}\Hom_{\bfA^o}(C, C');
			\]
			since the family $\{ \ba\tau^{o}_{C, C', w} \}_{\substack{w\in \Waff \\ d(C, w^{-1}C') = n}}$ is free over $\bfS$ by~\autorefitem{lemm:basisAo}{iii}, we see that $\frakd(C, w^{-1}C')f_w = 0$ for each $w$; hence $f_w = 0$. 
		\item
			Consider the following diagram:
			\[
				\begin{tikzcd}[column sep=10pt]
					0\arrow{r}& F^{\lg}_{ n-1}\Hom_{\bfA}(C, C') \arrow{d}{\phi'}\arrow{r}& F^{\lg}_{ n}\Hom_{\bfA}(C, C')\arrow{d}{\phi}\arrow{r} &  \gr^{\lg}_{n}\Hom_{\bfA}(C, C')\arrow{d}{\phi''}\arrow{r} & 0 \\
					0\arrow{r}& F^{\lg}_{ n-1}\Hom_{\bfA^o}(C, C')\arrow{r} &F^{\lg}_{ n}\Hom_{\bfA^o}(C, C')\arrow{r} & \gr^{\lg}_{n}\Hom_{\bfA^o}(C, C')\arrow{r} & 0.
				\end{tikzcd}
			\]
			The maps $\phi$ and $\phi'$ are the natural inclusions. By (i) and (iv), $\coker(\phi'')$ is isomorphic to
			\begin{align*}
				\bigoplus_{w} \ba\tau^{o}_{C, C', w}\bfS / \phi(\ba\tau_{C, C', w})\bfS \cong \bigoplus_{w} \bfS / \bfS\frakd(C, w^{-1}C')
			\end{align*}
			and $\phi''$ is injective. Since $\frakd(C, w^{-1}C')$ is not divisible by any element of $\Raff$, the cokernel $\coker(\phi'')$ is right $\alpha$-torsion-free for $\alpha\in \Raff$. The induction hypothesis implies that $\coker(\phi')$ is right $\alpha$-torsion-free; hence, $\coker(\phi)$ is also right $\alpha$-torsion-free by the snake lemma.
	\end{enumerate}
\end{proof}

\subsection{Translation symmetry}\label{subsec:symetrie}
The lattice $\frakP_{\bbZ}\subset \frakP$ acts on $\frakP_{\bbQ} \times \haff_{\bbQ}$ as follows: for $d\in \frakP_{\bbZ}$, let $t_d: \frakP_{\bbQ} \times \haff_{\bbQ}\to \frakP_{\bbQ}\times\haff_{\bbQ}$  be the translation $(c, z)\mapsto (\delta(z)d + c, z)$. This action extends $\bbR$-linearly to $\frakP_{\bbR}\times\haff_{\bbR}$ and restricts to $\fraka$.  The contragredient $\frakP_{\bbZ}$-action on the dual $\frakP^*_{\bbQ} \times \haff^*_{\bbQ}$ preserves the subsets $\Raff$ and $\Psi$. Consequently, the $\frakP_{\bbZ}$-action on $\fraka$ induces a $\frakP_{\bbZ}$-action on $\Ch(\fraka)$. \par

Similarly, the scalar extension of contragredient action induces a $\frakP_{\bbZ}$-action on $\bfS \cong \Sym(\frakP^* \times \haff^*) / (\delta - 1)$ by ring automorphisms, denoted by $t_d^*:\bfS\to \bfS$ for $d\in \frakP_{\bbZ}$. \par

We define a $\frakP_{\bbZ}$-action on $\bfA$ as follows: define $t_d: \bfA^o\to \bfA^o$ by setting 
\begin{align*}
	\Ch(\fraka)\ni C&\mapsto t_d(C)\in \Ch(\fraka)\\
	\Hom_{\bfA^o}(C, C')\ni a&\mapsto (t_d)_*(a):= (t_d^*)^{-1}\circ a\circ t_d^*\in \Hom_{\bfA^o}(t_d(C), t_d(C')).
\end{align*}
The automorphism $t_d$ preserves the subcategory $\bfA\subset \bfA^o$: indeed, we have $(t_d)_*\tau_G = \tau_{t_d(G)}$ for every gallery $G$,\; $(t_d)_*w = w$ in $\End_{\bbC}(\bfS)$ for every $w\in \Waff$ and $(t_d)_*(f) = (t_d^*)^{-1}(f)$ for every $f\in \bfS$. 
\subsection{Embedding of \texorpdfstring{$\bfH$}{H}}
Let $\kappa_0\in \Ch(\fraka)$ be the chamber which contains $\nu_0\times \left\{ 0 \right\}$, where $\nu_0\subseteq\haffd_{\bbR}$ is the fundamental alcove. For $d\in \frakP_{\bbZ}$, let $\kappa_d = t_d(\kappa_0) = \kappa_0 + (d, 0)\in \Ch(\fraka)$. 
\begin{theo}\label{theo:iota}
	For each $d\in \frakP_{\bbZ}$, there is an isomorphism of rings
	\begin{align*}
		\iota_d:\bfH&\xrightarrow{\sim} \End_{\bfA}(\kappa_d) \\
		f &\mapsto (t_d)_*f \quad \forall f\in \bfS\\
		1 - s_{\alpha} &\mapsto (\alpha - (\bfc_{\alpha} - d_{\alpha}))\vartheta_{\alpha} \quad \forall \alpha\in \Deltaaff,
	\end{align*}
	where $(t_d)_*f\in \bfS$ is given by $((t_d)_*f)(c, z) = f(c - d, z)$ for $(c, z)\in \frakP \times \haffd$. 
\end{theo}
\begin{proof}
	Let us consider first the case $d = 0$. It is well-known that the spherical polynomial representation $\bfH\otimes_{\bbC \Waff}\triv$ induced from the trivial character $\triv:\Waff\to \bbC^{\times}$ of the affine Weyl group is faithful, the underlying vector space is isomorphic to $\bfS$, and the $\bfH$-action on it is given by the homomorphism $\iota_0:\bfH\to \End_{\bbC}(\bfS)$ defined in the statement. In particular, $\iota_0$ is well-defined and injective. We shall prove the following:
	\begin{enumerate}
		\item
			$(\alpha - \bfc_{\alpha})\vartheta_\alpha$ lies in $\End_{\bfA}(\kappa_0)$ for each $\alpha\in \Deltaaff$;
		\item
			$\{(\alpha - \bfc_{\alpha})\vartheta_\alpha\}_{\alpha\in \Deltaaff}\cup \bfS$ generates $\End_{\bfA}(\kappa_0)$.
	\end{enumerate}
	Let $\alpha\in \Deltaaff$. We have $\frakd(\kappa_0, s_{\alpha}\kappa_0) = -\alpha - \bfc_{\alpha}$ and $\frakd(s_{\alpha}\kappa_0, \kappa_0) = \alpha - \bfc_{\alpha}$.  The intersection $\ba{\nu_0}\cap H_{\alpha}\subseteq\haff_{\bbR}$ is the closure of a face of the fundamental alcove $\nu_0$; let $z\in \ba{\nu_0}\cap H_{\alpha}$ be a point in its relative interior. Let $\epsilon = (\epsilon_*)\in \frakP_{\bbR}$ be such that $\epsilon_{\beta} > 0$ is small enough for every $\beta$, the point $(\epsilon, z)\in \fraka$ lies in a connected component of $\fraka \setminus \bigcup_{\psi\in \Psi}H_{\psi}$, denoted by $Z$, which is determined as follows:
	\[
		Z = \bigcap_{\beta\in \Phi^+\setminus\left\{ \alpha \right\}}\left\{ \beta - \bfc_\beta > 0 \right\}\cap \bigcap_{\beta\in \Phi^-\cup\left\{ \alpha \right\}}\left\{ \beta - \bfc_\beta < 0 \right\}.
	\]
	Here, $\Phi^{\pm}\subset \Phi$ is the set of positive/negative affine roots associated with the base $\Deltaaff$ (see \autoref{sssec:coxeter}).
	The complement $Z\setminus H_{\alpha}$ consists of two chambers, one of which is determined by the inequality $\alpha > 0$ and shares a face with $\kappa_0$ with support $H_{\alpha - \bfc_{\alpha}}$; denote this chamber by $C_{\alpha}\in \Ch(\fraka)$. \par

	Then, it follows that $C_{\alpha}\lceil_{\alpha - \bfc_{\alpha}}\kappa_0$, $s_{\alpha}(C_{\alpha})\|_{\alpha}C_\alpha$ and $s_{\alpha}(C_{\alpha})\lceil_{-\alpha - \bfc_{\alpha}}\kappa_d$ hold, and thus $G_{\alpha} = (\kappa_0, C_\alpha, s_{\alpha}C_{\alpha}, s_{\alpha}\kappa_0)$ is a minimal gallery from $\kappa_0$ to $s_{\alpha} \kappa_0$. We have $s_{\alpha}\tau_{G_{\alpha}} = s_{\alpha}(-\alpha - \bfc_{\alpha})\vartheta_{\alpha} = (\alpha - \bfc_{\alpha})\vartheta_{\alpha}$ and thus $(\alpha - \bfc_{\alpha})\vartheta_{\alpha}$ lies in $\End_{\bfA}(\kappa_0)$. \par
	Given any $w\in \Waff$, let $w = s_{\beta_1}\cdots s_{\beta_k}$ be a reduced decomposition, where $\beta_i\in \Deltaaff$ are simple roots. Then, the composite
	\[
		G_w := G_{\beta_1}(s_{\beta_1}G_{\beta_2})(s_{\beta_1}s_{\beta_2}G_{\beta_3}) \cdots (s_{\beta_1}\cdots s_{\beta_{k-1}}G_{\beta_k})
	\]
	is a minimal gallery from $\kappa_0$ to $w\kappa_0$: indeed, if we put $\gamma_i = s_{\beta_1}\cdots s_{\beta_{i-1}}\beta_i$, then we have
	\[
		\frakd(s_{\beta_1}\cdots s_{\beta_{i-1}}G_{\beta_i}) = s_{\beta_1}\cdots s_{\beta_{i-1}}\frakd(G_{\beta_i}) = \gamma_i - \bfc_{\gamma_i}
	\]
	and hence $\frakd(G_w) = \prod_{i = 1}^k(\gamma_i - \bfc_{\gamma_i}) = \frakd(\kappa_0, w\kappa_0)$ because $\Raff^+ \cap w \Raff^- = \left\{ \gamma_1, \cdots, \gamma_k \right\}$; similarly, we have $\frake(G_w)= \frake(\kappa_0, w\kappa_0)$ and the same holds for $G_w^{-1}$; it follows that the gallery $G_w$ is minimal by~\autoref{lemm:minimal}. Therefore, the element 
	\[
		w^{-1}\tau_{G_w} = s_{\beta_k} \cdots s_{\beta_1}\tau_{s_{\beta_1}\ldots s_{\beta_{k-1}}G_{\beta_k}}\cdots\tau_{s_{\beta_1}G_{\beta_2}}\tau_{G_{\beta_1}} = s_{\beta_k}\tau_{G_{\beta_k}}\cdots s_{\beta_2}\tau_{G_{\beta_2}}s_{\beta_1}\tau_{G_{\beta_1}}
	\]
	lies in the subalgebra of $\End_{\bbC}(\bfS)$ generated by $\{(\alpha - \bfc_{\alpha})\vartheta_\alpha = s_{\alpha}\tau_{G_{\alpha}}\}_{\alpha\in \Deltaaff}$. By~\autoref{theo:basis}, it follows that $\iota_0$ is surjective and thus an isomorphism. \par
	For $d\in \frakP_{\bbZ}$, we may write $\iota_d = (t_d)_*\circ\iota_0$, which is clearly an isomorphism.
\end{proof}
\subsection{Homogeneous version}\label{subsec:Adelta}
Just as $\bfH$ admits a homogeneous version $\bfH_{\delta}$ isomorphic to the Rees algebra of the canonical filtration $F^{\can}_{\bullet}\bfH$, we can define a homogeneous version $\bfA_{\delta}$ of $\bfA$ by replacing $\bfS$ with $\bfS_{\delta}$ in the constructions presented in \autoref{subsec:Ch}--\autoref{subsec:A}. The category $\bfA_{\delta}$ has for set of objects $\Ch(\fraka)$ and its hom-spaces are graded subspaces of $\gEnd_{\bbC}(\bfS_{\delta})$, where $\gEnd_{\bbC}(\bfS_{\delta})$ is the graded Hom-space, endowed with the grading induced from $\bfS_{\delta}$ (see \autoref{subsubsec:filtcan}). There are obvious analogues of \autoref{theo:basis} and \autoref{theo:iota} for $\bfA_{\delta}$ and $\bfH_{\delta}$ which can be proven using the same arguments. In particular, the hom-space $\Hom_{\bfA_{\delta}}(C, C')$ for $C, C'\in \Ch(\fraka)$ is graded-free over $\bbC[\delta]$, so it can be identified as the Rees spaces:
\[
	\Hom_{\bfA_{\delta}}(C, C') \cong \bigoplus_{n\in \bbZ} \delta^{n}F^{\can}_{n}\Hom_{\bfA}(C, C').
\]
For $d\in \frakP_{\bbZ}$, the translation $(t_d)_*:\bfS_{\delta} \to \bfS_{\delta}$ is given by $(t_{d})_*(\bfc_\alpha) = \bfc_\alpha - \delta d_{\alpha}$ for $\alpha\in \Raff$ and $(t_{d})_*(f) = f$ for $f\in \calO(\haff)$.

\section{Spectral completion of \texorpdfstring{$\bfA$}{A}}\label{sec:Aloc}
Fix $(c, \lambda)\in \frakP \times \haffd$. We define the following finite subsets
\[
	\Raff_{\lambda} = \left\{ \alpha\in \Raff\;;\; \alpha(\lambda) = 0 \right\},\quad \Psi_{c, \lambda} = \left\{ \mu\in \Psi\;;\; \mu(c, \lambda) = 0 \right\}. 
\]
Recall that $\Scom_{c, \lambda}$ is the completion of $\bfS = \calO(\frakP\times \haffd)$ at the defining ideal of the point $(c, \lambda)\in \frakP\times \haffd$.

\subsection{The category \texorpdfstring{$\scrA^{c, \lambda}$}{A lambda, c}}\label{subsec:Aloc}
We shall define a category $\scrA^{c, \lambda}$ analogous to $\bfA$ by repeating the constructions from~\autoref{sec:A} with the following replacements: $\Phi\mapsto \Phi_{\lambda}$, $\Psi\mapsto \Psi_{c, \lambda}$, $\Waff\mapsto \Waff_{\lambda}$ and $\bfS\mapsto \Scom_{c, \lambda}$. 
\subsubsection{}
The sets $\Phi_{\lambda}$ and $\Psi_{c, \lambda}$ define a finite hyperplane arrangement $\left\{ H_{\varphi} \right\}_{\varphi\in \Phi_{\lambda}\cup \Psi_{c, \lambda}}$ on $\fraka$. Chambers are defined to be the connected components of $\fraka \setminus \bigcup_{\varphi\in \Phi_{\lambda}\cup\Psi_{c, \lambda}}H_{\varphi}$. The set of chambers is denoted by $\Ch^{c, \lambda}(\fraka)$.  The notion of galleries is the same as in~\autoref{subsec:Ch}. \par

For each chamber $C\in \Ch^{c, \lambda}(\fraka)$, define 
\[
	\Raff_{\lambda, C}^{\pm} = \left\{ \alpha\in \Raff_{\lambda}\;;\; \pm \alpha(C)\subseteq \bbR_{>0} \right\},\quad \Psi_{c, \lambda, C}^{\pm} = \left\{ \mu\in \Psi_{c, \lambda}\;;\; \pm \mu(C)\subseteq \bbR_{>0} \right\}.
\] \par

\subsubsection{}
For each $C, C'\in \Ch^{c, \lambda}(\fraka)$, we define an invariant in $\Scom_{c,\lambda}$:
\[
	\frakd^{c, \lambda}(C, C') = \prod_{\mu\in \Psi_{c, \lambda, C'}^+\cap \Psi_{c, \lambda,  C}^-}\mu.
\]
For each $C, C'\in \Ch^{c, \lambda}(\fraka)$ such that $C$ and $C'$ lie in the same connected component of $\fraka \setminus \bigcup_{\alpha\in \Phi_{\lambda}}H_{\alpha}$, let 
\[
	\tau^{c, \lambda}_{C, C'}\in \End_{\bbC}(\Scom_{c, \lambda}),\quad f\mapsto f\;\frakd^{c, \lambda}(C, C').
\]
For each $C, C'\in \Ch^{c, \lambda}(\fraka)$ such that $C'\|_{\alpha} C$ for $\alpha\in \Raff_{\lambda}$, let 
\[
	\tau^{c, \lambda}_{C, C'}\in \End_{\bbC}(\Scom_{c, \lambda}),\quad f\mapsto \vartheta_{\alpha}(f) = \alpha^{-1}(f- s_{\alpha}(f)).
\]
For every gallery $G = (C_0, \cdots, C_n)$ in $\Ch^{c, \lambda}(\fraka)$, we put $\tau^{c, \lambda}_{G} = \tau^{c, \lambda}_{C_{n-1}, C_n}\cdots \tau^{c, \lambda}_{C_0, C_1}\in \End_{\bbC}(\Scom_{c, \lambda})$. 

\subsubsection{}\label{subsubsec:Acompl}
Let $\scrA^{c, \lambda}$ be the small category with set of objects $\Ch^{c, \lambda}(\fraka)$ and with morphisms defined by
\begin{equation}\label{equa:HomAloc}
	\Hom_{\scrA^{c, \lambda}}(C, C') = \sum_{w\in \Waff_{\lambda}}\sum_{G = (C, C_1, \ldots, C_{n-1}, w^{-1}C')} w\;\tau^{c, \lambda}_{G}\;\Scom_{c,\lambda} \subseteq \End_{\bbC}(\Scom_{c,\lambda})
\end{equation}
For $n\in \bbZ$, let $F^{\lg}_{n}\Hom_{\scrA^{c,\lambda}}(C, C')\subset \Hom_{\scrA^{c,\lambda}}(C, C')$ be the sum in~\eqref{equa:HomAloc} with the second summation taken over galleries $G$ with $\ell(G) \le n$. 
\subsubsection{}
Let $C, C'\in \Ch^{c, \lambda}(\fraka)$. For $w\in \Waff_{\lambda}$, choose a minimal gallery $G_w$ from $C$ to $w^{-1}C'$ and let 
\[
	\tau^{c, \lambda}_{C, C', w} = w\;\tau^{c, \lambda}_{G_w}\in F^{\lg}_{n}\Hom_{\scrA^{c,\lambda}}(C, C'),\quad n = d(C, w^{-1}C');
\]
moreover, let $\ba\tau^{c, \lambda}_{C, C', w}$ denote its image in the quotient:
\[
	\gr^{\lg}_{n}\Hom_{\scrA^{c,\lambda}}(C, C') = F^{\lg}_{n}\Hom_{\scrA^{c,\lambda}}(C, C') / F^{\lg}_{n -1}\Hom_{\scrA^{c,\lambda}}(C, C').
\]

\begin{theo}\label{theo:basisloc}
	For $C, C'\in \Ch^{c, \lambda}(\fraka)$ and $w\in \Waff_{\lambda}$, the element $\ba\tau^{c, \lambda}_{C, C', w}$ is independent of the choice of the minimal gallery $G_w$; moreover, $\{ \ba \tau^{c,\lambda}_{C, C', w}\}_{w\in \Waff_{\lambda}}$ forms a basis for $\gr^F\Hom_{\scrA^{c, \lambda}}(C, C')$ as left and right free $\Scom_{c,\lambda}$-module. 
\end{theo}
\begin{proof}
	The proof of~\autoref{theo:basis} applies \textit{mutatis mutandis}.
\end{proof}

\subsubsection{}\label{subsubsec:transeq}
The translation action of $d\in \frakP_{\bbZ}$ induces $t_d^*:\frakP^*_{\bbQ} \times \haff^*_{\bbQ}\to \frakP^*_{\bbQ} \times \haff^*_{\bbQ}$, which sends $\Psi_{c+d, \lambda}$ to $\Psi_{c, \lambda}$ and fixes $\Phi_{\lambda}$ pointwise. Hence, it induces a map $t_d: \Ch^{c, \lambda}(\fraka)\to \Ch^{c+d, \lambda}(\fraka)$ and an equivalence of categories $t_d: \scrA^{c, \lambda}\xrightarrow{\sim} \scrA^{c+d, \lambda}$.

\subsection{Comparison}
We shall define the spectral completion $\scrC^{c, \lambda}\bfA$ of the category $\bfA$ and compare it with $\scrA^{c,\lambda}$.
\subsubsection{}
Let $\scrS^{c, \lambda} = \bigoplus_{\lambda'\in [\lambda]}\Scom_{c, \lambda'}$, which is a non-unital ring. For $C, C'\in \Ch(\fraka)$, define the spectral completion $\scrC^{c, \lambda}\Hom_{\bfA}(C, C') = \scrS^{c, \lambda}\otimes_{\bfS}\Hom_{\bfA}(C, C')$. \par

The basis theorem (\autoref{theo:basis}) implies that $\Hom_{\bfA}(C, C')$ is the union of its $\bfS$-sub-bimodules which are finite as left and right $\bfS$-modules. Therefore, \autoref{lemm:BSSB} applies to those sub-bimodules and yields an isomorphism:
\begin{equation}\label{equa:CHomA}
	\scrS^{c, \lambda}\otimes_{\bfS}\Hom_{\bfA}(C, C')\xrightarrow{\eta}\Hom_{\bfA}(C, C')\otimes_{\bfS}\scrS^{c, \lambda};
\end{equation}
therefore, for $C,C',C''\in \Ch(\fraka)$, we can define a composition map:
\[
	\circ:\scrC^{c, \lambda}\Hom_{\bfA}(C', C'')\times \scrC^{c, \lambda}\Hom_{\bfA}(C, C')\to \scrC^{c, \lambda}\Hom_{\bfA}(C, C'').
\]

\subsubsection{}
The formula in \autoref{subsubsec:scrH} defines an action of the Demazure operators $\left\{ \vartheta_{\alpha} \right\}_{\alpha\in \Raff}$ on $\scrS^{c, \lambda}$. Therefore, we obtain an embedding $\Hom_{\bfA^o}(C, C')\hookrightarrow \End_{\bbC}(\scrS^{c, \lambda})$. \par
It follows that the embedding $\Hom_{\bfA}(C, C')\subseteq\Hom_{\bfA^o}(C, C')\hookrightarrow \End_{\bbC}(\scrS^{c, \lambda})$ extends to $\scrC^{c, \lambda}\Hom_{\bfA}(C, C')\hookrightarrow \End_{\bbC}(\scrS^{c, \lambda})$, compatible with the compositions.

\subsubsection{}
For each $\lambda'\in [\lambda]$, let $1_{\lambda'}\in \Scom_{c, \lambda}$ denote the unit element viewed as an idempotent in $\scrS^{c, \lambda}$. Let $\scrC^{c, \lambda}\bfA$ be the $\bbC$-linear category with set of objects $\Ch(\fraka)\times \Waff$ and with morphisms 
\[
	\Hom_{\scrC^{c, \lambda}\bfA}((C, w), (C', w')) = 1_{w'\lambda}\scrC^{c, \lambda}\Hom_{\bfA}(C, C')1_{w\lambda}\subseteq \Hom_{\bbC}(\Scom_{c, w\lambda}, \Scom_{c, w'\lambda}).
\]
The composition of morphisms is the usual composition of linear maps. 

\begin{theo}\label{theo:comparaison}
	For each chamber $C\in \Ch(\fraka)$, let $\widetilde C\in \Ch^{c,\lambda}(\fraka)$ denote the unique chamber satisfying $C\subseteq \til C$. Then, there is an equivalence of categories defined by the following formula:
	\[
		\rho: \scrC^{c, \lambda}\bfA \xrightarrow{\sim}\scrA^{c, \lambda},\quad \rho(C, w) = \widetilde{w^{-1}C},\quad  \rho(\varphi)= w'^{-1} \varphi w
	\]
	for $\varphi\in \Hom_{\scrC^{c, \lambda}\bfA}((C, w), (C', w'))$, where $w:\Scom_{c, \lambda}\xrightarrow{\sim}\Scom_{c, w\lambda}$ and $w':\Scom_{c, \lambda}\xrightarrow{\sim}\Scom_{c, w'\lambda}$ are induced by completion from the $\Waff$-action on $\bfS$.
\end{theo}
\begin{proof}
	The morphisms of the category $\scrA^{c, \lambda}$ are generated by the following set:
	\begin{align*}
		&\{\tau^{c,\lambda}_{C, C'}\in \Hom_{\scrA^{c, \lambda}}(C, C')\;;\; C, C'\in \Ch^{c,\lambda}(\fraka) \text{ adjacent}\}\cdot\Scom_{c, \lambda} \\
		&\cup \{ y \in \Hom_{\scrA^{c, \lambda}}(C, yC)\;;\; y\in \Waff, C\in \Ch^{c,\lambda}(\fraka) \}\cdot\Scom_{c, \lambda}.
	\end{align*}
	The morphisms of $\scrC^{c, \lambda}\bfA$ are generated by the following set:
	\begin{align*}
		&\bigcup_{w,w'\in \Waff}\{1_{w'\lambda}\tau_{C, C'}1_{w\lambda}\in \Hom_{\scrC^{c,\lambda}\bfA}((C, w), (C', w'))\;;\; C,C'\in \Ch(\fraka) \text{ adjacent}\}\cdot\Scom_{c, w\lambda} \\
		&\cup \bigcup_{w\in \Waff}\{1_{yw\lambda} y 1_{w\lambda} \in \Hom_{\scrC^{c, \lambda} \bfA}((C, w), (yC, yw))\;;\; y\in \Waff, C\in \Ch(\fraka) \}\cdot\Scom_{c, w\lambda}.
	\end{align*}
	We show that $\rho$ sends the generators of morphisms of $\scrC^{c, \lambda}\bfA$ into $\scrA^{c, \lambda}$.  \par
	Let $C, C'\in \Ch(\fraka)$ be a pair of chambers which share a common wall. We show that for each $w,w'\in \Waff$, the following holds:
	\begin{equation}\label{equa:rhotau}
		\rho(1_{w'\lambda}\tau_{C, C'}1_{w\lambda})\in \Hom_{\scrA^{c, \lambda}}(D, D'),\quad \text{where } D = \widetilde{w^{-1}C}, D' = \widetilde{w'^{-1}C'}.
	\end{equation}
	Suppose that $C\|_{\alpha}C'$ for some $\alpha\in \Phi$. Then, $\tau_{C,C'} = \vartheta_{\alpha}$. If $\alpha(w\lambda) = 0$, then $w^{-1}\alpha \in \Phi_{\lambda}$ and $D\|_{w^{-1}\alpha}D'$ in $\Ch^{c,\lambda}(\fraka)$  hold; therefore, 
	\[
		\rho(1_{w'\lambda}\tau_{C, C'}1_{w\lambda}) = \begin{cases} \vartheta_{w^{-1}\alpha} = \tau^{c,\lambda}_{D, D'} & \text{if $w\lambda = w'\lambda$}, \\ 0 & \text{otherwise.}\end{cases}
	\]
	If $\alpha(w\lambda)\neq 0$, then $w^{-1}\alpha$ is invertible in $\Scom_{c,\lambda}$ and $\Scom_{c,s_{\alpha}\lambda}$; therefore, the Demazure operator induces
	\[
		\vartheta_{\alpha}: \Scom_{c,w\lambda}\to \Scom_{c,w\lambda}\oplus \Scom_{c,s_\alpha w\lambda},\quad f\mapsto \alpha^{-1}f - \alpha^{-1}\pre{s_\alpha}f.
	\]
	and thus 
	\[
	\rho(1_{w'\lambda}\tau_{C, C'}1_{w\lambda}) = \begin{cases} (w^{-1}\alpha)^{-1} & \text{if $w'\lambda = w\lambda$}, \\  (w^{-1}\alpha)^{-1}w'^{-1}s_{\alpha}w & \text{if $w'\lambda = s_{\alpha}w\lambda$},\\ 0 & \text{otherwise.}\end{cases}
	\]
	Note that in the second case above, we have $w'^{-1}s_{\alpha}w\in \Waff_{\lambda}$. It follows that~\eqref{equa:rhotau} holds in each of the cases. Similar arguments show~\eqref{equa:rhotau} also holds in the case $C\lceil_{\psi}C'$ or $C'\lceil_{\psi}C$ for some $\psi\in \Psi$ and that $\rho(1_{yw'\lambda}y1_{w\lambda})$ lies in $\Hom_{\scrA^{c, \lambda}}(D, D')$ for all $w, w', y\in \Waff$. Since the formula for $\rho$ respects the composition, it follows that $\rho$ is a well-defined functor. \par
	The same kind of arguments as above shows that $\rho$ is full. The faithfulness of $\rho$ is due to the fact that $\varphi\mapsto w'^{-1}\varphi w$ gives an isomorphism from $\Hom_{\bbC}(\Scom_{c,w\lambda}, \Scom_{c,w'\lambda})$ to $\End_{\bbC}(\Scom_{c, \lambda})$ and that $\Hom_{\scrC^{c,\lambda}\bfA}((C, w), (C', w'))\subseteq \Hom_{\bbC}(\Scom_{c,w\lambda}, \Scom_{c,w'\lambda})$ is a subspace. The essential surjectivity is straightforward.
\end{proof}

\subsection{Duality}\label{subsec:duality}
\subsubsection{}
Let $C, C'\in \Ch^{c, \lambda}(\fraka)$. We say that $C$ and $C'$ are \emph{antipodal} if the following conditions are met:
\[
	\Phi^+_{\lambda, C}= \Phi^-_{\lambda, C'},\quad \Psi^+_{c, \lambda,C}= \Psi^-_{c, \lambda, C'}.
\]
Suppose that $C, C'$ are antipodal. Since the stabiliser $\Waff_{\lambda}$ is generated by reflections about affine hyperplanes from $\Phi_{\lambda}$, there exists a unique element $\omega\in \Waff_{\lambda}$ such that $\omega\Phi^+_{\lambda, C} = \Phi^-_{\lambda, C}$ and $\omega^2 = 1$ holds. By~\autoref{theo:basisloc}, $d(C, C')$ is the maximal degree of $\gr^{\lg}\Hom_{\scrA^{c, \lambda}}(C, C')$ and the graded piece of degree $d(C, C')$ is $\ba\tau^{c, \lambda}_{C, C'}\Scom_{c, \lambda}$. \par

Define the trace map $\tr:\Hom_{\scrA^{c, \lambda}}(C, C')\to \Scom_{c, \lambda}$ to be the following composition:
\[
	\Hom_{\scrA^{c, \lambda}}(C, C')\to \gr^{\lg}_{d(C, C')}\Hom_{\scrA^{c, \lambda}}(C, C') = \ba\tau^{c, \lambda}_{C, C'}\Scom_{c, \lambda}\xrightarrow{\ba\tau^{c, \lambda}_{C, C'}f\mapsto f}  \Scom_{c, \lambda}.
\]
From the relation $\vartheta_\alpha a = \vartheta_\alpha(a) + s_{\alpha}(a)\vartheta_\alpha\in \End_{\bbC}(\Scom_{c, \lambda})$ satisfied by the Demazure operator for $\alpha\in \Phi^+_{\lambda,C}$ and $a\in \Scom_{c, \lambda}$, it is not hard to show that 
\[
	\tr(a\varphi) = \omega(a)\tr(\varphi),\quad \tr(\varphi a) = \tr(\varphi)a
\]
for $a\in \Scom_{c, \lambda}$ and $\varphi\in \Hom_{\scrA^{c, \lambda}}(C, C')$.
\begin{lemm}\label{lemm:pairing}
	Suppose that $C, C'\in \Ch^{c, \lambda}(\fraka)$ are antipodal chambers. Then, given any chamber $C''\in \Ch^{c, \lambda}(\fraka)$, the trace map induces a $(\Scom_{c, \lambda},\omega)$-sesquilinear perfect pairing:
	\[
		\Hom_{\scrA^{c, \lambda}}(C'', C')\times \Hom_{\scrA^{c, \lambda}}(C, C'')\to \Scom_{c, \lambda},\quad (\phi, \psi)\mapsto \tr(\phi\psi).
	\]
\end{lemm}
\begin{proof}
	Since $C$ and $C'$ are antipodal, it is easy to see that $[C, C'] = \Ch^{c,\lambda}(\fraka)$; namely, $d(C, C') = d(C, C'') + d(C'', C')$ holds for every $C''\in \Ch^{c,\lambda}(\fraka)$. In particular, for each $w\in \Waff$, there exists a minimal gallery from $C$ to $C'$ passing through $w^{-1}C''$. It follows from~\autoref{theo:basisloc} that
	\[
		\ba\tau_{C'', C', w^{-1}}\ba\tau_{C, C'', w} = \ba\tau_{w^{-1}C'', C'}\ba\tau_{C, w^{-1}C''} = \ba\tau_{C, C'}
	\]
	holds. \par
	Pick any liftings $\tau_{C'', C', w}$ and $\tau_{C, C'', w}$ for each $w\in \Waff$; then,~\autoref{theo:basisloc} implies that $\left\{ \tau_{C'', C', w} \right\}_{w\in \Waff}$ and $\left\{ \tau_{C, C'', w} \right\}_{w\in \Waff}$ are free bases over $\Scom_{c, \lambda}$. If we order the elements $w$ of $\Waff$ by the distance $d(C, w^{-1}C'')$, then the matrix $\left(\tr(\tau_{C'', C', y^{-1}}\tau_{C, C'', w})\right)_{w, y\in \Waff}$ is triangular with 1s on the diagonal. Hence, the pairing is perfect.
\end{proof}
\subsubsection{}\label{subsubsec:duality}
Recall the ring $\scrZ^{c, \lambda} = (\Scom_{c, \lambda})^{\Waff_\lambda}$ from \autoref{subsubsec:Z}. The inclusion $\scrZ^{c, \lambda}\hookrightarrow \Scom_{c, \lambda}$ induces an action of $\scrZ^{c, \lambda}$ on $\End_{\scrA^{c, \lambda}}(C, C')\subseteq \End_{\bbC}(\Scom_{c, \lambda})$ by multiplication for each $C, C'\in \Ch^{c,\lambda}(\fraka)$. This yields $\scrZ^{c, \lambda} \to \End(\id_{\scrA^{c, \lambda}})$. It is not hard to show that this is an isomorphism. \par
It is known that $\Scom_{c, \lambda}$ is a symmetric algebra over $\scrZ^{c, \lambda}$ with trace form $\tr_{\Scom_{c, \lambda}}:\Scom_{c, \lambda}\to \scrZ^{c, \lambda}$ defined by $\tr_{\Scom_{c, \lambda}} = \vartheta_{\beta_1}\cdots\vartheta_{\beta_r}$, where $w_0 = s_{\beta_1}\cdots s_{\beta_r}$ is a reduced expression for the longest element $w_0\in \Waff_{\lambda}$ with respect to any base for the finite root system $(\haff_{\bbR},\Raff_{\lambda})$. In particular, the hom-spaces in $\scrA^{c, \lambda}$ are free $\scrZ^{c, \lambda}$-modules due to \autoref{theo:basisloc}.
\begin{prop}\label{prop:duality}
	The sesquilinear pairing from \autoref{lemm:pairing} composed with $\tr_{\bfS^{\wedge}_{c, \lambda}}$ is a $\scrZ^{c,\lambda}$-bilinear perfect pairing
	\[
		\Hom_{\scrA^{c, \lambda}}(C'', C')\times \Hom_{\scrA^{c, \lambda}}(C, C'')\to \scrZ^{c, \lambda},\quad (\phi, \psi)\mapsto \tr_{\Scom_{c,\lambda}}\tr(\phi\psi).
	\] 
\end{prop}
\begin{proof}
	The pairing is $\scrZ^{c,\lambda}$-bilinear because the involution $\omega$ fixes the elements of $\scrZ^{c,\lambda}$ and $\tr_{\Scom_{c,\lambda}}$ is $\scrZ^{c,\lambda}$-linear. The pairing is perfect by \autoref{lemm:pairing} and the fact that $\tr_{\Scom_{c,\lambda}}$ induces an isomorphism $\Scom_{c,\lambda}\cong\Hom_{\scrZ^{c, \lambda}}(\Scom_{c,\lambda}, \scrZ^{c, \lambda})$. \end{proof}

\subsection{Generic chambers}
Given $d\in \frakP_{\bbR}$, let $\fraka_d = \left\{ d \right\}\times\haffd_{\bbR}\subseteq \fraka$ and let $\Ch^{c, \lambda}(\fraka_d)\subseteq \Ch^{c, \lambda}(\fraka)$ be the subset of chambers $C\in \Ch^{c, \lambda}(\fraka)$ such that $C\cap \fraka_d \neq \emptyset$. 
We say that a chamber $C\in \Ch^{c, \lambda}(\fraka_d)$ is \emph{generic} if no non-constant affine function on $\fraka_d$ is bounded on $C\cap \fraka_d$.

\begin{lemm}\label{lemm:generic}
	If $C\in \Ch^{c, \lambda}(\fraka_d)$ is a generic chamber, then $C\in \Ch^{c, \lambda}(\fraka_{d'})$ for every $d'\in \frakP_{\bbR}$ (in other words, $C$ is mapped surjectively onto $\frakP_{\bbR}$ by the projection $\fraka\to \frakP_{\bbR}$).
\end{lemm}
\begin{proof}
	Let $\ba\Psi^\pm_{c, \lambda, C}$ (resp. $\ba\Phi^\pm_{\lambda, C}$) be the image of $\Psi^\pm_{c, \lambda, C}$ (resp. $\Phi^\pm_{\lambda, C}$) in $\frakP_{\bbQ}^*\times \hfin_{\bbQ}^*$. The genericity implies there exists $\rho\in \hfin_{\bbQ}$ such that $\pm\varphi(\rho) > 0$ for every $\varphi\in \ba\Psi^\pm_{c, \lambda, C}\cup \ba\Phi^\pm_{\lambda, C}$. It follows that for every pair $(d', z)\in \frakP_{\bbR} \times \haffd_{\bbR}= \fraka$, there exists $r \gg 0$ such that $(d', z + r\rho) \in C$. Hence, $C\in\Ch^{c, \lambda}(\fraka_{d'})$ for every $d'\in \frakP_{\bbR}$.
\end{proof}

\section{Translation bimodules}\label{sec:transl}

\subsection{Definition}\label{subsec:shift}
For $d\in \frakP_{\bbZ}$, set $\bfB^{\shift{d}} :=  \Hom_{\bfA}(\kappa_0, \kappa_d)$. 
In particular, we have $\bfB^{\shift{0}} = \End_{\bfA}(\kappa_0)\cong \bfH$ by \autoref{theo:iota}. Given $d, d'\in \frakP_{\bbZ}$, we define the multiplication $\star:\bfB^{\shift{d'}}\otimes \bfB^{\shift{d}} \to \bfB^{\shift{d+d'}}$ via
\begin{align*}
	&\bfB^{\shift{d'}}\otimes \bfB^{\shift{d}} = \Hom_{\bfA}(\kappa_0, \kappa_{d'})\otimes\Hom_{\bfA}(\kappa_0, \kappa_{d})  \\
	&\xrightarrow{(t_d)_*\otimes \id} \Hom_{\bfA}(\kappa_{d}, \kappa_{d+d'})\otimes\Hom_{\bfA}(\kappa_0, \kappa_{d}) \xrightarrow{\circ}\Hom_{\bfA}(\kappa_0, \kappa_{d+d'}) = \bfB^{\shift{d+d'}}.
\end{align*}
It is easy to show that $\star$ is associative. In particular, $\bfB^{\shift{d}}$ has an $\bfH$-bimodule structure and $\star$ factorises through $\bfB^{\shift{d'}}\otimes \bfB^{\shift{d}}\to \bfB^{\shift{d'}}\otimes_{\bfH} \bfB^{\shift{d}}$. 
Note that $\bfB^{\shift{d}}$ satisfies 
\begin{equation}\label{equa:shift}
	(\bfc_\alpha + d_\alpha)\star a = a\star \bfc_\alpha \quad\text{ for $a\in\bfB^{\shift{d}}$ and $\alpha\in \Raff$}.
\end{equation} 

\begin{defi}
	For $d\in \frakP_{\bbZ}$, the $\bfH$-bimodule $\BB{}{d}$ is called the translation bimodule.
\end{defi}

Let $c, c'\in \frakP$ be such that $d = c - c'\in \frakP_{\bbZ}$. By~\eqref{equa:shift}, we have $\BB{}{d}\frakm_c = \frakm_{c'}\BB{}{d}\subset \BB{}{d}$, where $\frakm_c\subset \calO(\frakP)$ is the defining ideal of $c$. Put $\BB{c'}{c} = \BB{}{d} / \BB{}{d}\frakm_c$. It is an $(\bfH_{c'},\bfH_c)$-bimodule. 
\subsection{Length filtration}\label{subsec:FlA}
Recall the length filtration $F^{\lg}_{\bullet}\bfB^{\shift{d}} = F^{\lg}_{\bullet}\Hom_{\bfA}(\kappa_0, \kappa_d)$ defined in~\autoref{subsubsec:A}. By~\autoref{theo:basis}, the associated graded space has a free basis:
	\[
		\gr^{\lg}\bfB^{\shift{d}}  = \bigoplus_{w\in \Waff}\bfS\ba\tau_{\kappa_0, \kappa_d, w}, \quad \deg \ba\tau_{\kappa_0, \kappa_d,w} = d(\kappa_0, w^{-1} \kappa_d).
	\]
	For simplicity of notation, denote $\gamma^{d, w}:=\ba\tau_{\kappa_0, \kappa_d, w}\in \gr^{\lg}\bfB^{\shift{d}}$.
	\begin{lemm}\label{lemm:gamma}
		The multiplication $\star$ descends to the associated graded $\gr^{\lg}\bfB^{\shift{d}}$ and satisfies the following rule: 
		\begin{equation}\label{equa:multgamma}
			\gamma^{d, w}\star \gamma^{e, y} = \begin{cases}\gamma^{d+e, wy}&  y^{-1}\kappa_{d}\in [\kappa_0, (wy)^{-1}\kappa_{d+e}] \\ 0 & \text{otherwise} \end{cases},
		\end{equation}
		for $w,y\in \Waff$ and $d,e\in \frakP_{\bbZ}$. 
	\end{lemm}
	\begin{proof}
		The condition $y^{-1}\kappa_{d}\in [\kappa_0, (wy)^{-1}\kappa_{d+e}]$ (see~\autoref{lemm:minimal}) is equivalent to the following: if $G$ (resp. $G'$) is a minimal gallery from $\kappa_0$ to $y^{-1}\kappa_{d}$ (resp. from $y^{-1}\kappa_{d}$ to $(wy)^{-1}\kappa_{d+e}$), then $GG'$ is a minimal gallery from $\kappa_0$ to $(wy)^{-1}\kappa_{d+e}$. The statement follows immediately from~\autoref{theo:basis}.
	\end{proof}

\subsection{Harish-Chandra property}
The translation bimodule $\bfB^{\shift{d}}$ is equipped with the canonical filtration $F^{\can}_{\bullet}\bfB^{\shift{d}}$ from~\autoref{subsubsec:A} making it filtered $\bfH$-bimodule. It is the motivating example of our definition of Harish-Chandra bimodules:
\begin{prop}\label{prop:BHC}
	For each $d\in \frakP_{\bbZ}$, the canonical filtration $F^{\can}_{\bullet}\BB{}{d}$ is a HC filtration, so we have $\BB{}{d}\in \HC(\bfH)$.
\end{prop}
\begin{proof}
	Recall the homogeneous version $\bfB^{\shift{d}}_{\delta} = \Hom_{\bfA_{\delta}}(\kappa_0, \kappa_d)$ defined in \autoref{subsec:Adelta}. Put $\ba\bfB^{\shift{d}} := \bfB^{\shift{d}}_{\delta} / (\delta) \cong \gr^{\can}\bfB^{\shift{d}}$. 
	The inclusion $\gr^{\can}\End_\bfA(\kappa_0)\hookrightarrow \gr^{\can}\End_{\bfA^o}(\kappa_0)$ becomes an isomorphism after inverting the images of elements of $\Psi$ in $\bfS_{\delta} / (\delta)$. In particular, it preserves the centre: 
	\[
		\rmZ(\gr^{\can}\End_{\bfA}(\kappa_0)) = \rmZ(\gr^{\can}\End_{\bfA^o}(\kappa_0))\cap \gr^{\can}\End_{\bfA}(\kappa_0).
	\]
	Similarly, we have $\ba\bfB^{\shift{d}}\hookrightarrow \gr^{\can}\Hom_{\bfA^o}(\kappa_0, \kappa_d)\cong \gr^{\can}\End_{\bfA^o}(\kappa_0)$. Therefore, for any $x\in \rmZ(\gr^{\can}\End_{\bfA}(\kappa_0))$ and $b\in \ba\bfB^{\shift{d}}$, we have $x \star b = ((t_d)_*x) b = xb = bx = b\star x$; the second equation is due to the fact that $(t_d)_*$ is reduced to the identity map on $\gr^{\can}\bfS = \calO(\frakP\times \frakh)$. This together with \autoref{lemm:B-fg} below implies that the canonical filtration on $\BB{}{d}$ is a HC filtration.
\end{proof}
\begin{lemm}\label{lemm:B-fg}
	$\bfB^{\shift{d}}_{\delta}$ is finitely generated as a graded left $\bfH_{\delta}$-module.
\end{lemm}
\begin{proof}
	It suffices to show that $\gr^{\lg}\bfB^{\shift{d}}_{\delta}$ is finitely generated as a graded left $\gr^{\lg}\bfH_{\delta}$-module. We first show that for each $w\in \Waff$ such that $\ell(w) \gg 0$, there exists $\alpha\in \Deltaaff$ such that $\gamma^{0, s_\alpha}\star\gamma^{d, s_{\alpha}w} = \gamma^{d, w}$. This will imply the required statement because $\gr^{\lg}\bfB^{\shift{d}}_{\delta}$ is spanned by $\left\{\gamma^{d, w} \right\}_{w\in \Waff}$ as a graded left $\bfS_{\delta}$-module. \par
	Let $w\in \Waff$ and suppose that $\ell(w) \gg 0$. By the decomposition $\Waff = Q^{\vee}\rtimes  W $, we can write $w = t_{\mu}w_1$ for $w_1\in  W $ and $\mu\in Q^{\vee}$. Moreover, there exists positive integers $\{e_{\alpha}\}_{\alpha\in \Deltaaff}$ satisfying $\sum_{\alpha\in \Deltaaff} e_{\alpha}\alpha = \delta$, so that
	\[
		\sum_{\alpha\in \Deltaaff} e_{\alpha}\langle\alpha, \mu\rangle = \langle \delta, \mu\rangle = 0;
	\]
	it implies that when $\ell(w)\gg 0$, we have $\ell(t_{\mu}) \ge \ell(w) - \# W \gg 0$, and there exists $\alpha\in \Deltaaff$ such that $\langle \alpha, \mu\rangle \ll 0$; it follows that for such $\alpha$, we have
	\[
	((t_{-\mu})_*\alpha - \bfc_{\alpha})(w_1 \kappa_d) = (\alpha - \bfc_{\alpha})(w_1\kappa_d) + \langle \alpha, \mu\rangle \subseteq \bbR_{<0}
	\]
	and thus $\alpha - \bfc_{\alpha}\in (t_{\mu})_*\Psi^-_{w_1\kappa_d} = \Psi^-_{w\kappa_d}$; similarly, $-\alpha - \bfc_{\alpha}\in \Psi^+_{w\kappa_d}$. Let $w' = s_{\alpha}w$. We have thus
	\begin{align*}
		\Psi^+_{\kappa_0}\cap \Psi^+_{w'\kappa_d} = (\Psi^+_{\kappa_0} \cap s_{\alpha}\Psi^+_{w'\kappa_d})\sqcup \left\{ \alpha - \bfc_{\alpha} \right\}, \\
	(\Psi^-_{\kappa_0}\cap \Psi^-_{w'\kappa_d})\sqcup \left\{ -\alpha - \bfc_{\alpha} \right\} = \Psi^-_{\kappa_0} \cap s_{\alpha}\Psi^-_{w'\kappa_d}.
	\end{align*}
	The first equation implies $\Psi^+_{\kappa_0} \cap \Psi^+_{w\kappa_d} \subseteq \Psi^+_{w'\kappa_d}$ and the second implies $\Psi^-_{\kappa_0} \cap \Psi^-_{w\kappa_d}\subseteq \Psi^-_{w'\kappa_d}$, or equivalently $\Psi^+_{w'\kappa_d}\subseteq \Psi^+_{\kappa_0} \cup \Psi^+_{w\kappa_d}$. Similar arguments show that $\Phi^+_{\kappa_0} \cap \Phi^+_{w\kappa_d}\subseteq \Phi^+_{w'\kappa_d}\subseteq \Phi^+_{\kappa_0} \cup \Phi^+_{w\kappa_d}$ holds. By~\autoref{lemm:minimal}, these combined imply that the condition $w'\kappa_d\in [\kappa_0, w\kappa_d]$ in~\autoref{lemm:gamma} is satisfied, and therefore $\gamma^{0, s_{\alpha}}\star \gamma^{d, s_{\alpha}w} = \gamma^{d, w}$ holds. 
\end{proof}

\subsection{Spectral completion of translation bimodules}\label{subsec:Bcompl}
The translation bimodule has been defined in terms of the category $\bfA$. We shall describe its spectral completion in terms of the completed category $\scrA^{c, \lambda}$, introduced in~\autoref{subsec:Aloc}, and establish a relation with the algebraic KZ functor $\bfV$.
\subsubsection{}
Let $[\lambda]\in \haffd / \Waff$ be a $\Waff$-orbit let $c,c'\in \frakP$ be such that $d = c - c'\in \frakP_{\bbZ}$. Set $\Bscr{c'}{c}{\lambda} = \scrC^{\lambda}_{c'}(\BB{c'}{c})$.

\begin{prop}
	There is a natural isomorphism
	\begin{equation}\label{equa:Bscr}
		\Bscr{c'}{c}{\lambda}\cong\bigoplus_{[w], [w']\in \Waff / \Waff_{\lambda}}\Hom_{\scrA^{c, \lambda}}(\widetilde{w^{-1}\kappa_0}, \widetilde{w'^{-1}\kappa_d})/(\bfc - c).
	\end{equation}
	In particular, when $c' = c$, we have
	\begin{equation}\label{equa:Hscr}
		\scrH_{c}^{\lambda}\cong\bigoplus_{[w], [w']\in \Waff / \Waff_{\lambda}}\Hom_{\scrA^{c, \lambda}}(\widetilde{w^{-1}\kappa_0}, \widetilde{w'^{-1}\kappa_0})/(\bfc - c).
	\end{equation}
\end{prop}
\begin{proof}
 The spectral completion can be described as follows:
\[
	\Bscr{c'}{c}{\lambda} = \calO^{\wedge}_{[\lambda],\haffd}\otimes_{\calO(\haffd)}\BB{c'}{c},\quad \calO^{\wedge}_{[\lambda],\haffd} = \bigoplus_{\lambda'\in [\lambda]}\calO^{\wedge}_{\lambda',\haffd}.
\]
Moreover, there is an isomorphism $\BB{c'}{c}\otimes_{\calO(\haffd)}\calO^{\wedge}_{[\lambda],\haffd}\cong\Bscr{c'}{c}{\lambda}$ derived from~\eqref{equa:CHomA}.  By~\autoref{theo:comparaison}, there is an isomorphism for each $w, w'\in \Waff$:
\[
	1_{w'\lambda}\cdot\Bscr{c'}{c}{\lambda}\cdot 1_{w\lambda} \cong \Hom_{\scrA^{c', \lambda}}(\widetilde{w^{-1}\kappa_{-d}}, \widetilde{w'^{-1}\kappa_0}) / (\bfc - c').
\]
Applying the translation equivalence $t_d:\scrA^{c', \lambda}\xrightarrow{\sim}\scrA^{c, \lambda}$ from~\autoref{subsubsec:transeq}, we obtain
\[
	\Hom_{\scrA^{c', \lambda}}(\widetilde{w^{-1}\kappa_{-d}}, \widetilde{w'^{-1}\kappa_0}) / (\bfc - c')\cong \Hom_{\scrA^{c, \lambda}}(\widetilde{w^{-1}\kappa_0}, \widetilde{w'^{-1}\kappa_d}) / (\bfc - c).
\]
Taking summation over cosets $[w],[w']\in \Waff / \Waff_{\lambda}$, we obtain the isomorphism \eqref{equa:Bscr}. 
\end{proof}

\subsubsection{}
Let $(c, \lambda)\in \frakP\times \haffd$. Recall the algebraic KZ functor $\bfV$ and the completed affine Hecke algebra $\scrK^{\ell}_t$ from~\autoref{subsec:KZ}. Let $\scrP_c^{\lambda} = \bfV(\scrH^{\lambda}_c)$; it is a $(\scrK^{\ell}_t, \scrH^{\lambda}_c)$-bimodule.
\begin{prop}\label{prop:PBP}
	For $\lambda\in \haffd$ and $c,c'\in \frakP$ such that $d = c - c'\in \frakP_{\bbZ}$, there is a natural isomorphism
	\[
		\alpha_{c', c}: \scrP_{c'}^{\lambda}\otimes_{\scrH_{c'}^{\lambda}}\Bscr{c'}{c}{\lambda}\cong \scrP_{c}^{\lambda}.
	\]
	Moreover, if we are given $c''\in \frakP$ such that $d' = c' - c''\in \frakP_{\bbZ}$, then the following square is commutative:
	\begin{equation}\label{equa:PBB}
		\begin{tikzcd}
			\scrP^{\lambda}_{c''}\otimes_{\scrH^{\lambda}_{c''}}\Bscr{c''}{c'}{\lambda}\otimes_{\scrH^{\lambda}_{c'}}\Bscr{c'}{c}{\lambda} \arrow{r}{1\otimes \star}\arrow{d}[swap]{\alpha_{c'', c'}\otimes 1}& \scrP^{\lambda}_{c''}\otimes_{\scrH^{\lambda}_{c''}}\Bscr{c''}{c}{\lambda} \arrow{d}{\alpha_{c'',c}}\\
			\scrP^{\lambda}_{c'}\otimes_{\scrH^{\lambda}_{c'}}\Bscr{c'}{c}{\lambda} \arrow{r}{\alpha_{c', c}}& \scrP^{\lambda}_{c}
		\end{tikzcd}.
	\end{equation}
\end{prop}
\begin{proof}
	The algebraic KZ functor is given by $\bfV: M\mapsto 1_{\bfV} M$, where $1_{\bfV}\in \scrH^{\lambda}_{c}$ is the sum $1_{\bfV}  = \sum_{w\in \Sigma} 1_{c, w\lambda}$ for some finite set $\Sigma\subseteq \Waff$ such that $w^{-1}\nu_0$ is in a generic $(c,\lambda)$-clan for each $w\in \Sigma$, or equivalently, $\left\{ 0 \right\}\times w^{-1}\nu_0\subseteq \fraka$ lies in a generic chamber, denoted by $C(w^{-1})$. Thus, \eqref{equa:Hscr} implies that we can express
	\[
		\scrP^{\lambda}_c = \bigoplus_{[w]\in \Waff / \Waff_{\lambda}}\bigoplus_{y\in \Sigma} \Hom_{\scrA^{c, \lambda}}(\widetilde{w^{-1}\kappa_0}, \widetilde{y^{-1}\kappa_0}) / (\bfc - c).
	\]
	On the other hand, \autoref{lemm:generic} implies that $C(y^{-1})\cap (\left\{ d \right\}\times \haffd_{\bbR})\neq \emptyset$ for each $y\in \Sigma$. We may assume that $\left\{ d \right\}\times y^{-1}\nu_0\subseteq C(y^{-1})$ for $y\in \Sigma$ by choosing $\Sigma$ such that each $y^{-1}\nu_0$ is far in the Weyl chamber. Hence, $\widetilde{y^{-1}\kappa_d} = C(y^{-1}) = \widetilde{y^{-1}\kappa_0}$ holds for each $y\in \Sigma$ and \eqref{equa:Bscr} yields isomorphisms
	\begin{align*}
		\scrP_{c'}^{\lambda}\otimes_{\scrH_{c'}^{\lambda}}\Bscr{c'}{c}{\lambda} \cong \bigoplus_{[w]\in \Waff / \Waff_{\lambda}}\bigoplus_{y\in \Sigma} \Hom_{\scrA^{c, \lambda}}(\widetilde{w^{-1}\kappa_0}, \widetilde{y^{-1}\kappa_d}) / (\bfc - c) \\
		= \bigoplus_{[w]\in \Waff / \Waff_{\lambda}}\bigoplus_{y\in \Sigma} \Hom_{\scrA^{c, \lambda}}(\widetilde{w^{-1}\kappa_0}, \widetilde{y^{-1}\kappa_0}) / (\bfc - c) \cong \scrP_{c}^{\lambda}.
	\end{align*}
	The isomorphism $\alpha_{c',c}$ is then defined as the composite. The commutativity of~\eqref{equa:PBB} results from the associativity of the composition of morphisms in $\scrA^{c, \lambda}$.
\end{proof}
\subsection{Double centraliser property}

Let $c, c'\in \frakP$ be such that $c - c'\in \frakP_{\bbZ}$.
\begin{prop}\label{prop:bicommutante}
	The following ring homomorphisms induced by the bimodule structure
	\begin{equation}\label{equa:func}
		\bfH_c^{\op} \to \End_{\bfH_{c'}}(\BB{c'}{c}),\quad \bfH_{c'} \to \End_{\bfH_{c}^{\op}}(\BB{c'}{c})
	\end{equation}
	are isomorphisms.
\end{prop}
\begin{proof}
	We will only show that the first map of~\eqref{equa:func} is an isomorphism, the proof for second map being similar.
	Since the $\bfS$-algebra $\scrS$ is faithfully flat by~\autoref{lemm:platitude}, it suffices to show that~\eqref{equa:func} becomes an isomorphism after applying the spectral completion $\scrC$. Since both sides of~\eqref{equa:func}, viewed as left $\calO(\frakP)$-modules, are punctually supported on $c\in \frakP$, it suffices to show that the map~\eqref{equa:func} becomes an isomorphism after applying $\scrC^{\lambda}_c$ for each $[\lambda]\in \haffd / \Waff$. By~\autoref{prop:HCExt}, the spectral completion becomes 
	\begin{equation}\label{equa:func-i}
		(\scrH^{\lambda}_c)^{\op} \to \ho_{\scrH^{\lambda}_{c'}}(\Bscr{c'}{c}{\lambda}, \Bscr{c'}{c}{\lambda}).
	\end{equation}
	This map is an isomorphism due to~\autoref{lemm:bicomm} below.
\end{proof}
\begin{lemm}\label{lemm:bicomm}
	Let $[\lambda]\in \haffd / \Waff$. Then, given compact projective modules $P, Q\in \scrH^{\lambda}_c\proj$, the following natural map
	\begin{equation}\label{equa:func-ii}
		\Hom_{\scrH^{\lambda}_c}(P, Q)\to \Hom_{\scrH^{\lambda}_{c'}}(\Bscr{c'}{c}{\lambda}\otimes_{\scrH^{\lambda}_c}P, \Bscr{c'}{c}{\lambda}\otimes_{\scrH^{\lambda}_c}Q)
	\end{equation}
	is an isomorphism.
\end{lemm}
\begin{proof}
	For simplicity of notation, we denote $H = \scrH^{\lambda}_c$, $H' = \scrH^{\lambda}_{c'}$ and $B = \Bscr{c'}{c}{\lambda}$. We apply the algebraic KZ functor $\bfV$ from~\autoref{theo:KZ} and consider the composition:
	\begin{equation}\label{equa:VB}
		\Hom_{H}(P, Q)\xrightarrow{\text{\eqref{equa:func-ii}}} \Hom_{H'}(B\otimes_{H}P, B\otimes_{H}Q) \xrightarrow{\bfV_*} \Hom_{K}(\bfV(B\otimes_{H}P), \bfV(B\otimes_{H}Q)).
	\end{equation}
	By~\autoref{prop:PBP}, there is a natural isomorphism $\bfV(B\otimes_H P)\cong \bfV(P)$ and the composite of~\eqref{equa:VB} can be identified with the map of $\bfV$ on the hom-space: 
	\[
		\Hom_{H}(P, Q)\to \Hom_{H}(\bfV(P), \bfV(Q)),
	\]
	which by the double centraliser property of $\bfV$ (see~\autoref{theo:KZ}), is an isomorphism; consequently, the map $\bfV_*$ in~\eqref{equa:VB} is surjective. To prove that~\eqref{equa:func-ii} is an isomorphism, it remains to show that $\bfV_*$ is injective. \par

	Recall the ring $\scrZ_c^{\lambda} = (\bfS^\wedge_{c, \lambda}/(\bfc - c))^{\Waff_\lambda}$ from \autoref{subsubsec:Z} Given a $\scrZ_c^{\lambda}$-module $M$, let $M_F = M\otimes_{\scrZ_c^{\lambda}}F$ denote the base change to the field of fractions $F := \Frac\scrZ_c^{\lambda}$. We have a commutative square
	\[
		\begin{tikzcd}
			\Hom_{H}(B\otimes_H P, B\otimes_H Q)\arrow{r}{\bfV_*} \arrow[hookrightarrow]{d} & \Hom_{K}(\bfV(B\otimes_H P), \bfV(B\otimes_H Q))\arrow[hookrightarrow]{d} \\
			\Hom_{H_F}((B\otimes_H P)_F, (B\otimes_H Q)_F)\arrow{r}{\sim} & \Hom_{K_F}((\bfV(B\otimes_H P))_F, (\bfV(B\otimes_H Q))_F) 
		\end{tikzcd}
	\]
	where the vertical maps are injective (because $B$, $H$, $P$ and $Q$ are free over $\scrZ_c^{\lambda}$, see \autoref{subsubsec:duality}); the lower horizontal map is an isomorphism because $(\bfV(B\otimes_H P))_F = \Hom_{H_F}(H_F1_{\bfV}, P)$ and $H_F 1_{\bfV}$ is a compact projective generator of $H_F\Mod$, which is equivalent to $F\Mod$ by~\autoref{lemm:FracZ}; hence, the upper horizontal map $\bfV_*$ is injective. This completes the proof. 
\end{proof}

\subsection{Translation bimodules along a semigroup}\label{subsec:translation-semigroup}
We provide here some results on the $\bbN$-graded ring $\bfB^{\otimes}$ formed by the translation bimodules along the semigroup generated by an element $d\in \frakP_{\bbZ}$. We prove the noetherianity of $\bfB^{\otimes}$. It serves as preparation for the key technical \autoref{lemm:finitude} regarding the generic freeness over the parameter space $\frakP$ of certain Harish-Chandra bimodules.

\subsubsection{}
Fix $d\in \frakP_{\bbZ}$. Form the following ring:
\[
	\bfB^{\otimes} = \bigoplus_{n \ge 0}\bfB^{\shift{nd}},
\]
whose multiplication is given by the map $\star$ from~\autoref{subsec:shift}. It has a $\bbZ^2$-filtration given by
\[
	F_{n, m}\bfB^{\otimes} = \bigoplus_{k \ge 0}F^{\lg}_n F^{\can}_{n+m}\bfB^{\shift{kd}},
\]
where the first filtration on the right-hand side is the length filtration (\autoref{subsec:FlA}) and the second is the canonical filtration (\autoref{subsec:shift}). Note that $F_{n,m} \bfB^{\otimes} = 0$ whenever $n < 0$ or $m < 0$. 

\begin{prop}\label{prop:Bnoeth}
	The graded ring $\bfB^{\otimes}$ is left and right noetherian.
\end{prop}
\begin{proof}
	Consider the $\bbN^2$-graded Rees ring
	\[
		\bfB^{\otimes}_{\zeta,\delta} = \bigoplus_{n,m \ge 0}\zeta^n \delta^m F_{n, m}\bfB^{\otimes}.
	\]
	By the technique of associated graded~\cite[Appendix D]{hotta2007d}, it suffices to show that the associated quotient $\ba\bfB^{\otimes}:=\bfB^{\otimes}_{\zeta,\delta} / (\zeta,\delta)$ is left and right noetherian. \autoref{lemm:fg} below establishes this property.
\end{proof}
\begin{lemm}\label{lemm:fg}
	Let $\ba \bfB^{\otimes}$ be as in the proof of~\autoref{prop:Bnoeth}. Then, $\ba \bfB^{\otimes}$ contains a commutative graded subring $A$ of finite type over $\bbC$ such that $\ba\bfB^{\otimes}$ is finitely generated as left and right $A$-module.
\end{lemm}
\begin{proof}
	By~\autoref{theo:basis}, the ring $\ba \bfB^{\otimes}$ admits the following decomposition:
	\[
		\ba \bfB^{\otimes} = \bigoplus_{n\in \bbN} \bigoplus_{w\in \Waff}\ba\bfS \ba\gamma^{nd, w},
	\]
	where $\ba\bfS = \bfS_{\delta} / (\delta)$ and $\ba\gamma^{nd, w}$ is the image of $\gamma^{nd, w}\in \gr^F \bfB^{\otimes}$ in $\ba\bfB^{\otimes}$.
	 We claim that the subring
	\[
		A = \bigoplus_{n\in \bbN}\bigoplus_{\mu\in Q^{\vee}} \ba\bfS \ba\gamma^{nd, t_{\mu}}\subset \ba \bfB^{\otimes}
	\]
	satisfies the requirements. \par

	Indeed, for each chamber $U\subset \bbR_{\ge 0} d\times \hfin_{\bbR}$ of the hyperplane arrangement given by the finite set of linear functions $\left\{ \psi\mid_{\bbR d\times \hfin_{\bbR}}\in (\bbR d\times \hfin_{\bbR})^* \;;\; \psi\in \Psi\right\}$, let $\Sigma_{U}$ be a finite generating subset of the semigroup $(\bbN d\times Q^{\vee})\cap \ba U$. It is easy to show that for pairs $(jd, \mu), (kd, \nu)\in (\bbN d\times Q^{\vee})\cap \ba U$, we have $t_{\mu}\kappa_d\in [\kappa_0, t_{\mu+\nu}\kappa_{(j+k)d}]$, so that 
	\[
		\ba\gamma^{jd, t_{\mu}} \star \ba\gamma^{kd, t_{\nu}} = \ba\gamma^{(j+k)d, t_{\mu+\nu}}
	\]
	holds by~\autoref{lemm:gamma}.  Put $\Sigma = \bigcup_{U} \Sigma_{U}$. Then, $\left\{ \ba\gamma^{kd, t_{\mu}} \right\}_{(kd, \mu)\in \Sigma}$ is a generating subset of the ring $A$. The proof of finite generation of $\ba \bfB^{\otimes}$ over $A$ is more sophisticated and will be completed in~\autoref{subsubsec:fg}.
\end{proof}
\subsubsection{}
Let $\fraka' := \bbR d\times \haffd_{\bbR}\subseteq \fraka$. The hyperplane arrangement $\{H_{\psi}\}_{\psi\in \Phi\cup \Psi}$ introduced in~\autoref{subsec:Ch} can be restricted to $\fraka'$. Let $\Ch(\fraka')$ denote the set of chambers of this restriction. \par
We choose any euclidean metric $\|\cdot \|$ for $\haffd_{\bbR}$ and extend it to $\fraka'$ by setting $\|(rd, y) - (sd, y')\|^2 = \|y - y'\|^2 + |r - s|^2$. Let $Q^{\vee}\subset \hfin_{\bbR}$ be the coroot lattice of the finite root system $\Rfin$. The lattice $\bbZ d\times Q^{\vee}$ acts on $\fraka'$ by translation and the quotient $\frakQ = \fraka' / \bbZ d\times Q^{\vee}$ is compact. 
\begin{lemm}\label{lemm:interval1}
	Let $C_0, C_1\in \Ch(\fraka')$ be chambers. Let $\varepsilon > 0$, $p_0\in C_0$, $p_1\in C_1$ be such that the open ball $B_{\varepsilon}(p_i)\subseteq \fraka'$ is contained in $C_i$ for $i \in \left\{ 0,1 \right\}$. Suppose that the following inequality holds: 
	\[
		\frac{\|p_1 - p_0\|}{\varepsilon} \ge \frac{\Vol(\frakQ)}{\Vol(B_{\varepsilon/2}(p_0))} + 1.
	\]
	Then, there exists $\mu\in Q^{\vee}$ and $l\in \bbZ$ such that $(l, \mu)\neq (0, 0)$ and there exists $p_2\in C_0 + (ld, \mu)$ which lies in the interval $[p_0,p_1]$. 
\end{lemm}
\begin{proof}
	Let $m\in \bbN$ be an integer satisfying $\Vol(\frakQ) / \Vol(B_{\varepsilon/2}(p_0)) < m \le  \|p_1 - p_0\| / \varepsilon$. Put $\vec v = \varepsilon (p_1 - p_0) / \|p_1 - p_0\|$. Consider the family of open balls $\{B_{\varepsilon/2}(p_0 + k\vec v)\}_{k=0}^m$. Since the sum of the volume of these balls exceeds the volume $\Vol(\frakQ)$, there exists $n_1, n_2 \in \bbN$ satisfying $m\ge n_1 > n_2 \ge 0$ such that the images of $B_{\varepsilon/2}(p_0 + n_1\vec v)$ and $B_{\varepsilon/2}(p_0 + n_2\vec v)$ in the quotient $\frakQ$ overlap. In other words, there exists $\mu\in Q^{\vee}$ and $l\in \bbZ$ such that $\|(n_1 - n_2)\vec v - (ld, \mu)\| < \varepsilon$. Consequently, the point $p_0 + (n_1 - n_2)\vec v - (ld, \mu)\in B_{\varepsilon}(p_0)$ lies in the interior of $C_0$, or equivalently, the point $p_2:= p_0 + (n_1 - n_2)\vec v$ lies in the interior of $C_0 + (ld, \mu)$. The point $p_2$ satisfies the requirement: we have $(ld, \mu)\neq (0, 0)$ because $\|(n_1 - n_2)\vec v\| > \varepsilon$; in addition, $p_2$ lies in the interval $[p_0, p_1]$ because $n_1 - n_2 \le \|p_1 - p_0\|/ \varepsilon$. 
\end{proof}

\begin{lemm}\label{lemm:interval2}
	Let $p_0, p_1, p_2\in \fraka$ such that $p_2\in [p_0,p_1]$. For $i\in \left\{ 0, 1, 2 \right\}$, put 
	\[
		R_i =  \{\alpha\in \Phi \;;\; \alpha(p_i) > 0\},\quad S_i = \{\mu\in \Psi \;;\; \mu(p_i) > 0\}.
	\]
	Then, we have $R_0\cap R_1\subseteq R_2 \subseteq R_0 \cup R_1$ and $S_0\cap S_1\subseteq S_2 \subseteq S_0 \cup S_1$.
\end{lemm}
\begin{proof}
	The proof is straightforward and left to the reader.
\end{proof}

\subsubsection{}\label{subsubsec:fg}

Let's show that $\ba\bfB^{\otimes}$ is finitely generated left $A$-module. Let $\varepsilon > 0$ be a number such that for each chamber $C\in \Ch(\fraka')$, there exists $p_C\in C$ satisfying $B_\varepsilon(p_C)\subseteq C$ --- such $\varepsilon$ exists due to the finiteness of $(Q^{\vee}\times \bbZ d)$-orbits in $\Ch(\fraka')$. Set $N' =  \left(\Vol(\frakQ)/\Vol(B_{\varepsilon/2}(p_0)) + 1\right)\varepsilon$ and let $N \ge N'$ be an integer satisfying the following property: for each $w\in \Waff$ such that $\ell(w) > N$ and for each $(p,p')\in \nu_0 \times w\nu_0$, the condition $\|p' - p\| > N'$ holds. We claim that $\Upsilon = \left\{ \ba\gamma^{kd, w}\;;\; w\in \Waff,\;\ell(w) \le N ,\; 0\le k \le N \right\}$ is a generating set of $\ba\bfB^{\otimes}$ as both a left and right $A$-module. \par

Let $w\in \Waff$ and $n\in \bbN$. We prove by induction on $\ell(w) + n$ that $\ba\gamma^{nd, w}\in A\Upsilon$ holds. Suppose $\max(\ell(w), n) \le N$, then $\ba\gamma^{nd, w}\in \Upsilon$ holds trivially. Suppose then $\max(\ell(w), n) > N$. Let $C_0 = \kappa_0\cap \fraka'\subset \fraka'$ be the chamber containing $\left\{ 0 \right\}\times \nu_0$ and let $C_1\subset \fraka'$ be the chamber containing $\left\{ nd \right\}\times w\nu_0$. Let $p_0\in C_0$ and $p_1\in C_1$ be such that $B_{\varepsilon}(p_i)\subseteq C_i$ for $i\in \left\{ 0, 1 \right\}$. By our choice of $N$, it follows that $\|p_0 - p_1\| > N'$. Therefore, the hypothesis of~\autoref{lemm:interval1} is satisfied, which yields a pair $(0,0)\neq (l, \mu)\in \bbZ\times Q^{\vee}$ and $p_2\in [p_0, p_1] \cap C_0 + (ld, \mu)$; moreover, $0\le l\le n$ holds. Let $C_2 = C_0 + (ld, \mu)\in \Ch(\fraka')$. \autoref{lemm:interval2} applied to $(p_0, p_1, p_2)$ yields $S_0\cap S_1\subseteq S_2\subseteq S_0\cup S_1$ and $R_0\cap R_1\subseteq R_2\subseteq R_0\cup R_1$, which implies that $C_2\in [C_0, C_1]$ holds by \autorefitem{lemm:minimal}{i}.

It allows us to apply \autoref{lemm:gamma} and to obtain that $\ba\gamma^{ld, t_{\mu}}\star\ba\gamma^{(n-l)d, t_{-\mu}w} = \ba\gamma^{nd, w}$. Since $\ell(t_{\mu}) + \ell(t_{-\mu}w) = \ell(w)$ holds by the condition $C_2\in [C_0, C_1]$, we have $\ell(t_{-\mu}w) + (n - l) < \ell(w) + n$. By induction hypothesis, $\ba\gamma^{(n-l)d, t_{-\mu}w}\in A \Upsilon$ holds and thus $\ba\gamma^{nd, w}\in A \Upsilon$ also holds. This completes the induction step. Similarly, one can show that $\ba\gamma^{nd, w}\in \Upsilon A$ by interchanging $C_0$ and $C_1$ in the arguments. \par

Finally, as the family $\left\{ \ba\gamma^{nd, w} \right\}_{w\in \Waff, n\in \bbN}$ spans $\ba\bfB^{\otimes}$ as both a left and right $\ba\bfS$-module and $\ba\bfS\subseteq A$, we have $A\Upsilon = \ba\bfB^{\otimes} =\Upsilon A$. This completes the proof of~\autoref{lemm:fg}.

\section{Derived equivalences}\label{sec:equiv}
Let $c, c'\in \frakP$ be such that $d:= c - c'\in \frakP_{\bbZ}$ holds.  The $(\bfH_{c'}, \bfH_c)$-bimodule $\BB{c'}{c}$ from \autoref{subsec:shift} induces the following \emph{translation functor}:
\[
	\trans{c'}{c} = \BB{c'}{c}\otimes^{\rmL}_{\bfH_c}\relbar:\Db(\bfH_{c})\to \Db(\bfH_{c'}). 
\]
Since $\bfH$ is flat over $\calO(\frakP)$, we can also express $\trans{c'}{c} = \BB{}{d}\otimes^{\rmL}_{\bfH}\relbar$. The aim of this section is to establish the following theorem:
\begin{theo}\label{theo:main}
	The translation functor $\trans{c'}{c}$ is an equivalence of categories between $\Db(\bfH_{c})$ and $\Db(\bfH_{c'})$.
\end{theo}

\subsection{Stratification of \texorpdfstring{$\frakP$}{P}}
\subsubsection{}
We shall introduce a family $\frakM$ of linear functions on $\frakP_{\bbQ}$. Let $\ba\Psi$ be the image of $\Psi\subseteq \haff_{\bbQ}^*\times\frakP_{\bbQ}^*$ under the projection $\haff_{\bbQ}^*\times\frakP_{\bbQ}^*\to \hfin_{\bbQ}^*\times\frakP_{\bbQ}^*$; concretely, it is given by 
\[
	\ba\Psi = \left\{ \alpha - \bfc_{\alpha}\;;\; \alpha\in \Rfin_{\red} \right\}\cup \left\{ \alpha/2 - \bfc_{\alpha}\;;\; \alpha\in \Rfin \setminus \Rfin_{\red} \right\}. 
\]
For $\mu\in \ba\Psi$, let $\dot\mu\in R_{\red}$ denote its $\hfin^*_{\bbQ}$-component. A subset $\sigma\subseteq \ba\Psi$ is called \emph{circuit} if $\left\{ \dot\mu \right\}_{\mu\in \sigma}$ is a minimal $\bbQ$-linearly dependent family. If $\sigma\subseteq \ba\Psi$ is a circuit, then there is a family of non-zero rational numbers $\{ d_{\mu}\}_{\mu\in \sigma}$, unique up to scaling by $\bbQ^{\times}$, such that
\[
	\mu_{\sigma} := \sum_{\mu\in \sigma}d_\mu \mu\in \frakP^*_{\bbQ};
\]
namely, the $\hfin^*_{\bbQ}$-component of $\mu_{\sigma}$ vanishes. In other words, if $\mu_{\sigma}\neq 0$, then $[\mu_{\sigma}] := \bbQ^{\times}\mu_{\sigma}$ is an element of the projectivisation $\bbP(\frakP^*_{\bbQ})$ independent of the choice of $\left\{ d_{\mu} \right\}_{\mu\in \sigma}$. Define 
\[
	\frakM = \left\{ [\mu_{\sigma}]\in \bbP(\frakP^*_{\bbQ})\;;\; \sigma \subseteq \ba\Psi \text{ is a circuit},\; \mu_{\sigma}\neq 0\right\}.
\]

Recall the notation for root systems from \autoref{exam:BC}--\ref{exam:G2}.
\begin{prop}
	When $(\frakh_{\bbR}, \Rfin)$ is simply-laced, we have $\frakP^*_{\bbQ} = \bbQ \bfc_{\natural}$ and
	\begin{align*}
		\frakM = \left\{ [\bfc_{\natural}] \right\}.
	\end{align*}
	When $(\frakh_{\bbR}, \Rfin)$ is of type $BC_1$, we have $\frakP^*_{\bbQ} = \bbQ \bfc_{\sharp}\oplus \bbQ \bfc_{\flat}$ and
	\begin{align*}
		\frakM = &\left\{ [\bfc_{\sharp}], [\bfc_{\flat}] \right\}\cup \left\{[\bfc_{\sharp} + \epsilon \bfc_{\flat}]  \right\}_{\epsilon\in \left\{ \pm 1 \right\}}.
	\end{align*}
	When $(\frakh_{\bbR}, \Rfin)$ is of type $BC_n$ ($n\ge 2$), we have $\frakP^*_{\bbQ} = \bbQ \bfc_{\natural}\oplus \bbQ \bfc_{\sharp}\oplus \bbQ \bfc_{\flat}$ and
	\begin{align*}
		\frakM = &\left\{ [\bfc_{\natural}], [\bfc_{\sharp}], [\bfc_{\flat}] \right\} \cup \left\{ [\bfc_{\sharp} + \epsilon \bfc_{\flat}]  \right\}_{\epsilon\in \left\{ \pm 1 \right\}}\\
		& \cup \left\{[k\bfc_{\natural} + 2\epsilon \bfc_{\sharp}],  [k\bfc_{\natural} + 2\epsilon \bfc_{\flat}]\right\}_{\epsilon\in \left\{ \pm 1 \right\},k\in \left\{ 1, \ldots, n \right\}}  \\
		&\cup \left\{ [k\bfc_{\natural} + \epsilon_1 \bfc_{\sharp} + \epsilon_2 \bfc_{\flat}] \right\}_{\epsilon_1, \epsilon_2\in \left\{ \pm 1 \right\},k\in \left\{ 1, \ldots, n-1 \right\}}.
	\end{align*}
	When $(\frakh_{\bbR}, \Rfin)$ is of type $F_4$, we have $\frakP^*_{\bbQ} = \bbQ \bfc_{\natural}\oplus \bbQ \bfc_{\sharp}$ and
	\begin{align*}
	\frakM = &\left\{ [\bfc_{\natural}], [\bfc_{\sharp}] \right\} \cup \left\{ [\bfc_{\natural} + \epsilon \bfc_{\sharp}], [\bfc_{\natural} + 2\epsilon \bfc_{\sharp}], [\bfc_{\natural} + 3\epsilon \bfc_{\sharp}]\right\}_{\epsilon\in \left\{ \pm 1 \right\}} \\
	&\cup\left\{[\bfc_{\natural} + 4\epsilon \bfc_{\sharp}], [2\bfc_{\natural} + 3\epsilon \bfc_{\sharp}], [3\bfc_{\natural} + 2\epsilon \bfc_{\sharp}], [3\bfc_{\natural} + 4\epsilon \bfc_{\sharp}]\right\}_{\epsilon\in \left\{ \pm 1 \right\}}.
	\end{align*}
	When $(\frakh_{\bbR}, \Rfin)$ is of type $G_2$, we have $\frakP^*_{\bbQ} = \bbQ \bfc_{\natural}\oplus \bbQ \bfc_{\sharp}$ and 
	\begin{align*}
		\frakM = &\left\{ [\bfc_{\natural}], [\bfc_{\sharp}] \right\} \cup \left\{ [\bfc_{\natural} + \epsilon \bfc_{\sharp}], [\bfc_{\natural} + 2\epsilon \bfc_{\sharp}], [\bfc_{\natural} + 3\epsilon \bfc_{\sharp}] , [2\bfc_{\natural} + 3\epsilon \bfc_{\sharp}]\right\}_{\epsilon\in \left\{ \pm 1 \right\}}.
	\end{align*}
\end{prop}
\begin{proof}
	The proof is elementary and left to the reader.
\end{proof}

\subsubsection{}\label{sssec:c-clans}
For each parameter $c\in \frakP$ we introduce a sub-family of $\frakM$:
\[
	\frakM_{c} = \left\{ [\mu]\in \frakM\;;\; \mu(c)\in \bbQ \right\}.
\]
Then, the zero loci of the elements of $\frakM_c$ define a hyperplane arrangement on $\frakP_{\bbR}$, which decomposes $\frakP_{\bbR}$ into a union of facets, called \emph{$c$-facets}. $c$-Facets of maximal dimension are called \emph{$c$-chambers}. \par

For $c\in \frakP$ and $\lambda\in \haffd$, consider the subfamily $\Psi_{c, \lambda}\subseteq \Psi$ introduced in~\autoref{subsec:Aloc}. Let $\ba\Psi_{c,\lambda}$ be its image in $\frakP_{\bbQ}^* \times \hfin_{\bbQ}^*$. A connected component of the complement of the hyperplane arrangement $\left\{ H_{\mu} \right\}_{\mu\in \ba\Psi_{c, \lambda}}$ on $\fraka^0$ is called a \emph{$(c,\lambda)$-clan} of $\fraka^0$. Let $\Cl^{c,\lambda}(\fraka^0)$ denote the set of $(c, \lambda)$-clans of $\fraka^0$. 

\begin{lemm}\label{lemm:strate}
	Let $\frakC\in \Cl^{c, \lambda}(\fraka^0)$ be a $(c,\lambda)$-clan of $\fraka^0$. Then, the image of $\frakC$ under the projection $p:\fraka^0\to \frakP_{\bbR}$ is a union of $c$-facets in $\frakP_{\bbR}$. 
\end{lemm}
\begin{proof}
	Since the closure $\ba \frakC$ is a rational convex polyhedral cone of maximal dimension in $\fraka^0$, the same is true for its image $p(\ba \frakC)$ in $\frakP_{\bbR}$; moreover, $p(\frakC)$ is the interior of $p(\ba \frakC)$. Let $F$ be the closure of a face (facet of codimension $1$) of $p(\ba \frakC)$. Since $p(\ba \frakC)$ is a rational convex polyhedral cone, there exists a non-zero linear form $\mu_F\in \frakP_{\bbQ}^*$ such that $p(\ba \frakC)\subseteq \mu_F^{-1}(\bbR_{\ge 0})$ and $F = \mu_F^{-1}(0)\cap p(\ba \frakC)$. It is unique up to scaling by $\bbQ^{\times}$. It suffices to show that $[\mu_F] = \bbQ^{\times}\mu_F$ lies in $\frakM_c$ --- indeed, this will imply that $p(\frakC)$ is the intersection of the open half-spaces $\mu_F > 0$ over such $F$. \par

	The intersection $(\mu_F\circ p)^{-1}(0)\cap \ba \frakC$ is a union of facets of $\ba \frakC$ because $\ba \frakC\subseteq (\mu_F\circ p)^{-1}(\bbR_{\ge 0})$. Let 
	\[
		\sigma_0 = \left\{ a\in \ba\Psi_{c, \lambda}\;;\; (\mu_F\circ p)^{-1}(0)\cap \ba \frakC\subseteq a^{-1}(0)\right\}.
	\]
	It is clear that the intersection $\bigcap_{a\in \sigma_0}a^{-1}(0)$ is the subspace spanned by the convex polyhedral cone $(\mu_F\circ p)^{-1}(0)\cap \ba \frakC$. In particular, the linear subspace $p(\bigcap_{a\in \sigma_0}a^{-1}(0)) \subseteq \frakP_{\bbR}$ coincides with the hyperplane $\mu_F^{-1}(0)$. Let $\sigma\subseteq \sigma_0$ be a minimal subset satisfying $p(\bigcap_{a\in \sigma}a^{-1}(0))  =\mu_F^{-1}(0)$. Then, it follows that $\sigma$ is a circuit --- indeed, let $V = \bigcap_{a\in \sigma}a^{-1}(0)\subseteq \fraka$; then, $p(V) = \mu_F^{-1}(0)$ and 
	\[
		\bigcap_{a\in \sigma}\dot a^{-1}(0) = p^{-1}(0) \cap V
	\]
	imply
	\[
		\dim\left(\bigcap_{a\in \sigma}\dot a^{-1}(0)\right) = \dim V - \dim p(V) = \dim \frakh_{\bbR}- \dim \frakP + 1 \ge \dim \frakh_{\bbR} - \# \sigma + 1,
	\]
	so that $\sigma$ is linearly dependent; 
	with similar arguments, the minimality of $\sigma$ implies that $\left\{ \dot a \right\}_{a\in \sigma'}$ is linearly independent for every proper subset $\sigma'\subsetneq \sigma$; therefore, $\sigma$ is a circuit. This implies that $[\mu_F] = [\mu_{\sigma}]$, where $\mu_{\sigma} = \sum_{a\in \sigma}d_a a$ for some $d_a\in \bbQ^\times$. For each $a\in \sigma$, let $\til a\in \Psi_{c, \lambda}$ be a lifting of $a$ along the projection $\Psi_{c, \lambda}\to \ba\Psi_{c, \lambda}$ and set $\til\mu_{\sigma}= \sum_{a\in \sigma}d_a \til a$. It follows that $\mu_{\sigma}(c) \in \til\mu_{\sigma}(c,\lambda) + \bbQ$ and $\til\mu_{\sigma}(c,\lambda) = \sum_{a\in \sigma}d_a \til a(c, \lambda)  = 0$; in other words, $[\mu_F] = [\mu_{\sigma}]\in \frakM_c$ holds. 
\end{proof}

\subsubsection{}\label{subsubsec:strat}
 Given a $c$-facet $F\subseteq \frakP_{\bbR}$, we define its complexification $F_{\bbC} := F + i\frakP_{\bbR}$ and call it \emph{$c$-facet} of $\frakP$. Then, $\left\{ F_{\bbC} \right\}_F$ forms a facet complex on $\frakP$, where $F$ runs over the $c$-facets of $\frakP_{\bbR}$. 
\subsection{Criteria of equivalence}
\subsubsection{}
Let $c, c'\in \frakP$ be such that $d:= c - c'\in \frakP_{\bbZ}$. 
\begin{lemm}\label{lemm:crit}
	The following conditions are equivalent:
	\begin{enumerate}
		\item\label{lemm:crit-i}
			The functor $\trans{c'}{c}: \Db(\bfH_c)\to \Db(\bfH_{c'})$ is an equivalence of categories.
		\item\label{lemm:crit-ii}
			The left $\bfH_{c'}$-module $\BB{c'}{c}$ is tilting, namely, $\BB{c'}{c}$ generates $\Db_{\perf}(\bfH_{c'})$ as thick triangulated subcategory and $\Ext^{>0}_{\bfH_{c'}}(\BB{c'}{c}, \BB{c'}{c}) = 0$ holds.
		\item\label{lemm:crit-iii}
			The following adjoint unit and counit are quasi-isomorphisms
			\begin{align}
				\epsilon:\bfH_c \to \RHom_{\bfH_{c'}}(\BB{c'}{c}, \BB{c'}{c}) \label{equa:adjunit} \\
				\eta:\BB{c'}{c}\otimes^{\rmL}_{\bfH_c}\RHom_{\bfH_{c'}}(\BB{c'}{c}, \bfH_{c'})\to \bfH_{c'}.  \label{equa:adjcounit}
			\end{align}
		\item[(i')\namedlabel{lemm:crit-i'}{(i')}]
			For each $\lambda\in \haffd$, the functor 
			\[
				\trans{c'}{c}^{\lambda} = \Bscr{c'}{c}{\lambda}\otimes^{\rmL}_{\scrH^{\lambda}_c}\relbar: \Db(\scrH^{\lambda}_c)\to \Db(\scrH^{\lambda}_{c'})
			\]
			is an equivalence of categories.
		\item[(ii')\namedlabel{lemm:crit-ii'}{(ii')}]
			For each $\lambda\in \haffd$, the left $\scrH^{\lambda}_{c'}$-module $\Bscr{c'}{c}{\lambda}$ is tilting.
	\item[(iii')\namedlabel{lemm:crit-iii'}{(iii')}]
			For each $\lambda\in \haffd$, the following morphisms are quasi-isomorphisms
			\begin{align}
				\scrH^{\lambda}_c \to \Rhom_{\scrH^{\lambda}_{c'}}(\Bscr{c'}{c}{\lambda}, \Bscr{c'}{c}{\lambda}) \label{equa:adjunit-} \\
				\Bscr{c'}{c}{\lambda}\otimes^{\rmL}_{\scrH^{\lambda}_c}\Rhom_{\scrH^{\lambda}_{c'}}(\Bscr{c'}{c}{\lambda}, \scrH^{\lambda}_{c'})\to \scrH^{\lambda}_{c'}.  \label{equa:adjcounit-}
			\end{align}

	\end{enumerate}
\end{lemm}
\begin{proof}
	The equivalence~\ref{lemm:crit-i} $\Leftrightarrow$ \ref{lemm:crit-ii} follows from~\cite[1.6]{happel87}, the condition $\bfH_{c'}\cong\End_{\bfH_c^{\op}}(\BB{c'}{c})$ being ensured by~\autoref{prop:bicommutante}. The assertion~\ref{lemm:crit-i} $\Rightarrow$ \ref{lemm:crit-iii} is trivial. The assertion~\ref{lemm:crit-iii} $\Rightarrow$ \ref{lemm:crit-ii} follows from the fact that $\bfH_{c'}$ generates $\Db_{\perf}(\bfH_{c'})$ as thick triangulated subcategory due to~\autoref{prop:dimglob}. The equivalences~\ref{lemm:crit-i'} $\Leftrightarrow$ \ref{lemm:crit-ii'} $\Leftrightarrow$ \ref{lemm:crit-iii'} can be proven similarly. The equivalence \ref{lemm:crit-iii} $\Leftrightarrow$ \ref{lemm:crit-iii'} follows from the faithful flatness of $\scrS$ over $\bfS$ proven in~\autoref{lemm:platitude} and the formul\ae~\autorefitem{prop:HCTor}{ii} and \autorefitem{prop:HCExt}{ii}.
\end{proof}
\subsubsection{}
Consider the following condition for $c, c'\in \frakP$ such that $c - c'\in \frakP_{\bbZ}$:
\begin{quote}
	$\EQ{c}{c'}$:\quad\begin{minipage}[t]{.8\linewidth}
		the multiplication $\star$ induces isomorphisms 
		\[
		\BB{c}{c'}\otimes_{\bfH_{c'}}\BB{c'}{c}\xrightarrow{\cong}\bfH_{c},\quad \BB{c'}{c}\otimes_{\bfH_{c}}\BB{c}{c'}\xrightarrow{\cong}\bfH_{c'}.
	\]
		\end{minipage}
\end{quote}
\begin{lemm}\label{lemm:t-exact}
Suppose that $\EQ{c}{c'}$ holds, then 
	\begin{enumerate}
		\item 
			$\trans{c'}{c}$ and $\trans{c}{c'}$ are t-exact equivalences and are inverse to each other;
		\item
			$\trans{c'}{c}\circ\trans{c}{c''} \cong\trans{c'}{c''}$ and $\trans{c''}{c'}\circ\trans{c'}{c} \cong\trans{c''}{c}$ hold for each $c''\in c + \frakP_{\bbZ}$.
	\end{enumerate}
\end{lemm}
\begin{proof}
	The conditions imply that the functors
	\[
		\BB{c'}{c}\otimes_{\bfH_c}\relbar:\bfH_c\Mod \to\bfH_{c'}\Mod,\quad\BB{c'}{c}\otimes_{\bfH_c}\relbar:\bfH_{c'}\Mod \to\bfH_{c}\Mod
	\]
	are equivalences of categories and inverse to each other. In particular, $\BB{c'}{c}$ is projective as left $\bfH_{c'}$-module and as right $\bfH_c$-module and subsequently, we have $\BB{c}{c'}\otimes^{\rmL}_{\bfH_{c'}}\BB{c'}{c}\xrightarrow{\sim}\bfH_{c}$ and $\BB{c'}{c}\otimes^{\rmL}_{\bfH_{c}}\BB{c}{c'}\xrightarrow{\sim}\bfH_{c'}$.
	
	The associativity of the multiplication $\star$ yields a commutative diagram
	\[
		\begin{tikzcd}
			\BB{c}{c'}\otimes^{\rmL}_{\bfH_{c'}}\BB{c'}{c}\otimes^{\rmL}_{\bfH_{c}}\BB{c}{c''} \arrow{r} \arrow{d}& \BB{c}{c'}\otimes^{\rmL}_{\bfH_{c'}}\BB{c'}{c''} \arrow{d}\\
			\BB{c}{c}\otimes^{\rmL}_{\bfH_{c}}\BB{c}{c''} \arrow{r}{\cong} & \BB{c}{c''} 
		\end{tikzcd},
	\]
	We have seen that the left vertical arrow is an isomorphism; in particular, the upper horizontal arrow admits a left inverse and the right vertical arrow admits a right inverse. Applying $\BB{c'}{c}\otimes^{\rmL}_{\bfH_c}\relbar$ to the diagram, we see that the morphism 
	\[
		m = 1\otimes \star:\BB{c'}{c}\otimes^{\rmL}_{\bfH_{c}}\BB{c}{c'}\otimes^{\rmL}_{\bfH_{c'}}\BB{c'}{c''} \to \BB{c'}{c}\otimes^{\rmL}_{\bfH_{c}}\BB{c}{c''}
	\]
	admits a right inverse. Interchanging $c$ and $c'$, we deduce that the multiplication $m$ also admits a left inverse. Therefore, $m$ is a quasi-isomorphism. Applying $\BB{c}{c'}\otimes^{\rmL}_{\bfH_{c'}}\relbar$ to $m$, and using the isomorphism $\BB{c}{c'}\otimes^{\rmL}_{\bfH_{c'}}\BB{c'}{c}\xrightarrow{\sim}\bfH_{c}$, we conclude that $\BB{c}{c'}\otimes^{\rmL}_{\bfH_{c'}}\BB{c'}{c''} \cong\BB{c}{c''}$ and thus $\trans{c}{c'}\circ\trans{c'}{c''}\xrightarrow{\sim}\trans{c}{c''}$. The other isomorphism is similar.
\end{proof}
By the associativity of the multiplication $\star$ and the above lemma, $\EQ{\relbar}{\relbar}$ defines an equivalence relation.
\begin{lemm}\label{lemm:crit'}
	The condition $\EQ{c}{c'}$ holds if and only if for each $\lambda\in \haffd$, the following maps induced from multiplication $\star$ are isomorphisms: 
	\[
		\Bscr{c}{c'}{\lambda}\otimes_{\scrH_{c'}^{\lambda}}\Bscr{c'}{c}{\lambda}\xrightarrow{\sim}\scrH_{c}^{\lambda},\quad \Bscr{c'}{c}{\lambda}\otimes_{\scrH_{c}^{\lambda}}\Bscr{c}{c'}{\lambda}\xrightarrow{\sim}\scrH_{c'}^{\lambda}.
	\]
\end{lemm}
\begin{proof}
	This follows immediately from the formula~\autorefitem{prop:HCTor}{ii} and the faithful flatness of $\scrS$ over $\bfS$ due to~\autoref{lemm:platitude}.
\end{proof}

\subsection{Equivalence in a facet}
Given $(c,\lambda)\in \frakP\times \haffd$, similarly to \autoref{sssec:c-clans}, a connected component of the complement of $\bigcup_{a\in \Psi_{c,\lambda}}a^{-1}(0)$ in $\fraka$ is called a \emph{$(c,\lambda)$-clan} of $\fraka$. The set of $(c, \lambda)$-clans of $\fraka$ is denoted by $\Cl^{c,\lambda}(\fraka)$. \par

\begin{lemm} \label{lemm:clan-isom}
	Let $(c,\lambda)\in \frakP\times\haffd$ and suppose that $C, C'\in \Ch^{c,\lambda}(\fraka)$ are $(c, \lambda)$-chambers lying in the same $(c, \lambda)$-clan of $\fraka$. Then, $C$ and $C'$ are isomorphic as objects in $\scrA^{c,\lambda}$. 
\end{lemm}
\begin{proof}
Since $C$ and $C'$ lie in the same $(c, \lambda)$-clan, they are separated only by walls from $\Phi_{\lambda}$. The wall $\alpha^{-1}(0)\subseteq \fraka$ attached to an affine root $\alpha\in \Phi_{\lambda}$ separates pairs of chambers conjugate under the reflection $s_{\alpha}\in \Waff_{\lambda}$. It follows that $C$ and $C'$ are $\Waff_{\lambda}$-conjugate and thus isomorphic in $\scrA^{c,\lambda}$ (see \autoref{subsubsec:Acompl}). 
\end{proof}
Let $(c_{\bbR}, \lambda_{\bbR})\in \fraka = \frakP_{\bbR} \times \haffd_{\bbR}$ denote the real part of $(c, \lambda)$. Then, $\psi(c_{\bbR}, \lambda_{\bbR}) = 0$ for every $\psi\in \Psi_{c,\lambda}$. The translation by $(c_{\bbR}, \lambda_{\bbR})$ induces an affine isomorphism $\fraka^0\cong \fraka$ identifying the hyperplane arrangements defined by $\Psi_{c, \lambda}$ on $\fraka$ and $\ba\Psi_{c, \lambda}$ on $\fraka^0$. Therefore, the translation by $(c_{\bbR}, \lambda_{\bbR})$ induces a bijection $\Cl^{c,\lambda}(\fraka^0) \cong \Cl^{c,\lambda}(\fraka)$.  \par

\begin{prop}\label{prop:intra}
	Suppose that $c, c'\in \frakP$ satisfy $d = c - c'\in \frakP_{\bbZ}$ and lie in the same $c$-facet of $\frakP$. Then, the condition $\EQ{c}{c'}$ holds. 
\end{prop}
\begin{rema}
	The converse may not be true in general. 
\end{rema}

\begin{proof}[Proof of~\autoref{prop:intra}]
	By~\autoref{lemm:crit'}, it suffices to verify the condition therein. Let $\lambda\in \haffd$. Consider the category $\scrA^{c, \lambda}$ introduced in~\autoref{subsec:Aloc}. 
	Recall the isomorphisms~\eqref{equa:Bscr} and~\eqref{equa:Hscr} from~\autoref{subsec:Bcompl}.
	Let $F\subseteq \frakP_{\bbR}$ be the $c$-facet whose complexification $F_{\bbC}\subseteq \frakP$ (defined in \autoref{subsubsec:strat}) contain $c$ and $c'$. \par
	Let $c_{\bbR}, c'_{\bbR}\in \frakP_{\bbR}$ and $\lambda_{\bbR}\in \frakP_{\bbR}$ be the real parts of $c, c'$ and $\lambda$. We have $c'_{\bbR} = c_{\bbR} - d$. By hypothesis, $-c_{\bbR}$ and $-c'_{\bbR}$ lie in the $c$-chamber $-F$. For each $w\in \Waff$, let $\frakC^0_w\in \Cl^{c, \lambda}(\fraka^0)$ denote the $(c,\lambda)$-clan which contains the subset $w^{-1}\kappa_d - (c_{\bbR}, \lambda_{\bbR})\subset \fraka^0$. The image of $\frakC^0_w$ under the projection $\fraka^0\to \frakP_{\bbR}$ contains $d - c_{\bbR} = -c'_{\bbR}$. By~\autoref{lemm:strate}, this image contains $-F$, to which $-c_{\bbR}$ belongs. By translation, the image of the clan $\frakC_w := \frakC^0_w + (c_{\bbR},\lambda_{\bbR})\in \Cl^{c, \lambda}(\fraka)$ under the projection $\fraka \to \frakP_{\bbR}$ contains $0$. Therefore, there exists $w'\in \Waff$ such that $w'^{-1} \kappa_0$ lies in $\frakC_w$. Since $w^{-1}\kappa_d$ also lies in $\frakC_w$, \autoref{lemm:clan-isom} implies that $\widetilde{w^{-1} \kappa_d}$ and $\widetilde{w'^{-1} \kappa_0}$ are isomorphic as objects in $\scrA^{c,\lambda}$. It follows that the direct factor
	\[
		\bigoplus_{[y]\in \Waff/\Waff_{\lambda}}\Hom_{\scrA^{c, \lambda}}(\widetilde{w^{-1} \kappa_d}, \widetilde{y^{-1}\kappa_0}) / (\bfc - c) \subseteq  \Bscr{c}{c'}{\lambda}.
	\]
	is also a direct factor of $\scrH^{\lambda}_c$, which is projective. This implies that $\Bscr{c}{c'}{\lambda}$ is a projective $\scrH^{\lambda}_c$-module. Similar arguments show that $\Bscr{c'}{c}{\lambda}$ is a projective right $\scrH^{\lambda}_c$-module. Consider the following map:
	\[
		\Bscr{c'}{c}{\lambda}\otimes_{\scrH^{\lambda}_c} \Bscr{c}{c'}{\lambda}\to \scrH^{\lambda}_{c'},\quad a\otimes b\mapsto a\star b.
	\]
	We have seen that both sides are projective $\scrH^{\lambda}_{c'}$-modules. \autoref{prop:PBP} implies that this map becomes an isomorphism when $\scrP^{\lambda}_{c'}\otimes_{\scrH^{\lambda}_c}\relbar$ is applied to both sides. The double centraliser property \autoref{theo:KZ} then implies that the map in question is indeed an isomorphism. Exchanging $c$ and $c'$, we obtain the other isomorphism; hence, the condition of~\autoref{lemm:crit'} holds.
\end{proof}

\subsection{Antipodal wall-crossing}
Given $c\in \frakP$, two $c$-facets $F,F'\subseteq \frakP_{\bbR}$ are called \emph{antipodal} if $F' =-F$. In this case, the complexified $c$-facets $F_{\bbC}$ and $F'_{\bbC}$ (in the notation of \autoref{subsubsec:strat}) are also called \emph{antipodal}. On the other hand, given $(c,\lambda)\in \frakP\times \haffd$, we say that two $(c,\lambda)$-clans $\frakC,\frakC'\in\Cl^{c,\lambda}(\fraka)$ are antipodal if $\frakC' - (c_{\bbR}, \lambda_{\bbR}) = - (\frakC - (c_{\bbR}, \lambda_{\bbR}))$. 
\begin{lemm}\label{lemm:Tfixte-op}
	Let $\lambda\in \haffd$ and $c,c'\in \frakP$ such that $d:= c - c'\in \frakP_{\bbZ}$ and $c,c'$ belong to antipodal $c$-facets of $\frakP$. Then, for each $w\in \Waff$, there exists $y\in \Waff$ such that the chambers $\widetilde{w^{-1}\kappa_0}, \widetilde{y^{-1}\kappa_d}\in \Ch^{c, \lambda}(\fraka)$ lie in antipodal $(c,\lambda)$-clans.
\end{lemm}
\begin{proof}
	Let $c_{\bbR}, c'_{\bbR}$ and $\lambda_{\bbR}$ be the real parts of $c, c'$ and $\lambda$. By hypothesis, $-c_{\bbR}$ and $-c'_{\bbR}$ lie in antipodal $c$-facets: $-c_{\bbR}\in F$ and $-c'_{\bbR}\in -F$ for a $c$-facet $F\subseteq \frakP_{\bbR}$. 
	For each $w\in \Waff$, the image of the clan $\frakC^0_{w'}\in \Cl^{c,\lambda}(\fraka^0)$ containing $w'^{-1} \kappa_d - (c_{\bbR}, \lambda_{\bbR})$ under the projection $\fraka^0\to \frakP_{\bbR}$ contains $-c'_{\bbR}$. By~\autoref{lemm:strate}, this image should contain $F$. Consequently, the image of the antipodal clan $-\frakC^0_{w'}\in \Ch^{c, \lambda}(\fraka^0)$ contains the antipodal facet $-F$. By translation, the image of the clan $-\frakC^0_{w'} + (c_{\bbR}, \lambda_{\bbR}) \in \Cl^{c, \lambda}(\fraka)$ under the projection $\fraka \to \frakP_{\bbR}$ contains $0$. Therefore, there exists $y\in \Waff$ such that $\widetilde{y^{-1} \kappa_0}$ lies in $-\frakC^0_{w'} + (c_{\bbR}, \lambda_{\bbR})$, which is antipodal to $\frakC^0_{w'} + (c_{\bbR}, \lambda_{\bbR})$. 
\end{proof}
Recall the centre $\scrZ_c^{\lambda}$ from \autoref{subsubsec:Z}. We say that a compact module $M\in \scrH^{\lambda}_c\mof$ is \emph{relatively injective} if 
\begin{enumerate}
	\item the $\scrZ_c^{\lambda}$-module structure on $M$ introduced in~\autoref{subsubsec:Z} is free, and
	\item $\Ext^{>0}_{\scrH^{\lambda}_c}(M, N) = 0$ for every $N\in \scrH^{\lambda}_c\mof$ satisfying the condition (i).
\end{enumerate}
The functor $M\mapsto \hom_{\scrZ_c^{\lambda}}(M, \scrZ_c^{\lambda}) =\bigoplus_{\lambda\in [\lambda]}\Hom_{\scrZ_c^{\lambda}}(M 1_{c, \lambda}, \scrZ_c^{\lambda})$ induces a bijection between the isomorphism classes of compact projective right $\scrH^{\lambda}_c$-modules and compact relatively injective left $\scrH^{\lambda}_c$-modules.
\begin{prop}\label{prop:antipode}
	Let $c,c'\in \frakP$ such that $d:= c - c'\in \frakP_{\bbZ}$ and $c,c'$ belong to antipodal $c$-chambers of $\frakP_{\bbQ, c}$. Then, the functor $\trans{c'}{c}$	is an equivalence of categories.
\end{prop}
\begin{proof}
	By~\autoref{lemm:crit}, it suffices to show that $\trans{c'}{c}^{\lambda}$ is an equivalence for $\lambda\in \haffd$.

	Let $w\in \Waff$. By~\autoref{lemm:Tfixte-op}, there exists $y\in \Waff$ such that $\widetilde{w^{-1} \kappa_0}$ and  $\widetilde{y^{-1} \kappa_d}$ lie in antipodal clans. Set
	\begin{align*}
		&\scrI_{w} = \bigoplus_{w'\in \Waff / \Waff_{\lambda}}\Hom_{\scrA^{c, \lambda}}(\widetilde{w'^{-1} \kappa_d}, \widetilde{w^{-1} \kappa_0}) / (\bfc - c) \\
		&\scrP_{w} = \bigoplus_{w'\in \Waff / \Waff_{\lambda}}\Hom_{\scrA^{c, \lambda}}(\widetilde{y^{-1} \kappa_d}, \widetilde{w'^{-1} \kappa_d}) / (\bfc - c).
	\end{align*}
	Set $C' = \widetilde{w^{-1} \kappa_0}$ and let $C\in \Ch^{c,\lambda}(\fraka)$ be the $(c,\lambda)$-chamber antipodal to $C'$. Then, $\widetilde{y^{-1} \kappa_d}$ and $C$ lie in the same $(c,\lambda)$-clan, so that we have $C\cong \widetilde{y^{-1}\kappa_d}$ in $\scrA^{c,\lambda}$ by \autoref{lemm:clan-isom}. 
	Then, the composition in the category $\scrA^{c,\lambda}$ yields an $\scrH^\lambda_{c'}$-bilinear map
	\[
		\scrI_{w} \times \scrP_w\to \Hom_{\scrA^{c,\lambda}}(\widetilde{y^{-1} \kappa_d}, C') / (\bfc - c)\cong\Hom_{\scrA^{c,\lambda}}(C, C') / (\bfc - c).
	\]
	Since $C$ and $C'$ are antipodal, \autoref{prop:duality} yields a trace map 
	\[
		\Hom_{\scrA^{c,\lambda}}(C, C')\to \scrZ^{c,\lambda}
	\]
	which induces a perfect pairing for each $\lambda'\in [\lambda]$:
	\[
		1_{c',\lambda'}\scrI_{w} \times \scrP_w1_{c',\lambda'}\to \scrZ^{c,\lambda} / (\bfc - c) = \scrZ_c^{\lambda}
	\]
 and thus an isomorphism $1_{c',\lambda'}\scrI_w \cong \Hom_{\scrZ_c^{\lambda}}(\scrP_w1_{c',\lambda'}, \scrZ_c^{\lambda})$. Summing over $\lambda'\in [\lambda]$, we obtain an isomorphism of $\scrH_{c'}^{\lambda}$-modules $\scrI_w \cong \hom_{\scrZ_c^{\lambda}}(\scrP_w, \scrZ_c^{\lambda})$. The right $\scrH^{\lambda}_{c'}$-module $\scrP_w$ being a direct factor of $\scrH^{\lambda}_{c'}$, is projective. It follows that $\scrI_w$ is a relatively injective left $\scrH^{\lambda}_{c'}$-module. Similar arguments show that the sum $\Bscr{c}{c'}{\lambda} = \bigoplus_{w\in \Waff / \Waff_{\lambda}} \scrI_w$ contains every indecomposable relatively injective left $\scrH_{c'}^{\lambda}$-module as direct summand. Since the direct sum of indecomposable relatively injective modules is a tilting module, the condition \ref{lemm:crit-ii'} of \autoref{lemm:crit} is satisfied.
\end{proof}
\subsection{Simple wall-crossing}
Let $c, c'\in \frakP$ such that $d = c - c'$ lies in $\frakP_{\bbZ}$ and $c, c'$ belong to different $c$-chambers $\frakP$ whose real parts are separated by a single wall, say $[\mu]\in \frakM_c$. 
\begin{prop}\label{prop:simple}
	Under the above hypothesis, the functor $\trans{c'}{c}$ is an equivalence of categories.
\end{prop}
\begin{proof}
	We prove this proposition by the degeneration technique from~\cite{BL21}. Let $D\subset \frakP$ be the complex hyperplane parallel to $h = 0$ such that $c\in D$ and put $\bfB_D = \BB{}{d}\otimes_{\calO(\frakP)}\calO(D)$. Let $U\subseteq D$ be the subset of $\bbC$-points of $D$ defined by
	\[
		U = \left\{ u\in D\;;\; \mu'(u) \notin\bbQ,\;\forall [\mu']\in \frakM_c \setminus \left\{ [\mu] \right\}\right\}.
	\]
	It follows that $\frakM_u = \left\{ [\mu] \right\}$ consists of a single wall, so that the pair $(u, u')$ belongs to antipodal $c$-chambers. Applying~\autoref{prop:antipode} to $u\in U$ and setting $u' = u-d$, we see that $\BB{u'}{u}\otimes^{\rmL}_{\bfH_u}\relbar$ induces an equivalence of categories. Consequently, the following adjoint unit and counit are quasi-isomorphisms:
	\[
		\varepsilon_u:\bfH_u \to \RHom_{\bfH_{u'}}(\BB{u'}{u}, \BB{u'}{u}),\quad \eta_u:\BB{u'}{u}\otimes^{\rmL}_{\bfH_u}\RHom_{\bfH_{u'}}(\BB{u'}{u}, \bfH_{u'})\to \bfH_{u'}.
	\]
	Set $D' = D - d\subseteq \frakP$, $\bfH_{D} = \bfH\otimes_{\calO(\frakP)}\calO(D)$, $\bfH_{D'} = \bfH\otimes_{\calO(\frakP)}\calO(D')$ and
	\begin{align*}
		K &= \Cone(\bfH_D \to \RHom_{\bfH_{D'}}(\bfB_D, \bfB_D)), \\
		K' &= \Cone(\bfB_D\otimes^{\rmL}_{\bfH_D}\RHom_{\bfH_{D'}}(\bfB_D, \bfH_{D'})\to \bfH_{D'}).
	\end{align*}
	It follows that $K\otimes^{\rmL}_{\calO(D)}\kappa(u) \cong \Cone(\varepsilon_u) = 0$ and $K'\otimes^{\rmL}_{\calO(D)}\kappa(u) \cong \Cone(\varepsilon_u) = 0$ hold for each $u\in U$ and hence $\rmH^\bullet(K)_u = 0 = \rmH^\bullet(K')_u$. \par
	Applying the lemma of generic freeness \autoref{lemm:gen-libre} below to $\rmH^\bullet(K)$ viewed as an $(\bfH, \calO(D))$-bimodule, there exists a Zariski-open dense subset $V\subset D$ containing $U$ such that $K\mid_V = 0$. Similarly, we may assume that $K'\mid_V = 0$. Since the subset
	\[
		D'=\left\{z\in D\cap (c+ \frakP_{\bbZ})\;;\; (c, z) \text{ same $c$-facet},\; (c', z-d) \text{ same $c$-facet} \right\}
	\]
	is Zariski-dense in $D$, we may choose $z\in D'\cap V\neq \emptyset$, so that $K_z = 0$ and $K'_z = 0$. Set $z' = z-d$, the translation functor $\trans{z'}{z}$ is an equivalence of categories. By~\autoref{prop:intra}, the conditions $\EQ{c}{z}$ and $\EQ{c'}{z'}$ hold, and it follows from~\autoref{lemm:t-exact} that $\trans{c'}{c} \cong \trans{c'}{z + d}\circ\trans{z'}{z}\circ\trans{z}{c}$ and it is also an equivalence of categories.
\end{proof}

\begin{lemm}\label{lemm:gen-libre}
	Let $E$ be a noetherian integral commutative $\bbC$-algebra and let $M$ be a finitely generated $(\bfH, E)$-bimodule. Then, there exists a non-zero element $f\in E$ such that the localisation $M_f = M\otimes_{E}E_f$ is free as $E_f$-module, where $E_f = E[f^{-1}]$.
\end{lemm}
\begin{proof}
	The proof is modelled on~\cite[5.7, 5.8]{BL21}. We extend the filtration $F^{\can}_{\bullet} \bfH$ to $\bfH\otimes E$ by setting $F^{\can}_n(\bfH\otimes E) = (F^{\can}_n\bfH)\otimes E$. Pick a good filtration $\left\{ F_k M \right\}_{k\in \bbZ}$ for $M$ with respect to $F^{\can}_{\bullet} (\bfH\otimes E)$. Then, the associated graded $\gr^F M$ is finitely generated over $(\gr^{\can}\bfH)\otimes E$ and is also finitely generated over the centre $Z(\gr^{\can}\bfH)\otimes E$. The lemma of generic freeness~\cite[6.9.2]{EGAIV2} implies that there exists a non-zero element $f\in E$ such that $(\gr^{F}M)_f$ is free over $E_f$. Since the filtration $\left\{ F_k M \right\}_{k\in \bbZ}$ is bounded below, $M_f$ is isomorphic to $(\gr^{F}M)_f$ as free module over $E_f$. 
\end{proof}

\subsection{General wall-crossing}
\begin{prop}\label{prop:general}
	Let $c, c'\in \frakP$ be such that $d = c - c'\in \frakP_{\bbZ}$. Suppose that $c$ and $c'$ lie in $c$-chambers of $\frakP$. Then, the functor $\trans{c'}{c}$	is an equivalence of categories.
\end{prop}
\begin{proof}
	Let $F,F'\subseteq \frakP_{\bbR}$ be the $c$-chambers such that $c\in F_{\bbC}$ and $c'\in F'_{\bbC}$. We prove the statement by induction on the number $k\in \bbN$ of walls in $\frakM_c$ separating $F$ and $F'$. The case $k = 0$ is covered by~\autoref{prop:intra}. Let $k > 0$ and assume that the statement holds for smaller $k$. Let $F''$ be an $c$-chamber which is separated from $F$ by one wall and from $F'$ by $k - 1$ walls. Let $[\mu_0]\in \frakM_c$ be the unique element such that $\mu_0^{-1}(0)$ separates $F$ and $F''$. Choose $c''\in F''_{\bbC} \cap (c + \frakP_{\bbZ})$ and put $d'' = c - c''$. Let $D = \left\{ \mu_0 - \mu_0(c) = 0 \right\}\subset \frakP$ be the complex hyperplane parallel to $\mu_0 = 0$ such that $c\in D$. Choose an element $d'\in \frakP_{\bbZ} \cap F'$ and put $\bfB^{\otimes} = \bigoplus_{n\ge 0}\BB{}{nd'}$ and $\bfC^{\otimes} = \bigoplus_{n\ge 0}\BB{}{nd'+d''}$. A variant of the proof of~\autoref{lemm:fg} shows that $\bfC^{\otimes}$ is a finitely generated left $\bfB^{\otimes}$-module. The multiplication $\star$ yields a morphism
	\[
		m:\bfB^{\otimes}\otimes^{\rmL}_{\bfH}\BB{}{d''}\to \bfC^{\otimes},
	\]
	which induces
	\[
		m_D:\bfB^{\otimes}\otimes^{\rmL}_{\bfH}\BB{}{d''}\otimes^{\rmL}_{\calO(\frakP)}\calO(D)\to \bfC^{\otimes}\otimes^{\rmL}_{\calO(\frakP)}\calO(D).
	\]
	Put $K = \Cone(m_D)\in \Db(\bfB^{\otimes}\bimod \calO(D))$. Let $U\subset D$ be the subset of $\bbC$-points of $D$ defined by
	\[
		U = \left\{ u\in D\;;\; \mu(u) \notin\bbQ,\;\forall [\mu]\in \frakM_c \setminus \left\{ [\mu_0] \right\}\right\}.
	\]
	It is a dense subset of $D$. We have $K\otimes^{\rmL}_{\calO(D)}\kappa(u) = 0$ for each $u\in U$ by~\autorefitem{prop:intra}{ii}, because $\frakM_u = \left\{ [\mu] \right\}$ and $u-d''$ and $u-d'' - nd'$ belong to the same $c$-chamber for all $n\in \bbN$. Since $U$ is dense in $D$, the generic freeness~\autoref{lemm:finitude} below implies that there exists $f\in \calO(D) \setminus\left\{ 0 \right\}$ such that $K\otimes^{\rmL}_{\calO(D)}\calO(D)[f^{-1}] = 0$. In other words, there exists a Zariski-open dense subset $V\subset D$ containing $U$ such that $K\mid_V = 0$. \par
	On the other hand, the subset 
	\[
		D' = \left\{ z\in D\cap(c + \frakP_{\bbZ})\;;\;z\in F_{\bbC}, z - d''\in F''_{\bbC} \right\}
	\]
	is Zariski-dense in $D$ because $F$ and $F''$ are separated by a single wall $\{\mu_0 = 0\}$. The density implies that there exists $z\in D'\cap V\neq \emptyset$, so that $z\in F_{\bbC}$, $z'' := z - d''\in F''_{\bbC}$ and $K \otimes^{\rmL}_{\calO(D)}\kappa(z) = 0$. Pick $n_0 \gg 0$ such that $z':= z - d''- n_0 d' \in F'_{\bbC}$. The condition $K \otimes^{\rmL}_{\calO(D)}\kappa(z) = 0$ implies that $\BB{z'}{z''}\otimes^{\rmL}_{\bfH_{z''}}\BB{z''}{z} \xrightarrow{\sim}\BB{z'}{z}$ and hence $\trans{z'}{z''}\circ\trans{z''}{z} \xrightarrow{\sim} \trans{z'}{z}$. 
	Since $F''$ and $F'$ are separated by $k - 1$ walls in $\frakM_c$, the induction hypothesis implies that $\trans{z'}{z''}$ is an equivalence of categories. On the other hand, \autoref{prop:simple} implies that $\trans{z''}{z}$ is also an equivalence of categories; hence, so is $\trans{z'}{z}$ an equivalence of categories. \autoref{prop:intra} implies that the conditions $\EQ{c}{z}$ and $\EQ{c'}{z'}$ hold, so by \autoref{lemm:t-exact} there is an isomorphism:
	\[
		\trans{c'}{z'}\circ\trans{z'}{z}\circ\trans{z}{c} \xrightarrow{\sim} \trans{c'}{c}.
	\]
	 Therefore, $\trans{c'}{c}$ is also an equivalence of categories.
\end{proof}
 
\begin{lemm}\label{lemm:finitude}
	Let $E$ be a noetherian integral commutative $\bbC$-algebra. For every finitely generated $(\bfB^{\otimes}, E)$-bimodule $M$, there exists $f\in E \setminus \{0\}$ such that $M_f = M\otimes_E E_f$ is free over $E_f$.
\end{lemm}
\begin{proof}
	The arguments are similar to the proof of~\autoref{lemm:gen-libre}:  find a good filtration $\{F_{j,k}M\}_{j,k\in \bbZ}$ with respect to the $\bbN^2$-filtration $F_{\bullet,\bullet}\bfB^{\otimes}\otimes E$; pass to the associated quotient $\gr^{F}M$; the restriction of $\gr^{F}M$ is finitely generated over $A\otimes E$ due to~\autoref{lemm:fg}, where $A$ is a commutative $\bbC$-algebra of finite type; apply the generic freeness to deduce that when $E$ is integral, there exists a non-zero element $f\in E$ such that $\gr^F M_f$ is a free $E_f$-module. It follows that $M_f$ is also free over $E_f$.
\end{proof}

\subsection{Exceptional case}

In this subsection, we deal with the case where $c$ lies in a $c$-facet of $\frakP$ which may not be a chamber.
\subsubsection{}
Consider the most exceptional case, where the parameter $c\in \frakP$ lies in the intersection of all the hyperplanes appearing in $\frakM_c$. The following proposition shows that in this case, all the hyperplanes in $\frakM_c$ become \emph{irrelevant} and the translation functors are abelian equivalences.
\begin{prop}\label{prop:except}
	Suppose that $c\in \frakP$ satisfies $\mu(c) = 0$ for every $[\mu]\in \frakM_c$. Then, $\EQ{c}{c'}$ holds for each $c'\in c + \frakP_{\bbZ}$.
\end{prop}
This will result from the following lemma:
\begin{lemm}\label{lemm:KZequiv}
	Let $c$ and $c'$ be as in~\autoref{prop:except}. Then, the algebraic KZ functor $\bfV: \scrH_{c'}^{\lambda}\Mod\to \scrK_t^{\ell}\Mod$ is an equivalence of categories for each $\lambda = (\lambda^0, 1)\in \haffd$ and $\ell = \exp(2\pi i\lambda^0)\in T$.
\end{lemm}
\begin{proof}
	Consider first the case where $c' = c$. We show that:
	\begin{quote}
		for each $\lambda\in \haffd$, every $(c, \lambda)$-clan is generic (see \autoref{defi:generic}). 
	\end{quote}
	We prove this assertion by contradiction. Suppose there exists a non-generic $(c,\lambda)$-clan $\frakC\in \Cl^{c,\lambda}(\haffd)$. Then, the set 
	\[
		\Sigma = \{\alpha\in \Psi_{c, \lambda}\cup -\Psi_{c, \lambda}\;;\; \alpha(\frakC) \subseteq \bbR_{>0}\}
	\]
	contains a positive $\bbQ$-linear relation $\sum_{\alpha\in \Sigma'} d_{\alpha}\ba\alpha = 0$, where $\Sigma'\subseteq \Sigma$ is a non-empty subset, $\ba\alpha = \alpha\mid_{\frakh_{\bbQ}\times \frakP_{\bbQ}}$  and $d_\alpha\in \bbQ_{>0}$ for each $\alpha\in \Sigma'$. We assume that $\Sigma'$ is minimal with this property and let $\sigma = \left\{ \ba \alpha\in \frakh^*_{\bbQ}\times \frakP_{\bbQ}^*\;;\; \alpha\in \Sigma' \right\}$. Let $h\in \frakC$, so that $\alpha(h)>0$ for each $\alpha\in \Sigma'$. It follows that
	\[
		0< \sum_{\alpha\in \Sigma'} d_{\alpha}\alpha(h) = \sum_{\alpha\in \Sigma'} (d_{\alpha} c_{\alpha} + \langle d_{\alpha}\ba\alpha, h - \lambda\rangle) = \sum_{\alpha\in \Sigma'} d_{\alpha}c_{\alpha}.
	\]
	By the minimality of $\Sigma'$, the set $\sigma$ is a circuit, so that $[\mu_{\sigma}]\in \frakM_c$ holds for $\mu_{\sigma} = \sum_{\alpha\in \sigma}d_\alpha \bfc_{\alpha}$. This contradicts the assumption that $\mu_{\sigma}(c) = 0$. The claim is proven.\par
	Since $\scrH^{\lambda}_c 1_{\bfV}$ is a direct sum of such $\scrH^{\lambda}_c 1_{w \lambda}$ for $w\in \Waff$ such that $w^{-1} \nu_0$ lies in a generic $(c,\lambda)$-clan, the previous paragraph shows that $\scrH^{\lambda}_c1_{\bfV}$ is a compact projective generator of $\scrH^{\lambda}_c\Mod$, so the functor $\bfV = \Hom_{\scrH^{\lambda}_c}(\scrH^{\lambda}_c1_{\bfV}, \relbar)$ is an equivalence of categories. Consequently, the category $\scrK^{\ell}_h\Mod$ has finite global dimension. Since $\scrK^{\ell}_h$ is a Frobenius algebra over the centre $\scrZ^{\ell}_h = \rmZ(\scrK^{\ell}_h)\cong \scrZ^{\lambda}_c$ (see~\cite[9.3]{liu22}), all $\scrK^{\ell}_h$-modules which are projective over $\scrZ^{\ell}_h$ must be projective.\par
	Now, let $c'\in c + \frakP_{\bbZ}$. Since $\bfV$ sends $\scrZ^{\lambda}_c$-projective objects to $\scrZ^{\ell}_h$-projective objects, it preserves projective objects. The double centraliser property (see~\autoref{theo:KZ}) implies that $\bfV$ is an equivalence of categories.
\end{proof}
\begin{proof}[Proof of~\autoref{prop:except}]
	By~\autoref{lemm:crit'}, it suffices to verify the conditions therein. Let $\lambda\in \haffd$. \autoref{lemm:KZequiv} implies that it suffices to show that 
	\[
		\scrP^{\lambda}_c\otimes_{\scrH^{\lambda}_c}\Bscr{c}{c'}{\lambda}\otimes_{\scrH^{\lambda}_{c'}}\Bscr{c'}{c}{\lambda}\xrightarrow{\scrP^{\lambda}_c\otimes\star} \scrP^{\lambda}_c\quad \text{and } 
		\scrP^{\lambda}_{c'}\otimes_{\scrH^{\lambda}_{c'}}\Bscr{c'}{c}{\lambda}\otimes_{\scrH^{\lambda}_{c}}\Bscr{c}{c'}{\lambda}\xrightarrow{\scrP^{\lambda}_{c'}\otimes\star} \scrP^{\lambda}_{c'}
	\]
	are isomorphisms. By~\autoref{prop:PBP}, these maps can be factorised as $\alpha_{c', c}\circ(\alpha_{c, c'}\otimes \Bscr{c'}{c}{\lambda})$ and $\alpha_{c, c'}\circ(\alpha_{c', c}\otimes \Bscr{c}{c'}{\lambda})$, which are isomorphisms.
\end{proof}
\subsubsection{}
Let $c\in \frakP$ and let $F\subset \frakP_{\bbR}$ be the $c$-facet whose complexification $F_{\bbC}$ contains $c$. Let $F'\subset \frakP_{\bbR}$ be a $c$-chamber whose closure contains $F$ and let $c'\in F'_{\bbC}\cap (c+\frakP_{\bbZ})$. Set $d = c - c'$.
\begin{prop}\label{prop:except1}
	Under the above hypothesis, the condition $\EQ{c}{c'}$ holds.
\end{prop}
\begin{proof}
	Let $E\subseteq \frakP$ be the $\bbC$-linear span of $F$ and set $D = c' + E$, let $\frakM'_c = \left\{[\mu]\in \frakM_c\;;\; \mu(c) \neq 0  \right\}$ and let $U\subset D$ be the subset defined by
	\[
		U = \left\{ u\in D\;;\; \mu(u) \notin\bbQ,\;\forall [\mu]\in \frakM'_c\right\}.
	\]
	Put $\bfB_D = \BB{}{d}\otimes_{\calO(\frakP)}\calO(D)$. For each $u\in U$, we have $\frakM_u = \frakM_c \setminus \frakM'_{c}$, so the condition $\EQ{u}{u + d}$ holds by~\autoref{prop:except}. Again, the degeneration argument used in the proof of~\autoref{prop:simple} implies that there exists $z\in D\cap F_{\bbC}$ such that $z':= z - d \in F'_{\bbC}$ and $\EQ{z}{z'}$ holds. \autoref{prop:intra} implies that $\EQ{c}{z}$ and $\EQ{c'}{z'}$ hold. Hence, $\EQ{c}{c'}$ holds.
\end{proof}
\subsection{Proof of \autoref{theo:main}}
Choose $c$-chambers $F, F'\subset \frakP_{\bbR}$ such that $c\in \ba{F_{\bbC}}$ and $c'\in \ba{F'_{\bbC}}$. Let $z\in F_{\bbC}\cap (c+ \frakP_{\bbZ})$ and $z'\in F'_{\bbC}\cap (c+ \frakP_{\bbZ})$. Then, $\trans{z'}{z}$ is an equivalence of categories by~\autoref{prop:general}; $\trans{c'}{z'}$ and $\trans{c'}{z'}$ are equivalences of categories and $\trans{c'}{c} \cong \trans{c'}{z'}\circ \trans{z'}{z}\circ \trans{z}{c}$ holds by~\autoref{prop:except1} and~\autorefitem{lemm:t-exact}{ii}; hence, $\trans{c'}{c}$ is an equivalence of categories. \qed

\subsection{Derived equivalences on integrable modules}\label{subsec:equivO}
We prove that the derived equivalence established in \autoref{theo:main} preserves the integrable modules, their blocks and their support.
\subsubsection{}
Consider the subring $\calO(T^{\vee})\cong \bbC Q^{\vee} \subseteq \bfH$, where $T^{\vee}$ is the torus whose character lattice is the coroot lattice $Q^{\vee}$ of $(\hfin_{\bbR}, R)$. Given any $M\in \rmO(\bfH)$, by the triangular decomposition $\bfH = \bbC Q^{\vee}\otimes \bbC W\otimes \calO(\frakP\times \haffd)$, it is clear that $M$ is finitely generated over $\bbC Q^{\vee}$. Let $\supp_{T^{\vee}/W}M\subseteq T^{\vee} / W$ be the (set-theoretic) support of $M$ considered as $\calO(T^{\vee})^W$-module by restriction. 
\begin{lemm}\label{lemm:supp}
	Let $B\in \HC(\bfH)$ and $M\in \rmO(\bfH)$. Then, $\supp_{T^{\vee}/W}\Tor^{\bfH}_i(B, M)\subseteq \supp_{T^{\vee}/W} M$ and $\supp_{T^{\vee}/W}\Ext_{\bfH}^i(B, M)\subseteq \supp_{T^{\vee}/W} M$ hold for each $i\in \bbZ$. 
\end{lemm}
\begin{proof}
Choose a HC filtration $\left\{ F_{k} B \right\}_{k\in \bbZ}$ for $B$ and a good filtration $\left\{ F_{k} M \right\}_{k\in \bbZ}$ for $M$ with respect to the canonical filtration $F^{\can}_{\bullet}\bfH$. Note that the subspaces $F_k B$ and $F_k M$ are stable under multiplication by $\calO(T^{\vee})$ for each $k\in \bbZ$.  Since $\gr^F B$ is finitely generated over the centre $\rmZ(\bfH_{\delta} / (\delta))$ and the adjoint action of the latter is nilpotent on $\gr^F B$, there exists a number $m > 0$ such that $(\ad \rmZ(\bfH_{\delta} / (\delta)))^m\gr^F B = 0$. Note that $\calO(T^{\vee})^W\subseteq \rmZ(\bfH_{\delta} / (\delta)))$. Let $I = \ann_{\calO(T^{\vee})^W} M$ be the annihilator of $M$. Then, the ideal $I^m$ annihilates $\Ext^i_{\bfH}(\gr^F B, \gr^F M)$. \par
	As noted in the proof of~\autoref{prop:extHC}, there is an induced good filtration on $\Ext^i_{\bfH}(B, M)$, denoted by $F_{\bullet}$, such that $\gr^F\Ext^i_{\bfH}(B, M)$ is a subquotient of $\Ext^i_{\bfH}(\gr^F B, \gr^F M)$. Hence, $I^m$ annihilates $\gr^F\Ext^i_{\bfH}(B, M)$ and thus $\supp_{T^{\vee}/W} \Ext^i_{\bfH}(B, M)\subseteq V(\sqrt{I}) = \supp_{T^{\vee}/W}M$. The result for $\Ext$ follows. The statement for $\Tor$ can be proven in a similar way. 
\end{proof}

\subsubsection{}
For any closed subvariety $Z\subseteq T^{\vee}/W$ and $(c, [\lambda])\in \frakP\times \haffd/\Waff$, let $\rmO_{\lambda, Z}(\bfH_c)$ be the full subcategory of $\rmO_{\lambda}(\bfH_c)$ consisting of objects $M$ satisfying $\supp_{T^{\vee}/W}(M)\subseteq Z$ and let $\Db_{Z}(\rmO_\lambda(\bfH_c))$ (resp. $\Db_{\rmO_{\lambda, Z}}(\bfH_c)$) be the full subcategory of $\Db(\rmO_{\lambda}(\bfH_c))$ (resp. of $\Db(\bfH_c)$) consisting of complexes $K^{\bullet}$ such that $\rmH^k(K^{\bullet})\in \rmO_{\lambda, Z}(\bfH_c)$ for every $k\in \bbZ$. 
\begin{theo}\label{theo:equivO}
	Under the hypothesis of~\autoref{theo:main}, the translation functor $\trans{c'}{c}$ induces to an equivalence $\Db_{Z}(\rmO_{\lambda}(\bfH_{c}))\xrightarrow{\sim} \Db_{Z}(\rmO_{\lambda}(\bfH_{c'}))$ for each $[\lambda]\in \haffd / \Waff$ and each closed subset $Z\subseteq T^{\vee}/W$. 
\end{theo}
\begin{proof}
	Since there is natural isomorphism $\calO^{\wedge}_{[\lambda], \haffd}\otimes_{\calO(\haffd)}\BB{c'}{c}\cong \BB{c'}{c}\otimes_{\calO(\haffd)}\calO^{\wedge}_{[\lambda], \haffd}$ (see~\autoref{subsec:Bcompl}), by \autoref{theo:DequivO}, \autoref{prop:HCTor}, the restriction of $\trans{c'}{c}$ to $\Db_{\rmO_{\lambda}}(\bfH_c)$ is given by 
	\[
		\left(\scrC^{\lambda}_c\right)^{-1}\circ \trans{c'}{c}^{\lambda}\circ \scrC^{\lambda}_c: \Db_{\rmO_{\lambda}}(\bfH_{c}) \to \Db_{\rmO_{\lambda}}(\bfH_{c'}).
	\]
	By \autoref{lemm:supp}, it restricts to  a functor from from $\Db_{\rmO_{\lambda, Z}}(\bfH_{c})$ to $\Db_{\rmO_{\lambda, Z}}(\bfH_{c'})$. Similarly, the right adjoint functor of $\trans{c'}{c}$ restricts to a functor from $\Db_{\rmO_{\lambda, Z}}(\bfH_{c'})$ to $\Db_{\rmO_{\lambda, Z}}(\bfH_{c})$. By~\autoref{theo:main}, they are inverse to each other, thus equivalences of categories. Applying \autoref{coro:DO}, we obtain
	\[
		\Db_{Z}(\rmO_{\lambda}(\bfH_{c}))\cong \Db_{\rmO_{\lambda, Z}}(\bfH_{c})\xrightarrow{T_{c'\leftarrow c}} \Db_{\rmO_{\lambda, Z}}(\bfH_{c'})\cong \Db_{Z}(\rmO_{\lambda}(\bfH_{c'})).
	\]
\end{proof}
\begin{rema}
	With slight modifications, all results and arguments in the present article work over an arbitrary field of characteristic $0$ instead of $\bbC$. 
\end{rema}
\begin{rema}
	In~\cite{losev17}, it is shown that certain wall-crossing derived equivalences obtained from translation bimodules are perverse in the sense of Rouquier~\cite{rouquier06}. We expect the derived equivalences in~\autoref{theo:equivO} to be perverse as well.
\end{rema}

\printbibliography

\end{document}